\documentclass[11pt]{article}

\usepackage{xypic}
\usepackage{amsgen}
\usepackage{amsmath}
\usepackage{amstext}
\usepackage{amsbsy}
\usepackage{amsopn}
\usepackage{amsfonts}
\usepackage{amssymb}
\usepackage{eepic}
\usepackage[pdftex]{graphicx}
\usepackage{epsf}
\usepackage{pstricks}
\xyoption{all}

\def\Box{\square}
\def\edg{\relbar\joinrel\relbar}
\def\mapright#1{\smash{\mathop{\longrightarrow}\limits^{#1}}}
\def\tra#1{\smash{\mathop{\mid\kern
-1pt\joinrel\relbar\joinrel\relbar}\limits^{*}_{#1}}}
\def\longtra#1{\smash{\mathop{\mid\kern
-1pt\joinrel\relbar\joinrel\relbar\joinrel\relbar}\limits^{*}_{#1}}}
\def\vlongtra#1{\smash{\mathop{\mid\kern
-1pt\joinrel\relbar\joinrel\relbar\joinrel\relbar\joinrel\relbar}\limits^{*}_{#1}}}
\def\vvlongtra#1{\smash{\mathop{\mid\kern
-1pt\joinrel\relbar\joinrel\relbar\joinrel\relbar\joinrel\relbar\joinrel\relbar}\limits^{*}_{#1}}}
\def\vvvlongtra#1{\smash{\mathop{\mid\kern
-1pt\joinrel\relbar\joinrel\relbar\joinrel\relbar\joinrel\relbar\joinrel\relbar\joinrel\relbar}\limits^{*}_{#1}}}
\def\etra#1{\smash{\mathop{\mid\kern
-1pt\joinrel\relbar\joinrel\relbar}\limits_{#1}}}

\def\vlongrightarrow{\relbar\joinrel\longrightarrow}
\def\vvlongrightarrow{\relbar\joinrel\vlongrightarrow}
\def\vvvlongrightarrow{\relbar\joinrel\vvlongrightarrow}
\def\vvvvlongrightarrow{\relbar\joinrel\vvvlongrightarrow}
\def\vvvvvlongrightarrow{\relbar\joinrel\vvvvlongrightarrow}

\def\longmapright#1{\smash{\mathop{\vlongrightarrow}\limits^{#1}}}

\def\vvvvlongmapright#1{\smash{\mathop{\vvvvvlongrightarrow}\limits^{#1}}}

\def\id{\mbox{Id}}
\def\iff{\Leftrightarrow}
\def\Rw{\Rightarrow}
\def\oo{\overline}

\def\wt{\widetilde}
\def\wh{\widehat}
\def\B{{\cal{B}}}

\def\E{\mbox{Rh}}

\def\L{\mathrel{{\mathcal L}}}
\def\M{{\cal{M}}}
\newcommand{\N}{{\rm I}\kern-2pt {\rm N}}

\def\im{\mbox{im}}

\def\deg{\mbox{deg}}

\def\dom{\mbox{dom}}
\def\ell{\mbox{Ell}}

\def\lm{\mbox{lm}}
\def\ray{\mbox{Ray}}
\def\mray{\mbox{MRay}}
\def\geo{\mbox{Geo}}
\def\vert{\mbox{Vert}}
\def\edge{\mbox{Edge}}

\def\stab{\mbox{Stab}}

\def\con{\mbox{Con}}
\def\dep{\mbox{dep}}
\def\endo{\mbox{End}}

\def\max{\mbox{max}}
\def\min{\mbox{min}}
\def\son{\mbox{Sons}}
\def\sup{\mbox{sup}}

\def\P{{\cal{P}}}

\def\R{\mathrel{{\mathcal R}}}
\def\H{\mathrel{{\mathcal H}}}
\def\D{\mathrel{{\mathcal D}}}

\def\p{\varphi}
\def\pinv{{\p \inv}}
\def\inv{^{-1}}

\def\J{\mathrel{{\mathcal J}}} 
\def\bi{\begin{itemize}}
\def\ei{\end{itemize}}
\def\beq{\begin{equation}}
\def\eeq{\end{equation}}

\newtheorem{T}{Theorem}[section]
\newcommand{\bt}{\begin{T}}
\newcommand{\et}{\end{T}}
\newcommand{\ftd}{$\square$\end{T}}

\newtheorem{Proposition}[T]{Proposition}
\newcommand{\bp}{\begin{Proposition}}
\newcommand{\ep}{\end{Proposition}}
\newcommand{\fpd}{$\square$\end{Proposition}}

\newtheorem{Definition}[T]{Definition}
\newcommand{\bd}{\begin{Definition}}
\newcommand{\ed}{\end{Definition}}

\newtheorem{Lemma}[T]{Lemma}
\newcommand{\bl}{\begin{Lemma}}
\newcommand{\el}{\end{Lemma}}
\newcommand{\fld}{$\square$\end{Lemma}}

\newtheorem{Corol}[T]{Corollary}
\newcommand{\bc}{\begin{Corol}}
\newcommand{\ec}{\end{Corol}}
\newcommand{\fcd}{$\square$\end{Corol}}

\newtheorem{Result}[T]{Result}
\newcommand{\br}{\begin{Result}}
\newcommand{\er}{\end{Result}}
\newcommand{\frd}{$\square$\end{Result}}

\newtheorem{Example}[T]{Example}
\newcommand{\be}{\begin{Example}}
\newcommand{\ee}{\end{Example}}

\newtheorem{Problem}[T]{Problem}
\newcommand{\bq}{\begin{Problem}}
\newcommand{\eq}{\end{Problem}}

\newcommand{\proof}
   {\par\medbreak\noindent{\bf Proof}.\enspace}

\newcommand{\qed}{
$\Box$
\par\bigbreak}

\newlength{\lengtha} 
\newlength{\lengthb} 

\normalbaselines
\def\abstract#1{\par\bigskip
\begingroup\small
\baselineskip=12truept
\begin{center}ABSTRACT\end{center}
\par\medskip\par\noindent
\null\hfill\hbox{\vbox{\hsize=5truein\noindent#1}}
\hfill\null\par\endgroup\par}

\title{Further results on monoids acting on trees}
\author{{\bf John Rhodes}\\ 
 $ $\\ {\em Department of Mathematics, University of California, Berkeley,}\\ 
{\em California 94720, U.S.A.}\\
{\em email:} rhodes@math.berkeley.edu, BlvdBastille@aol.com\\
$ $\\
{\bf Pedro V. Silva}\\ $ $\\ {\em Centro de
Matem\'{a}tica, Faculdade de Ci\^{e}ncias, Universidade do
Porto,}\\ {\em R. Campo Alegre 687, 4169-007 Porto, Portugal}\\
{\em email:} pvsilva@fc.up.pt} \date{\today}

\begin{document}
\maketitle

\begin{center}\small
2010 Mathematics Subject Classification: 20M10, 20M20, 20M30, 20B07
\end{center}

\abstract{This paper further develops the theory  of arbitrary semigroups
acting on trees via elliptic mappings. A key tool is the Lyndon-Chiswell
length function $L$ for the semigroup $S$ which allows one to
construct a tree $T$ and an action of $S$ on $T$ via elliptic
maps. Improving on previous results, the length function of the action
will also be $L$.}

\section{Introduction}

This paper substantially  improves and extends the results in \cite{Rho1}. We
consider the case of expansions cut down to generators, which is more
compatible with geometric semigroup theory \cite{GST} and  also
allows the following major improvement over \cite{Rho1}. A
Lyndon-Chiswell function $L$ for the semigroup $S$ with generators $X$
allows one to construct a tree $T$ and an elliptic action of $S$ on
$T$. The action also gives a unique length function $L'$ on $S$. In
\cite{Rho1},$L$ and $L'$  need not be equal. However, in this paper,
by cutting to generators $X$ and performing a more refined
construction, one obtains that $L = L'$. Unfortunately, this makes the
proofs sometimes more difficult and longer than in \cite{Rho1}. Our
proofs here occasionally correct some minor errors and misprints in
\cite{Rho1} and also just refer to the arguments in \cite{Rho1} when
the proofs are the same. Applications of these results to the free
Burnside semigroups, see \cite{13McC, NF, GST}, are indicated in
Section 9. Full details of the elliptic actions of the free Burnside
semigroups will be given in a future  paper.

\section{Graphs and contractions}
\label{graphs}

Throughout the paper, morphisms and contractions shall be written on
the right. Other mappings will be written on the left.

Given a nonempty set $X$ and $n \in \N$, let $$\P_n(X) = \{ Y \subseteq
X: |Y| = n \}.$$ As usual, we identify $P_1(X)$ with $X$ to simplify
notation.

We define a {\em graph} to be an ordered pair of the form $G = (X,e)$
where
\bi
\item[(G1)]
$X$ is a nonempty set;
\item[(G2)]
$e:X \to \P_1(X) \cup \P_2(X)$ is a one-to-one mapping satisfying
$$\forall x,v \in X\; (v \in e(x) \Rw e(v) = v).$$
\ei

The elements of $\vert(G) = e\inv(P_1(X))$ are the {\em vertices} of $G$
and those of $\edge(G) = e\inv(P_2(X))$ are the {\em edges}. The
mapping $e$ fixes the vertices since
$$e(v) = w \in X \Rw e(v) = w = e(w) \Rw v = w = e(v)$$
by (G2) and injectivity of $e$,
and associates to each edge its two {\em
  adjacent} vertices. Note that this definition of (unordered) graph
excludes loops and multiple edges due to the fact of $e$ being
one-to-one.

A {\em path} in $G = (X,e)$ of length $n \in \N$ is a sequence $p = (v_0,
\ldots, v_n)$ in $\vert(G)$ such that $\{ v_{i-1}, v_i \} \in
e(\edge(G))$ for $i = 1, \ldots, n$. We say that $p$ is a path {\em
  from} $v_0$ {\em to} $v_n$. If $n = 0$ the path is said to be
{\em trivial}. The graph $G$ is said to be {\em connected} if, for all
$v,w \in \vert(G)$, there exists a path in $G$ from $v$ to $w$. 

A {\em cycle} in $G$ is a path of the form $(v_0,
\ldots, v_{n-1}, v_n)$ with $n \geq 3$, $v_n = v_0$ and $v_0,
\ldots, v_{n-1}$ all distinct. A connected graph with no cycles is
said to be a {\em tree}. 

Let $G_i = (X_i, e_i)$ be a graph  for $i = 1,2$.
A graph {\em morphism} $\p:G_1 \to G_2$ is a mapping
$\p:X_1 \to X_2$ such that 
\bi
\item[(GM1)]
$(\vert(G_1))\p \subseteq \vert(G_2)$;
\item[(GM2)]
$(e_1(x_1))\p = e_2(x_1\p)$ for every $x_1 \in X_1$.
\ei
Note that $\p$ can collapse vertices to edges: for
example, every graph has a morphism onto the trivial graph with a
single vertex.

Given a connected graph $G$, we define a distance $d$ on $\vert(G)$ by
taking $d(v,w)$ to be the length of the shortest path from $v$ to $w$
in $G$. Such a shortest path is said to be a {\em geodesic} from $v$
to $w$ and $d$ is the {\em geodesic distance} in $G$. We write
$\geo(G) = (\vert(G),d)$. If $G$ is  a tree, there is a unique
geodesic connecting $v$ and $w$ and $\geo(G)$ is a {\em hyperbolic}
metric space as considered in \cite{7Gro}.



Let $G_i = (X_i, e_i)$ be a graph  for $i = 1,2$ and let
$\geo(G_i) = (\vert(G_i),d_i)$. A {\em contraction}
$\psi: \geo(G_1) \to
\geo(G_2)$ is a 
mapping $\psi:\vert(G_1) \to
\vert(G_2)$ satisfying 
$$\forall v,w \in \vert(G_1) \; d_2(v\psi, w\psi) \leq d_1(v,w).$$

\bp
\label{mvc}
\cite[Fact 1.4]{Rho1}
Let $G_i = (X_i, e_i)$ be a graph  for $i = 1,2$. A mapping $\psi:\vert(G_1) \to
\vert(G_2)$ is a contraction if and only if $\psi$ can be extended to
a morphism $\oo{\psi}:G_1 \to G_2$. In that case, the extension is
unique.
\ep




Let $G$ be a graph. From now on, given a graph $G$, we shall identify
$G$ with its underlying set, and we shall assume that the one-to-one
mapping is denoted by $e$ and the geodesic distance by $d_G$.
We denote by $\endo(G)$ the monoid of all endomorphisms
of $G$ and by $\con(G)$ the monoid of all contractions of $\geo(G)$
into itself.

The following result is a straightforward consequence of Proposition
\ref{mvc}.

\bc
\label{mvc2}
The mapping
$$\begin{array}{rcl}
\endo(G)&\to&\con(G)\\
\p&\mapsto&\p|_{\vert(G)}
\end{array}$$
is a monoid isomorphism.
\ec

\proof
By Proposition \ref{mvc}, this mapping  is a well-defined bijection. Since
$\p(\vert(G)) \subseteq \vert(G)$ for every $\p \in \endo(G)$, it
follows that
$$(\p\p')|_{\vert(G)} = \p|_{\vert(G)}\p'|_{\vert(G)}$$
for all $\p,\p' \in \endo(G)$. Since the restriction of the identity
endomorphism is the identity contraction, our mapping is indeed a
monoid isomorphism.
\qed

\section{Elliptic $M$-trees}
\label{ellmtrees}

Let $G$ be a graph and let $M$ be a monoid with identity 1. A {\em
  (right) action} of $M$ on $G$ is a monoid homomorphism 
$$\begin{array}{rcl}
\theta:M&\to&\endo(G)\\
m&\mapsto&\theta_m
\end{array}$$
The action is {\em faithful} if $\theta$ is one-to-one.

To simplify notation, we write $xm = x\theta_m$. With this notation,
the action can be equivalently defined through the axioms:
\bi
\item[(A1)]
$(\vert(G))M \subseteq \vert(G)$
\item[(A2)]
$(e(x))m = e(xm)$
\item[(A3)]
$x(mm') = (xm)m'$
\item[(A4)]
$x1 = x$
\ei
for all $x \in G$ and $m,m' \in M$.

Note that, in view of Corollary \ref{mvc2}, the action could be
equivalently defined as a monoid homomorphism $M \to \con(G)$.

We are interested in the case of $G$ being a tree, a rooted
tree to be more precise. A {\em rooted tree} is an ordered pair of
the form $(r_0,T)$, where $T$ is a tree and $r_0 \in \vert(T)$. A
rooted tree admits a natural representation by levels $0,1,2,\ldots$ where we
locate at level (or depth) $n$ those vertices lying at distance $n$ from
$r_0$. We write then
$$\dep(v) = d(r_0,v)$$
for $v \in \vert(T)$.
Let $\oo{\N} = \N \cup \{ \omega \}$.
The 
{\em depth} of a rooted tree is defined by
$$\dep(r_0,T) = \sup\{ \dep(v) \; ; v \in \vert(T) \} \in \oo{\N}.$$ 
Given $v \in \vert(T)$, we define the {\em degree} of $v$ in $(r_0,T)$
by
$$\deg(v) = \left\{
\begin{array}{ll}
|\{ x \in \edge(T) \mid v \in e(x) \}|&\mbox{ if } v = r_0\\
|\{ x \in \edge(T) \mid v \in e(x) \}|-1&\mbox{ otherwise,}
\end{array}
\right.$$
that is, we count the number of outgoing edges if we orient them {\em
  away} from the root. A vertex of degree 0 is called a {\em leaf}. If
two vertices $v$ and $w$ are connected by an edge, we say that
$$v\mbox{ is }\left\{
\begin{array}{ll}
\mbox{a }son \mbox{ of } w&\mbox{ if } \dep(v) = \dep(w)+1\\
\mbox{the }father \mbox{ of } w&\mbox{ if }\dep(v) = \dep(w)-1
\end{array}
\right.$$
Note that a father may have many sons, but the father is always
unique. All vertices but the root have a father.

We generalize this notion with the obvious terminology. If $v_i$ is a
son of $v_{i-1}$ for $i = 1,\ldots,k$, we say that $v_k$ is a
{\em descendant} of $v_0$ and $v_0$ an {\em ancestor} of $v_k$.

A very important example is given by {\em rooted uniformly branching
  trees}:

\be
Let $n_1, \ldots, n_l \geq 1$. Up to isomorphism, the rooted uniformly branching
  tree $(r_0,T(n_l, \ldots, n_1))$ is the rooted tree of depth $l$
  such that every vertex of depth $i-1$ has degree $n_{i}$ $(i = 1,
  \ldots, l$). For example, $(r_0,T(3,2))$ can be pictured by
$$
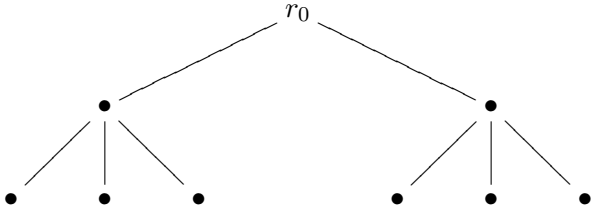
\xymatrix{
&&& r_0 \ar@{-}[dll] \ar@{-}[drr] &&&\\ 
&{\bullet} \ar@{-}[dl] \ar@{-}[d] \ar@{-}[dr] &&&&
{\bullet} \ar@{-}[dl] \ar@{-}[d] \ar@{-}[dr] & \\
{\bullet} & {\bullet}&{\bullet}&&{\bullet} & {\bullet}&{\bullet}
}$$
\ee

We can of course extend this definition to infinite cardinals in the
obvious way, as well as considering $T(\ldots,n_2,n_1)$ for an
infinite sequence. It is standard to represent $\vert(T(n_l,
\ldots,n_1))$ as
 $$\{ r_0 \} \cup (\bigcup_{i=1}^l X_i \ldots \times X_1)$$ with
$|X_i| = n_i$ for every $i$.

Let $(r_0,T), (r'_0,T')$ be rooted trees. 
An {\em elliptic} contraction $\p:(r_0,T) \to (r'_0,T')$ is a
depth-preserving contraction, that is, a contraction $\p:\vert(T) \to
\vert(T')$ satisfying
$$\forall v \in \vert(T)\; \dep(v\p) = \dep(v).$$
In view of Proposition \ref{mvc}, a bijective elliptic contraction
extends to an isomorphism of rooted trees.

\bl
\label{ell}
Let $(r_0,T)$ be a rooted tree and let $\p:\vert(T) \to \vert(T')$ be a
mapping.  Then $\p$ is an elliptic
contraction from $(r_0,T)$ into $(r'_0,T')$ if and only if
\bi
\item[(i)] $r_0\p = r'_0$;
\item[(ii)] if $v \in \vert(T)$ is the father of $w$, then $v\p$ is
  the father of $w\p$. 
\ei
\el

\proof
Assume that $\p$ is an elliptic
contraction. Then (i) holds trivially and (ii) follows from $\p$
preserving depth and being the 
restriction of a
tree morphism by Proposition \ref{mvc}.

Assume now that $\p$ satisfies conditions (i) and (ii). We extend $\p$
to $\oo{\p}: T \to T'$ as follows. Given $x \in \edge(T)$, we may
write $e(x) = \{ v,w \}$ and assume that $v$ is the father of $w$. By
(ii), it follows that $v\p$ is the father of $w\p$ and so there exists
some $x' \in \edge(T')$ such that $e(x') = \{ v\p,w\p \}$. We define
$x\oo{\p} = x'$. 

It follows from the definition that $\oo{\p}:T \to T'$ is a morphism. By
Proposition \ref{mvc}, $\p$ is a contraction. By (i), $\p$ preserves
depth 0. By (ii) and induction, $\p$ preserves
depth $n$ for each $n \in \{ 0, \ldots, \dep(r_0,T) \}$.
\qed

The set of all elliptic contractions on $(r_0,T)$ is denoted by
$\ell(r_0,T)$. This is a monoid under composition and is termed the 
{\em elliptic product} on $(r_0,T)$. 

Wreath products constitute as we shall see important examples of
elliptic products. A 
{\em partial transformation monoid} is an ordered pair of the form $(X,M)$,
where $X$ is a nonempty set and $M$ is a submonoid of the monoid
$P(X)$ of all partial transformations of 
$X$. If $M$ is a submonoid of the monoid
$M(X)$ of all full transformations of 
$X$, we say that $(X,M)$ is a {\em transformation monoid}.

Throughout the paper, given a direct product of the form $X = X_l
\times \ldots \times X_1$ 
and $i \in \{ 1,\ldots, l\}$, we shall denote by $\pi_i:X \to X_i$ the
projection on the $i$th component, and by $\pi_{[i,1]}:X \to X_i
\times \ldots \times X_1$ the 
projection on the last $i$ components.

Assume that $X = \cup_{i=1}^l (X_i \times \ldots \times 
X_1)$. For $i = 1, \ldots, l$, we define an equivalence relation
$\equiv_i$ on $X$ by
$$(x_j, \ldots,x_1) \equiv_i (x'_k, \ldots,x'_1)
\hspace{.5cm}\mbox{if}\hspace{.5cm} (i \leq 
j,k \hspace{.5cm}\mbox{and} \hspace{.5cm} x_i =
x'_i, \ldots, x_1 = x'_1).$$
Given $\p \in P(X)$, we denote by $\dom\p$ the domain of
$\p$.

A mapping $\p \in P(X)$ is said to be {\em sequential} if:
\bi
\item[(SQ1)]
$\forall i \in \{ 2, \ldots, l \} \; ( (x_i, \ldots,x_1) \in \dom\p
\Rw (x_{i-1}, \ldots,x_1) \in \dom\p)$;
\item[(SQ2)]
$\forall i \in \{ 1, \ldots, l \} \;\forall  (x_i, \ldots,x_1) \in
\dom\p \; (x_i, \ldots,x_1)\p \in X_i \times \ldots \times 
X_1$;
\item[(SQ3)]
$\forall i \in \{ 1, \ldots, l \} \; \forall x,x' \in \dom\p \; (x
\equiv_i x' \Rw x\p
\equiv_i x'\p).$
\ei
It is immediate that the composition of sequential
partial transformations of $X$ is still sequential. 

Adjoining a root $r_0$ provides a natural tree representation for 
$\cup_{i=1}^l X_i \times
\ldots \times 
X_1$.
For example, taking $X_2 = X_1 = \{ 0,1 \}$, we obtain the tree
$$
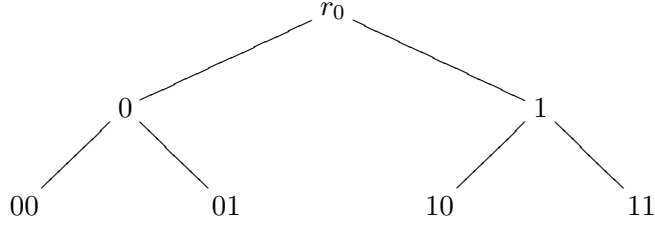
\xymatrix{
&&& r_0 \ar@{-}[dll] \ar@{-}[drr] &&&\\ 
&0 \ar@{-}[dl] \ar@{-}[dr] &&&&
1 \ar@{-}[dl] \ar@{-}[dr] & \\ 
00 && 01&&10 && 11
}$$

Given $(a_{i-1}, \ldots,a_1) \in X_{i-1} \times
\ldots \times 
X_1$ $(i \in \{ 1, \ldots, l \})$, we have $(\cdot,a_{i-1},
\ldots,a_1)\p\pi_i \in P(X_i)$.

Graphically, whenever $y\p = z$ for $y = (a_{i-1},
\ldots,a_1)$ and $X_i = \{ b_1, \ldots, b_m\}$, then we have

$$
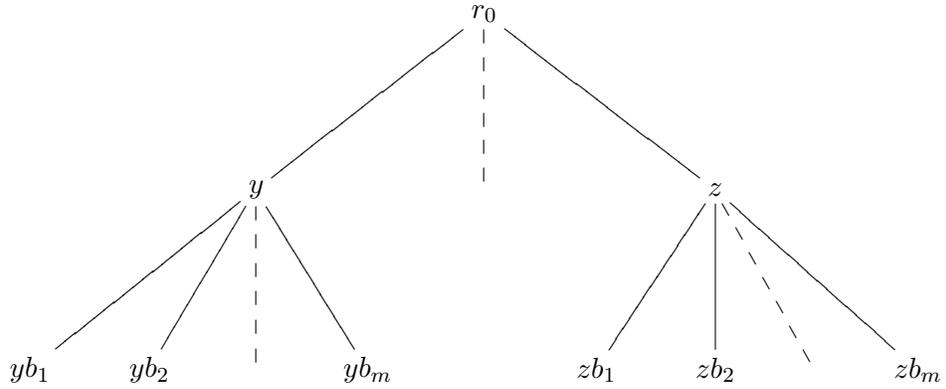
\xymatrix{
&&&& r_0 \ar@{-}[ddll] \ar@{--}[dd] \ar@{-}[ddrr] &&&&\\ &&&&&&&& \\
&&y \ar@{-}[ddll] \ar@{-}[ddl] \ar@{--}[dd] \ar@{-}[ddr] &&&&
z  \ar@{-}[ddl] \ar@{-}[dd] \ar@{--}[ddr] \ar@{-}[ddrr]
&& \\ &&&&&&&& \\ 
yb_1 & yb_2&&yb_m && zb_1 & zb_2 &&zb_m
}$$
in the tree representation and $\{ yb_1, \ldots,yb_m\}\p \subseteq \{
zb_1, \ldots,zb_m\}$. Then $(\cdot,a_{i-1},
\ldots,a_1)\p\pi_i$ is the induced partial mapping $\{ b_1, \ldots,b_m\} \to
\{ b_1, \ldots,b_m\}$ (not necessarily injective!). 

If $\p,\p' \in P(X)$ and $(a_{i-1},
\ldots,a_1)\p = (a'_{i-1},
\ldots,a'_1)$, it is easy to check \cite{6Eil,RS} that we have
\begin{equation}
\label{wp1}
(\cdot,a_{i-1},
\ldots,a_1)(\p\p'\pi_i) = ((\cdot,a_{i-1},
\ldots,a_1)\p\pi_i)((\cdot,a'_{i-1},
\ldots,a'_1)\p'\pi_i).
\end{equation}


Given partial transformation monoids $(X_l,M_l), \ldots, (X_1,M_1)$,
their {\em wreath product} is defined by
$$(X_l,M_l) \circ \ldots \circ (X_1,M_1) = (X_l \times \ldots \times
X_1, M_l \circ \ldots \circ M_1),$$
where 
$M_l \circ \ldots \circ M_1$ consists of all $\p \in P(X)$ satisfying
\bi
\item[(W1)] $\p$ is sequential;
\item[(W2)] $\p\pi_1 \in M_1$
\item[(W2)] $(\cdot,a_{i-1}, 
\ldots,a_1)\p\pi_i \in M_i$ for all $i \in \{2, \ldots, l\}$ and $(a_{i-1},
\ldots,a_1)\in \dom\p$.
\ei
More informally, $M_l \circ \ldots \circ M_1$ consists of those
partial 
self-maps of $X$ ``in sequential form 
with component action in the $M_i$'s''. Note that $M_l \circ \ldots
\circ M_1$ is a submonoid 
of $P(X)$ since the composition of sequential mappings is sequential
and by (\ref{wp1}): if $(\cdot,a_{i-1},
\ldots,a_1)\p\pi_i$ and $(\cdot,a'_{i-1},
\ldots,a'_1)\p'\pi_i$ are both in $M_i$, so is their composition
$(\cdot,a_{i-1},
\ldots,a_1)(\p\p'\pi_i)$.
Therefore $(X_l,M_l) \circ
\ldots \circ (X_1,M_1)$ is a well-defined partial transformation
monoid.

If $(X_l,M_l), \ldots, (X_1,M_1)$ are (full) transformation monoids,
their wreath product is a submonoid 
of $M(X)$.
In the case of a wreath product of two monoids with $X_1 = \{ a_1,
\ldots, a_m\}$, it is common to use
the notation $(\beta_1,\ldots,\beta_m )\alpha$ $(\alpha \in M_1, \beta_i \in
M_2)$ to denote the element of $M_2 \circ M_1$ defined by
$$(x_2,a_i)((\beta_1,\ldots,\beta_m)\alpha) = (x_2\beta_{i},a_i\alpha).$$

The wreath product of (partial) transformation monoids is associative, among
other properties. See \cite{1Arb,6Eil,RS} for more details about the
wreath product.

\bp
\label{ellwp}
For all nonempty sets $X_l, \ldots, X_1$, the monoids $M(X_l) \circ
\ldots \circ M(X_1)$ and $\ell(r_0,T(|X_l|,\ldots,|X_1|))$ are isomorphic.
\ep

\proof
We may write $T = T(|X_l|,\ldots,|X_1|)$  with $$\vert(T) = \{
r_0\} \cup (\bigcup_{i=1}^l X_i \times \ldots \times X_1).$$ 

We
consider
$$\begin{array}{rcl}
\eta:\ell(r_0,T)&\to&M(X_l) \circ
\ldots \circ M(X_1)\\
\p&\mapsto&\p\mid_{X_l \times \ldots \times X_1}.
\end{array}$$

Let $\p \in \ell(r_0,T)$. Since elliptic contractions
preserve depth, $\eta(\p) \in M(X_l 
\times \ldots \times X_1)$. It follows easily from Lemma \ref{ell}(ii)
that $\eta(\p)$ is 
sequential: if $x \equiv_i x'$, then $x,x'$ are descendants of
$x\pi_{[i,1]}$ and so $x\p,x'\p$ are descendants of
$x\pi_{[i,1]}\p$, yielding $x\p \equiv_i x'\p$.
Since (W2) and (W3) are trivially satisfied due to
$\bigcup_{i=1}^l X_i \times \ldots \times X_1 \subset \dom\p$, $\eta$ is well
defined. 

Also by Lemma \ref{ell}, the image of each $v\in \vert(T)$ by $\p$
determines the images of all its ancestors, hence $\p$ is determined
by its restriction to the leafs of $T$, i.e., $\eta(\p)$. Therefore
$\eta$ is one-to-one. 

Next let $\psi \in M(X_l) \circ
\ldots \circ M(X_1)$. We define $\p: \vert(T) \to \vert(T)$ by $r_0\p
= r_0$ and
$$(x_i, \ldots, x_1)\p = (x_l, \ldots, x_1)\psi\pi_{[i,1]},$$
the domain extension described before for a sequential map. 
If $(x_i, \ldots, x_1)$ is a son of $(x_{i-1}, \ldots, x_1)$ $(i = 2,
\ldots, l)$, then $(x_l, \ldots, x_1)\psi\pi_{[i,1]}$ is a son of
$(x_l, \ldots, x_1)\psi\pi_{[i-1,1]}$ and so  $(x_i, \ldots, x_1)\p$ is
a son of $(x_{i-1}, \ldots, x_1)\p$. Since $(x_1)\p = (x_l, \ldots,
x_1)\psi\pi_{1}$ is always a son of $r_0$, condition (ii) of Lemma
\ref{ell} holds and so $\p$ is an elliptic contraction. Since $\psi =
\eta(\p)$, we conclude that $\eta$ is onto and therefore a
bijection.

Since $(X_l \times \ldots \times X_1)\p \subseteq X_l \times \ldots
\times X_1$ for every elliptic contraction $\p$, it follows that
$\eta$ is a monoid homomorphism and therefore an isomorphism.
\qed

An {\em elliptic action} of a monoid $M$ on the rooted tree $(r_0,T)$
is a monoid homomorphism $\theta:M \to \ell(r_0,T)$.
The elliptic action is {\em faithful} if $\theta$ is one-to-one. 

We can generalize Proposition \ref{ellwp} to the case of arbitrary
wreath products of transformation monoids:

\bc
\label{ellawp}
For all transformation monoids $(X_l,M_l), \ldots, (X_1,M_1)$, the
monoid $M_l \circ 
\ldots \circ M_1$ embeds in
$\ell(r_0,T(|X_l|,\ldots,|X_1|))$. 
\ec

\proof
Write $T = T(|X_l|,\ldots,|X_1|)$.
We proved in Proposition \ref{ellwp} that 
$$\begin{array}{rcl}
\eta:\ell(r_0,T)&\to&M(X_l) \circ
\ldots \circ M(X_1)\\
\p&\mapsto&\p\mid_{X_l \times \ldots \times X_1}.
\end{array}$$
is a monoid isomorphism, its inverse being 
the mapping $\eta\inv$ that assigns to every $\psi\in M(X_l) \circ
\ldots \circ M(X_1)$ its natural domain extension $\oo{\psi}:
\vert(T) \to \vert(T)$. The
restriction of $\eta\inv$ to the submonoid $M_l \circ 
\ldots \circ M_1$ of $M(X_l) \circ 
\ldots \circ M(X_1)$ defines a faithful elliptic action of $M_l \circ 
\ldots \circ M_1$ on $\ell(r_0,T)$ and so $M_l \circ 
\ldots \circ M_1$ embeds in
$\ell(r_0,T(|X_l|,\ldots,|X_1|))$. 
\qed

However, not all submonoids of $\ell(r_0,T)$, where $T$ is a rooted
uniformly branching 
  tree, can be obtained via wreath products of transformation
  monoids, as the next example shows. 

\be
\label{ella}
Let $M$ be the submonoid of $M(\{1,2,3,4\})$ given by 
$$M = \{ (1234), (2244), (3434), (3444), (4444)\},$$
where $\p = (a_1a_2a_3a_4)$ is defined by $i\p =a_i$ for $i =
1,\ldots,4$. Then $M$ acts faithfully on $(r_0,T(2,2))$ by elliptic
contractions according to the labelling
$$
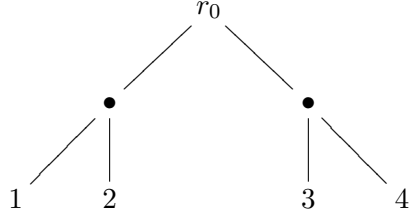
\xymatrix{
&&r_0 \ar@{-}[dl] \ar@{-}[dr] &&\\ 
&{\bullet} \ar@{-}[dl] \ar@{-}[d] &&
{\bullet} \ar@{-}[d] \ar@{-}[dr] &\\
1&2&&3&4
}$$
and so $M$ embeds in $\ell(r_0,T(2,2))$. 
However, $M$ cannot be obtained as $M_2 \circ M_1$ with $M_2,M_1$
monoids of
full transformations since  
$$|M_2 \circ M_1| = |M_2|^2|M_1|$$
and $M_1,M_2 \leq M(\{0,1\})$ implies $|M_1|,|M_2| \leq 4$.
\ee

\bigskip

Let $(r_0,T)$ be a rooted tree. 
Clearly, every $v \in (r_0,T)$ is determined by the geodesic
$$\alpha = (v = \alpha_l, \ldots, \alpha_1, \alpha_0 = r_0).$$ We call such a
geodesic a {\em ray} of $(r_0,T)$. An infinite path of the form 
$\alpha = ( \ldots, \alpha_1, \alpha_0 = r_0)$ is also said to be a ray if
$\dep(\alpha_i) = i$ for every $i \in \N$. An infinite ray is also called
an {\em end}. We denote by $\ray(r_0,T)$ the
set of all rays of $(r_0,T)$. 

Given $\alpha = (\alpha_l, \ldots, \alpha_1, \alpha_0) \in
\ray(r_0,T)$, we write $|\alpha| = l$ and 
$\dom\alpha = \{ 0, \ldots, l\}$. If $\alpha$ is infinite, we write
$|\alpha| = \omega$ and 
$\dom\alpha = \N$. In any case, given $\alpha \in
\ray(r_0,T)$ and $i \in \dom\alpha$, we denote by $\alpha_i$ the
vertex of depth $i$ in $\alpha$.

We define a partial order on $\ray(r_0,T)$ by 
$$\alpha \leq \beta \hspace{.7cm} \mbox{if} \hspace{.7cm} |\alpha|
\leq |\beta| \mbox{ and } \alpha_i = \beta_i \mbox{ for every } i
\in \dom\alpha.$$ 
We say that a ray of $(r_0,T)$ is {\em maximal} if it is maximal for
this partial order. Clearly, the maximal rays are either ends or
correspond to the leaves of the tree. We denote by $\mray(r_0,T)$ the
set of all maximal rays of $(r_0,T)$. 

We say that $(r_0,T)$ is {\em uniform} if  
all its maximal rays have the same length $l \in \oo{\N}$. In particular, 
if $(r_0,T)$ has finite depth $l$, it is uniform if all its
leaves have depth $l$. If $(r_0,T)$ has infinite depth, it is uniform if  
it has no leaves at all. The concept of maximal ray constitutes the possible 
generalization of the concept of leaf to uniform trees of infinite
depth.

Assume that $\p \in \ell(r_0,T)$. We extend $\p$ to a mapping
$\oo{\p}: \ray(r_0,T) \to \ray(r_0,T)$ by
$$\alpha\oo{\p} = (\ldots, \alpha_2\p, \alpha_1\p).$$
It follows from Lemma \ref{ell} that $\oo{\p}$ is
well defined. Identifying finite rays with vertices as usual,
$\oo{\p}$ can be seen as an extension of $\p$. 
We shall denote
$\oo{\p}$ by $\p$ when no confusion arises. 

As a particular case, if $M$
acts elliptically on $(r_0,T)$, we can extend this action to
$\ray(r_0,T)$ by
$$\alpha m = (\ldots, \alpha_2m, \alpha_1m) \hspace{1.5cm}(\alpha \in
\ray(r_0,T), \; m \in M).$$ Note that 
$$\alpha 1 = \alpha, \quad \alpha (mm') = (\alpha m)m'$$
for all $\alpha \in \ray(r_0,T)$ and $m,m' \in M$, hence we can
properly speak
of an action of $M$ on $\ray(r_0,T)$.

Let  $(r_0,T)$ be a uniform
rooted tree and let $\alpha \in \mray(r_0,T)$. An elliptic
action of $M$ on  $(r_0,T)$ is 
said to be $\alpha$-{\em transitive} if $$\vert(T) = \alpha M =
\cup_{i \in \dom\alpha} 
\alpha_iM.$$
If $(r_0,T)$ has finite depth $l$, it should be clear that the elliptic
action of $M$ on  $(r_0,T)$ is 
$\alpha$-transitive if and only if $\alpha_lM$ is the set of leaves of
$(r_0,T)$. Indeed, since the action of $M$ is depth-preserving, only
leaves can be sent to leaves. On the other hand, transitivity at the
deepest level clearly implies transitivity on the upper levels in view
of Lemma \ref{ell}.

An {\em elliptic $M$-tree} is a structure of the form $\chi =
(r_0,T,\alpha,\theta)$, where
\bi
\item[(E1)] $(r_0,T)$ is a  uniform
rooted tree;
\item[(E2)] $\alpha \in \mray(r_0,T)$;
\item[(E3)] $\theta:M \to \ell(r_0,T)$ is an $\alpha$-transitive action.
\ei
We say that $\chi$ is a {\em faithful} elliptic $M$-tree if $\theta$
is one-to-one. We say $\chi$ is a {\em strongly faithful} elliptic $M$-tree if
$$\alpha m = \alpha m' \Rw m = m' \hspace{.5cm} \mbox{for all }m,m' \in M.$$
We shall omit
$\theta$ from the representation of $\chi$ when no confusion arises
from doing so.

Let $\chi =
(r_0,T,\alpha)$, $\chi' =
(r'_0,T',\alpha')$ be elliptic $M$-trees. A {\em morphism} $\p:\chi \to 
\chi'$ of elliptic $M$-trees is an elliptic contraction $\p:(r_0,T) \to
(r'_0,T')$ such that:
\bi
\item[(EM1)] $\alpha\p = \alpha'\p$;
\item[(EM2)] $\forall v \in \vert(T)\; \forall m \in M$ $(vm)\p = (v\p)m$.
\ei
If $\p$ is bijective, we say it is an {\em isomorphism} of elliptic $M$-trees.


Given a transformation monoid $(X,M)$ and $x_0 \in X$, we say that $M$
acts transitively on $(X,x_0)$ if $x_0M = X$. A {\em pointed
  transformation monoid} is a triple of the form $(X,x_0,M)$, where
$M$ acts transitively on $(X,x_0)$.

\bc
\label{tellawp}
Let $\ldots,(X_2,x_2,M_2),
(X_1,x_1,M_1)$ be pointed transformation monoids. Then
$(r_0,T(\ldots,|X_2|,|X_1|),(\ldots,x_2,x_1))$ is a faithful elliptic
$(\ldots \circ M_2 \circ M_1)$-tree.
\ec

\proof
Axioms (E1) and (E2) are trivially verified. 
Let $\alpha = (\ldots,x_2, x_1)$. 
We observed in the proof of Corollary \ref{ellawp} that the
restriction of $\eta\inv$ as defined in Proposition \ref{ellwp} to
the submonoid $\ldots \circ M_2 \circ M_1$ of $\ldots \circ M(X_2)
\circ M(X_1)$ defines a faithful elliptic action of $ 
\ldots \circ M_2 \circ M_1$ on $\ell(T(\ldots,|X_2|,|X_1|))$. A straightforward
induction on $i$ proves that this action is $\alpha$-transitive: indeed, it is
enough to show that, given $(w_i, \ldots,w_1) \in X_i \times \ldots
\times X_1$, there exists $\p_i \in M_i \circ  
\ldots \circ M_1$ such that $(x_i, \ldots,x_1)\p_i = (w_i,
\ldots,w_1)$. The case $i = 1$ follows from 
$M_1$ acting transitively on $(X_1,x_1)$.
Assume that $(x_i, \ldots,x_1)\p_i = (w_i,
\ldots,w_1)$ for some $\p_i \in M_i \circ  
\ldots \circ M_1$. Since $x_{i+1}\xi = w_{i+1}$ for some $\xi
\in M_{i+1}$, we can define $\p_{i+1} \in M_{i+1} \circ  
\ldots \circ M_1$ by $$\p_{i+1} = (\xi, \ldots, \xi)\p_i.$$
It follows that 
$$(x_{i+1}, \ldots,x_1)\p_i = (x_{i+1}\xi,w_{i},
\ldots,w_1) = (w_{i+1},
\ldots,w_1)$$ and so (E3) holds as required.
\qed

Example \ref{ella} shows also that not all faithful elliptic
$M$-trees on a rooted uniformly branching 
  tree can be obtained via wreath products,
the action of $M$ on $(r_0,T(2,2))$ being obviously
$\alpha$-transitive for the ray defined by the leaf 1.

\section{Length functions}

Let $M$ be a monoid. Let $\oo{\N} = \N \cup \{ \omega \}$ have the
obvious ordering. A {\em length function} for $M$ is a function $D:
M \times M \to \oo{\N}$ satisfying the axioms
\bi
\item[(L1)]
$D(m,m') = D(m',m)$
\item[(L2)]
$D(m',m'') \leq D(m,m)$
\item[(L3)]
$D(m',m'') \leq D(m'm,m''m)$
\item[(L4)] 
$D(m,m'') \geq \min\{ D(m,m'), D(m',m'') \}$ (isoperimetric
inequality)
\ei
for all $m,m',m'' \in M$.

Note that, by (L2), $D$ has a maximum $l \in \oo{\N}$ and
\beq
\label{lfqu0}
D(m,m) = l \quad \mbox{for every } m \in M.
\eeq
Moreover, for any submonoid $M'$ of $M$, the restriction of $D$ to $M'
\times M'$ is a
length function for $M'$.

We recall that a {\em quasi-ultrametric} on a set $X$ is a function
$d:X \times X \to \mathbb{R}$ satisfying the axioms
\bi
\item[(Q1)]
$d(x,x') = d(x',x)$
\item[(Q2)]
$d(x,x) = 0$
\item[(Q3)]
$d(x,x'') \leq \max\{ d(x,x'), d(x',x'') \}$
\ei
for all $x,x',x'' \in X$.

Bounded length functions can be related to quasi-ultrametrics as follows:

\bp
\label{lfqu}
Let $M$ be a monoid and let $D:
M \times M \to \N$ be a bounded function with maximum $l \in \N$. Define
$d:M \times M \to \N$ by $d(m,m') = 2l - 2D(m,m')$. Then $D$ is
a length function for $M$ if and only if the following conditions
hold:
\bi
\item[(i)] $d$ is a quasi-ultrametric;
\item[(ii)] $d(m'm,m''m) \leq d(m',m'')$ for all $m,m',m'' \in M$.
\ei
\ep

\proof
It is immediate that (Q1) $\iff$ (L1), (Q2) $\iff$ (L2), (Q3) $\iff$
(L4) and (ii) $\iff$ (L3). 
\qed

Let $(r_0,T)$ be a rooted tree and consider the partial order $\leq$
defined on $\ray(r_0,T)$ in Section \ref{ellmtrees}.
It is immediate that $(\ray(r_0,T),
\leq)$ is a $\wedge$-semilattice and $\alpha \wedge \beta$ is defined by
$$\alpha \wedge \beta = \left\{
\begin{array}{ll}
\alpha&\mbox{ if } \alpha = \beta\\
(\alpha_k, \ldots, \alpha_0)&\mbox{ if } k = \max \{ i \in \N\mid
\alpha_i = \beta_i \}.
\end{array}
\right.$$
Thus $(\ray(r_0,T),\wedge)$ is a
semilattice with zero $(r_0)$ (in the semigroup theory
sense). Identifying vertices with finite rays, we can say that
$\vert(T)$ is a $\wedge$-subsemilattice of $(\ray(r_0,T),
\leq)$, considering the ancestor partial ordering on $\vert(T)$:
$$v \leq w \hspace{.7cm} \mbox{ if } \hspace{.7cm} v = w \mbox{ or } v
\mbox{ is an 
  ancestor of }w.$$ 


\bl
\label{finch}
Let $(r_0,T)$ be a rooted tree and $\alpha \in \ray(r_0,T)$. Write
$(\alpha ] = \{ \beta \in \ray(r_0,T) \mid \beta \leq \alpha \}$. Then:
\bi
\item[(i)] $(\alpha ]$ is a chain;
\item[(ii)] if $\alpha$ is finite, $(\alpha ]$ is finite.
\ei
\el

\proof
(i) and (ii) follow from 
$$(\alpha ] = \{ \alpha \} \cup \{ (\alpha_i, \ldots, \alpha_0)
\mid i = 0, \ldots, 
\dom\alpha \}.$$
\qed


\bl
\label{isoi}
Let $(r_0,T)$ be a rooted tree. Then
$$|\alpha\wedge \alpha''| \geq \min\{ |\alpha\wedge \alpha'|,
|\alpha'\wedge \alpha''| 
\}$$
for all $\alpha,\alpha',\alpha'' \in  \ray(r_0,T)$.
\el

\proof
We have $\alpha\wedge \alpha', \alpha'\wedge \alpha'' \leq
\alpha'$. Since $(\alpha' ]$
is a chain by Lemma \ref{finch}(i) and $\wedge$ is commutative, we
may assume that 
$\alpha\wedge \alpha' \geq \alpha'\wedge \alpha''$. Thus
$$\alpha\wedge \alpha'' \geq \alpha\wedge \alpha' \wedge \alpha'' =
(\alpha\wedge 
\alpha') \wedge (\alpha'\wedge \alpha'') = 
\alpha'\wedge \alpha''$$
and so $$|\alpha\wedge \alpha''| \geq |\alpha'\wedge
\alpha''| \geq \min\{ |\alpha\wedge 
\alpha'|, |\alpha'\wedge \alpha''| \}$$
as claimed.
\qed

\bp
\label{islf}
Let $\chi = (r_0,T,\alpha)$
be an elliptic $M$-tree and
define a mapping $D_{\chi}: M \times M \to \oo{\N}$ by
$$D_{\chi}(m,m') = |\alpha m \wedge \alpha m'|.$$
Then:
\bi
\item[(i)] $D_{\chi}$ is a length function for $M$;
\item[(ii)] if $\dep(r_0,T) = l \in \N$, the 
quasi-ultrametric $d$ associated to $D_{\chi}$ satisfies
$$d(m,m') = d_T(\alpha_lm,\alpha_lm');$$
\item[(ii)] if $\dep(r_0,T) = l \in \N$ and $\chi$ is strongly
  faithful, then $d$ is an ultrametric.
\ei
\ep

\proof
(i) Since $\wedge$ is commutative, axiom (L1) is trivially satisfied.

Since $\alpha, \alpha m$ are maximal rays of $(r_0,T)$, we have
$$\begin{array}{lll}
D_{\chi}(m',m'')&=&|\alpha m' \wedge \alpha m''| \leq \dep(r_0,T) =
|\alpha m|\\
&=&|\alpha m \wedge \alpha m| = D_{\chi}(m,m)
\end{array}$$
and so (L2) holds.

For all $\beta, \beta' \in \ray(r_0,T)$, we have that
$$\beta_i = \beta'_i\mbox{ for } i = 0, \ldots,k \hspace{.7cm}\Rw
\hspace{.7cm} \beta_im = \beta'_im\mbox{ for } i = 0, \ldots,k.$$
Thus $|\beta\wedge\beta'| \leq |\beta m\wedge\beta'm|$ and so
$$D_{\chi}(m',m'') = |\alpha m' \wedge \alpha m''| \leq |\alpha m'm
\wedge \alpha m''m| = D_{\chi}(m'm,m''m).$$
Thus (L3) holds. Since (L4) follows fromm Lemma \ref{isoi}, 
$D_{\chi}$ is a length function for $M$.

(ii) Assume that $\dep(r_0,T) = l \in \N$. By (L2), $l$ is the maximum
value of $D_{\chi}$. 
By Proposition \ref{lfqu}, the associated quasi-metric is defined by
$d(m,m') = 2l - 2D_{\chi}(m,m')$. Let $v$ be the deepest vertex of
$\alpha m \wedge \alpha m'$. Since $v$ lies in
the 
geodesics $r_0 \edg \alpha_lm$,  $r_0 \edg \alpha_lm'$ and $\alpha_lm
\edg \alpha_lm'$,  
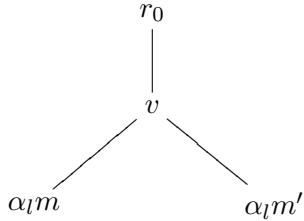
$$\xymatrix{
&r_0 \ar@{-}[d]  &\\
&{v} \ar@{-}[dl] \ar@{-}[dr] &\\ 
\alpha_lm&&\alpha_lm'
}$$
we
obtain
$$\begin{array}{lll}
d(m,m')&=&2l - 2D_{\chi}(m,m')\\
&=&d_T(r_0,\alpha_lm) + d_T(r_0,\alpha_lm')-2|\alpha m \wedge \alpha
m'|\\
&=&d_T(r_0,\alpha_l m) - d_T(r_0,v)
+ d_T(r_0,\alpha_l m') - d_T(r_0,v)\\
&=&d_T(\alpha_lm,v) + d_T(\alpha_lm',v) =
d_T(\alpha_lm,\alpha_lm').
\end{array}$$

(iii) Assume that $\dep(r_0,T) = l \in \N$ and $\chi$ is strongly
  faithful. Let $m,m' \in M$ be such that $d(m,m') = 0$. By (ii), we
  have $d_T(\alpha_lm,\alpha_lm') = 0$ and so $\alpha_lm =
  \alpha_lm'$. Hence $\alpha m =
  \alpha m'$ and so $m = m'$ since $\chi$ is strongly
  faithful. Therefore $d$ is an ultrametric.
\qed

In view of Corollary \ref{tellawp}, it is interesting to analyze the
particular case of wreath products. The canonical length function for
two mappings $\p,\p' \in M_l \circ  
\ldots \circ M_1$ measures the maximum number of components (from
right to left) where $(x_l,\ldots,x_1)\p$ and $(x_l,\ldots,x_1)\p'$
coincide: 

\bc
\label{etwp}
Let $\ldots, (X_2,x_2,M_2), 
(X_1,x_1,M_1)$ be pointed transformation monoids and let
$\chi = (r_0,T(\ldots,|X_2|,|X_1|),(\ldots,x_2,x_1))$ be the
corresponding faithful elliptic 
$(\ldots \circ M_2 \circ M_1)$-tree. Then
$$D_{\chi}(\p,\p') = \max\{ i :
(\ldots,x_2,x_1)\p\pi_{[i,1]} = (\ldots,x_2,x_1)\p'\pi_{[i,1]} \}.$$  
\ec

\proof
We prove the finite case for pointed transformation monoids 
$(X_l,x_l,M_l), \ldots,
(X_1,x_1,M_1)$.
Let $(x'_l,\ldots,x'_1),(x''_l,\ldots,x''_1) \in
\vert(T(|X_l|,\ldots,|X_1|))$ (we identify $r_0$ with the empty
sequence).   
Since 
$$(x'_l,\ldots,x'_1) \wedge (x''_l,\ldots,x''_1) = (x'_k,\ldots,x'_1)$$
where $$k = \max\{ i \in \{ 0,\ldots,l\}: (x'_i,\ldots,x'_1) =
(x''_i,\ldots,x''_1) \},$$ we have 
$$\begin{array}{lll}
D_{\chi}(\p,\p')&=&|(x_l,\ldots,x_1)\p \wedge (x_l,\ldots,x_1)\p'|\\
&=&\max\{ i \in \{ 0,\ldots,l\}: (x_l,\ldots,x_1)\p\pi_{[i,1]} =
(x_l,\ldots,x_1)\p'\pi_{[i,1]} \}. 
\end{array}$$
\qed

\bl
\label{isoet}
Let $\chi,\chi'$
be elliptic $M$-trees. Then $D_{\chi} = D_{\chi'}$ if and only if 
$\chi \cong \chi'$.
\el

\proof
Let $\chi = (r_0,T,\alpha)$, $\chi' = (r'_0,T',\alpha')$. Assume that
 $D_{\chi} = D_{\chi'}$.
Note that
$$\dep(r_0,T) = \max D_{\chi} = \max D_{\chi'} = \dep(r'_0,T').$$
Write $l = \dep(r_0,T)$.

We define a mapping
$$\begin{array}{rcl}
\p:\vert(T)&\to&\vert(T')\\
\alpha_im&\mapsto&\alpha'_im
\end{array}$$
where $i \in \dom\alpha$ and $m \in M$.
Since the action of $M$ on $(r_0,T)$ is $\alpha$-transitive, we have
$$\{ v \in \vert(T) \mid dep(v) = i \} = \alpha_iM.$$ 
If $\alpha_im = \alpha_im'$, then 
$$\begin{array}{lll}
|\alpha' m \wedge \alpha' m'|&=&D_{\chi'}(m,m') = D_{\chi}(m,m')\\
&=&|\alpha m \wedge \alpha m'| \geq i,
\end{array}$$
hence $\alpha'_im = \alpha'_im'$. Since 
$$\{ v \in \vert(T') \mid dep(v) = i \} = \alpha'_iM,$$ it
follows that $\p$ is well defined and onto. By symmetry, $\p$ is also
one-to-one.  

Clearly, $r_0\p = (\alpha_0\cdot 1)\p = \alpha'_0\cdot 1 =
r'_0$. Assume that $v$ is 
the father of $w = \alpha_im$. Then $v = \alpha_{i-1}m$ since the
action of $M$ on $(r_0,T)$ is elliptical, hence
$v\p = \alpha'_{i-1}m$ is the father of $w\p = \alpha'_{i}m$. By Lemma
\ref{ell}, $\p$ is an elliptic contraction from $(r_0,T)$ onto
$(r'_0,T')$ and therefore an isomorphism of rooted trees.

Since
$$\alpha\p = (\ldots,\alpha_1\cdot 1,\alpha_0\cdot 1)\p = 
(\ldots,\alpha'_1\cdot 1,\alpha'_0\cdot 1) = \alpha'$$
and
$$ ((\alpha_im)m')\p = (\alpha_imm')\p = \alpha'_imm' = ((\alpha_im)\p)m'$$
for all $m,m' \in M$, axioms (EM1) and (EM2) are satisfied and so
$\p:\chi \to \chi'$ is an isomorphism of elliptic $M$-trees. 

Conversely, assume that $\p:\chi \to \chi'$ is an isomorphism of
elliptic $M$-trees. For all $m,m' \in M$, we have 
$$\begin{array}{lll}
D_{\chi}(m,m')&=&|\alpha m \wedge \alpha m'| = |(\alpha m)\p \wedge
(\alpha m')\p|\\ 
&=&|(\alpha \p)m \wedge (\alpha \p)m'| = |\alpha' m \wedge \alpha' m'|\\
&=&D_{\chi'}(m,m'),
\end{array}$$
hence $D_{\chi} = D_{\chi'}$ and the lemma holds.
\qed
 
A proof for the following theorem can be found in \cite{Rho1}, but the
important role played by the Chiswell construction in it makes it
worthwhile to include it here.

\bt
\label{lftree}
\cite[Theorem 1.12]{Rho1}
Let $M$ be a monoid and let $D:M \times M \to \oo{\N}$ be a mapping. Then
the following conditions are equivalent: 
\bi
\item[(i)] $D$ is a length function for $M$;
\item[(ii)] $D = D_{\chi}$ for some elliptic $M$-tree $\chi$. 
\ei
Moreover, if the conditions hold, 
$\chi$ is unique up to isomorphism.
\et

\proof
Assume that $D$ is a length function for $M$.
We adapt the important {\em Chiswell construction} of \cite{5Chi} as
follows. By (L2), there exists 
$l = \max D \in \oo{\N}$. Let
$$P = \left\{
\begin{array}{ll}
\{ 0, \ldots, l \} \times M&\mbox{ if }l \in \N\\
\N \times M&\mbox{ if }l = \omega
\end{array}
\right.$$ 
and define a relation $\sim$ on $P$
by $$(k,m)\sim (k',m') \hspace{1cm} \mbox{if} \hspace{1cm} k
= k' \mbox{ and } 
D(m,m') \geq k.$$
We show that $\sim$ is an equivalence relation on $P$.

In fact, $\sim$ is reflexive by (\ref{lfqu0}), and symmetric by
(L1). Transitivity follows 
from the isoperimetric inequality (L4). Let $[k,m]$ denote the $\sim$
equivalence class of $(k,m)$. We define a graph $T$ by
\bi
\item[]
$\vert(T) = P/\sim$,
\item[]
$\edge(T) = \{ [k,m] \edg [k+1,m] ; \; (k,m) \in P, \; k < l \}$.
\ei

It follows from the definitions that
$$\forall m,m' \in M, \;  (0,m) \sim (0,m').$$
Let $r_0 = [0,1]$. Since
\begin{equation}
\label{lftree1}
[k,m] \edg  \ldots \edg [1,m] \edg[0,m] = r_0
\end{equation}
is a path in $T$ for every $(k,m) \in P$, $T$ is a connected graph.

We show next that
\begin{equation}
\label{lftree2}
[k,m] \edg[k+1,m'] \in \edge(T) \iff (k,m) \sim (k,m')
\end{equation}
holds for all $k \in \{ 0, \ldots, l-1 \}$ and $m,m' \in M$. 
Indeed, if $[k,m] \edg[k+1,m'] \in \edge(T)$ then $[k,m] =
[k,m'']$ and 
$[k+1,m'] = [k+1,m'']$ for some $m'' \in
M$. Hence $D(m,m'') \geq k$ and 
$D(m',m'') \geq k+1$, yielding
$$D(m,m') \geq \min\{ D(m,m''),
D(m',m'') \} \geq k$$ 
by (L4). Thus $(k,m) \sim (k,m')$. The converse implication is trivial,
therefore (\ref{lftree2}) holds.

We can prove now that $T$ is a tree. Assume that $T$ has a cycle $C$ and let
$[k,m]$ be a vertex in $C$ with $k$ maximum. Let $[k',m']$
and $[k'',m'']$ be 
its adjacent vertices in $C$. By maximality of $k$, we have $k' =
k'' = k-1$, 
hence $(k-1,m') \sim (k-1,m) \sim (k-1,m'')$ by
(\ref{lftree2}) and so $[k',m'] 
= [k'',m'']$, contradicting $C$ being a
cycle. Therefore $T$ is a tree and so 
$(r_0,T)$ is a rooted tree.

Clearly, (\ref{lftree1}) is a ray for every vertex $[k,m]$. If $l \in
\N$, then $(r_0,T)$ has finite depth $l$ and it is uniform since
$\mray(r_0,T)$ consists of all paths of the form
$([l,m], \ldots,[1,m],[0,m] = r_0)$
with $m \in M$. If $l = \omega$, all rays must have infinite length
since $[k+1,m] \edg [k,m]$ is an edge for every $(k,m) \in P$, hence
$(r_0,T)$ is uniform as well.

We define a mapping $\eta:\vert(T) \times M \to \vert(T)$ by
$$\eta([k,m],m') = [k,m]m' = [k,mm'].$$
Note that $$[k,m] = [k,n] \Rw D(m,n) \geq k
\Rw D(mm',nm') \geq k \Rw 
[k,mm'] = [k,nm']$$
by (L3) and so the mapping is well defined.  

Clearly, $r_0m = r_0$ for every $m \in M$. On the other hand, if
$[k-1,m']$ is the father of $[k,m']$, then $[k-1,m']m$ is the father
of $[k,m']m$ and so $\eta$ induces a mapping $$\begin{array}{rcl} 
\theta:M&\to&\ell(r_0,T)\\
m&\mapsto&\eta(\cdot,m)
\end{array}$$ 
by Lemma \ref{ell}. Since $\theta$ is a monoid homomorphism due to
$$[k,m]1 = [k,m],\quad [k,m](m'm'') = ([k,m]m')m'',$$
it follows that $\theta$ is an elliptic action of $M$ on $(r_0,T)$.

Let $\alpha \in \ray(r_0,T)$
be defined by 
$$|\alpha| = l,\quad \alpha_i = [i,1] \; (i \in \dom\alpha).$$ 
Since $[i,m] = \alpha_im$ for every $m \in M$, the action $\theta$ is
$\alpha$-transitive. Thus $\chi = (r_0,T,\alpha,\theta)$ is an elliptic
$M$-tree. We show that $D = D_{\chi}$. 


For all $m,m' \in M$, we have 
$$\begin{array}{lll}
D_{\chi}(m,m')&=&|\alpha m \wedge \alpha m'| = |([i,m])_i \wedge
([i,m'])_i| = \sup\{ i \in \dom\alpha: [i,m] = [i,m'] \}\\
&=&\sup\{ i \in \dom\alpha: D(m,m') \geq i \} = D(m,m')
\end{array}$$
and so $D = D_{\chi}$. Therefore (ii) holds.

(ii) $\Rw$ (i). By Proposition \ref{islf}.

The uniqueness of $\chi$ up to isomorphism follows from Lemma \ref{isoet}.
\qed

We consider now the case of strongly faithful elliptic
$M$-trees. A length function $D:M \times M \to \oo{\N}$ is said to be
{\em strict} if 
\bi
\item[(L5)] $D(m',m'') = D(m,m) \Rw m' = m''$ for all $m,m',m'' \in
  M$;
\ei

\bc
\label{flftree}
Let $M$ be a monoid and let $D:M \times M \to \oo{\N}$ be a mapping. Then
the following conditions are equivalent: 
\bi
\item[(i)] $D$ is a strict length function for $M$;
\item[(ii)] $D = D_{\chi}$ for some strongly faithful elliptic $M$-tree $\chi$. 
\ei
Moreover, if the conditions hold, 
$\chi$ is unique up to isomorphism.
\ec

\proof
(i) $\Rw$ (ii). Assume that (i) holds. By Theorem \ref{lftree}, $D =
D_{\chi}$ for the elliptic $M$-tree $\chi = (r_0,T,\alpha,\theta)$
defined in its proof. We show that $\chi$ is strongly  faithful. Indeed, let
$m,m' \in M$. Suppose that $\alpha m = \alpha m'$. Then $[k,m] =
[k,m']$ for every $k \in 
\dom\alpha$ and so $D(m,m') \geq k$ for every
$k \in
\dom\alpha$. It follows that $D(m,m') = l = \max D = D(m,m)$ and so $m
= m'$ by (L5). Thus $\chi$ is strongly  faithful.

(ii) $\Rw$ (i). Assume that (ii) holds for $\chi = (r_0,T,\alpha)$. By
Theorem \ref{lftree}, we 
only need to show that $D_{\chi}$ satisfies (L5). Suppose that
$D(m',m'') = D(m,m)$ for some $m,m',m'' \in
  M$. Hence
$$\begin{array}{lll}
|\alpha m' \wedge \alpha m''|&=&D_{\chi}(m',m'') =
D(m',m'') = D(m,m) = D_{\chi}(m,m)\\
&=&|\alpha m \wedge \alpha m| = |\alpha m| = \dep(r_0,T)
\end{array}$$ and so $\alpha m' = \alpha m''$. Since $\chi$ is
strongly faithful, we get $m' = m''$ and so (L5) holds.

The uniqueness of $\chi$ up to isomorphism follows from Lemma \ref{isoet}.
\qed

We end this section by associating a length function to any wreath
product of partial transformation monoids. To simplify notation, we
present just the infinite 
case, the finite one being absolutely similar.

\bp
\label{turq}
Let $\ldots, (X_2,M_2), (X_1,M_1)$ be partial transformation monoids
and let $$(X,M) = \ldots \circ (X_2,M_2) \circ (X_1,M_1) = (
\ldots \times X_2 \times X_1, \; \ldots \circ  M_2 \circ M_1)$$
be their wreath product. Let $D:M \times M \to \oo{\N}$ be defined by
$$D(\p,\psi) = \sup\{ j\in \N \mid \p|_{X_j \times \ldots \times X_1}
= \psi|_{X_j \times \ldots \times X_1}\}.$$
Then $D$ is a strict length function for $M$.
\ep

\proof
Axioms (L1), (L2) and (L5) hold trivially. 

Let $\p,\psi, \mu \in M$. Since $(X_j \times \ldots \times X_1)\theta
\subseteq X_j \times \ldots \times X_1$ for every $\theta \in M$,
$\p|_{X_j \times \ldots \times X_1} 
= \psi|_{X_j \times \ldots \times X_1}$ implies $(\p\mu)|_{X_j \times
  \ldots \times X_1} = (\psi\mu)|_{X_j \times \ldots \times
  X_1}$. Thus $D(\p,\psi) \leq D(\p\mu,\psi\mu)$
 and (L3) holds.

Finally, $\p|_{X_j \times \ldots \times X_1} \neq \mu|_{X_j \times
  \ldots \times X_1}$ implies either $\p|_{X_j \times \ldots \times
  X_1} \neq \psi|_{X_j \times 
  \ldots \times X_1}$ or $\psi|_{X_j \times \ldots \times X_1}$ $\neq
\mu|_{X_j \times 
  \ldots \times X_1}$, hence
$D(\p,\mu) \geq \min\{ D(\p,\psi), D(\psi,\mu) \}$ and (L4) holds.
\qed

\section{Expansions}

Let $\M$ denote the category of all monoids. A {\em monoid expansion}
is a functor $F:\M \to \M$ preserving surjective morphisms such that
there exists a natural 
transformation $\eta$ from the functor $F$ to the identity functor
with $\eta_M$ surjective for each $M \in \M$. 

That is, $F$ assigns to each monoid $M$ a monoid $F(M)$ and a
surjective morphism $\eta_M:F(M) \to M$, and to
each monoid homomorphism $\p:M \to N$ a monoid homomorphism
$F(\p):F(M) \to F(N)$ satisfying:
\bi
\item[(E1)] if $\p$ is surjective, so is $F(\p)$;
\item[(E2)] if $\p = \id_M$, then $F(\p) = \id_{F(M)}$;
\item[(E3)] if $\p:M \to M'$, $\p':M' \to M''$ are morphisms, then
  $F(\p\p') = F(\p)F(\p')$; 
\item[(E4)] if $\p:M \to N$ is a morphism, then the following diagram
  commutes:
$$\xymatrix{
F(M) \ar[rr]^{F(\p)} \ar[dd]_{\eta_M}&&F(N) \ar[dd]^{\eta_N} \\ &&\\
M \ar[rr]_{\p} && N
}$$
\ei
{\em Semigroup expansions} are defined analogously.

An element $a \in M$ is said to be {\em aperiodic} if $a^{n+1} = a^n$
for some $n \in \N$. A morphism $\p:M \to N$ is said to be {\em
  aperiodic} if, whenever $a \in 
N$ is aperiodic, all elements in $a\p\inv$ are also aperiodic.
The expansion $F$ is said to be {\em aperiodic} if the morphism
$\eta_M$ is aperiodic for every monoid $M$.

We define now the Rhodes expansion for monoids, omitting the expansion
of morphisms. The reader is referred to \cite{20Til, 15Rho, 4BR, 6Eil,
  14Rei, 16Rho, RS} for more details. 

The {\em $\L$-preorder} on a monoid $M$ is defined by
$$a \leq_{\L} b \hspace{.7cm} \mbox{if} \hspace{.7cm} a \in Mb.$$
This preorder is clearly compatible with multiplication on the right:
$$\forall a,b,m \in M \; (a \leq_{\L} b \Rw am \leq_{\L} bm).$$
The Green relation $\L$ can of course be defined by
$$a \L b \hspace{.7cm} \mbox{if} \hspace{.7cm} a \leq_{\L} b \mbox{ and
} b \leq_{\L} a.$$
The {\em strict $\L$-order} on $M$ is defined by
$$a <_{\L} b \hspace{.7cm} \mbox{if} \hspace{.7cm} a \in Mb \mbox{ and
} b \notin Ma,$$
i.e., $<_{\L} \; = \; \leq_{\L} \setminus \L$.

The $\R$- and $\J$-versions are defined similarly. In particular,
$$a <_{\J} b \hspace{.7cm} \mbox{if} \hspace{.7cm} a \in MbM \mbox{ and
} b \notin MaM.$$

Given a finite chain of the form
$$\sigma = (m_k \leq_{\L} \ldots \leq_{\L} m_1 \leq_{\L} m_0)$$
in $M$, we define a chain
$$\lm(\sigma) = (m_{i_l} <_{\L} \ldots <_{\L} m_{i_1} <_{\L} m_{i_0})$$
by keeping the leftmost term in each $\L$-class of terms of
$\sigma$. Thus
$$\sigma = (m_{i_l} \L m_{i_l-1} \L \ldots \L m_{i_{l-1}+1} <_{\L}
m_{i_{l-1}} \ldots m_{i_{0}+1} <_{\L}
m_{i_0} \L \ldots \L m_0).$$

We define the {\em Rhodes expansion} $\E(M)$ of $M$ to be the set of
all finite chains of the form
$$m_k <_{\L} \ldots <_{\L} m_1 <_{\L} m_0 = 1$$
with $k \geq 0$ and $m_i \in M$. The product of two chains 
$$\sigma = (m_k <_{\L} \ldots <_{\L} m_1 <_{\L} m_0 = 1), \quad
\tau = (m'_l <_{\L} \ldots <_{\L} m'_1 <_{\L} m'_0 = 1)$$
is defined by
$$\sigma\tau = \lm(m_km'_l \leq_{\L} \ldots \leq_{\L} m_1m'_l
\leq_{\L} m_0m'_l = m'_l <_{\L} \ldots <_{\L} m'_1 <_{\L} m'_0
= 1).$$
Note that the product is well defined since $\leq_{\L}$ is right
compatible. It turns out that $\E(M)$ is a monoid having the trivial
chain $(m_0 = 1)$ as identity.

The surjective morphisms
$\eta_M: \E(M) \to M$ are defined by $$(m_k <_{\L} \ldots <_{\L} m_1
<_{\L} m_0 = 1)\eta_M = m_k.$$

It follows from the definition of $\leq_{\L}$ that the elements of
$\E(M)$ are precisely the finite chains of the form
$$x_k\ldots x_2x_1 <_{\L} x_{k-1}\ldots x_2x_1 <_{\L} \ldots <_{\L}
x_2x_1 <_{\L} x_1 <_{\L} 1$$ with $x_1, \ldots, x_k \in M$. Moreover,
\beq
\label{gene}
(x_k\ldots x_2x_1 <_{\L}  \ldots <_{\L}
x_2x_1 <_{\L} x_1 <_{\L} 1) = (x_k <_{\L} 1)\ldots (x_2 <_{\L} 1)
(x_1 <_{\L} 1),
\eeq
hence $\E(M)$ is generated (as a monoid) by the chains $m <_{\L} 1$.

Given a set $Y$, we define a $Y$-{\em monoid} to be an ordered pair
of the form $(M,\p)$, where $M$ is a monoid and $\p:Y^* \to M$ is
a surjective morphism. Similarly, we define $Y$-{\em semigroup}.
A {\em morphism} from the $Y$-monoid $(M,\p)$
to the $Y$-monoid $(M',\p')$ is a monoid morphism $\theta:M \to M'$
such that the diagram
$$\xymatrix{
Y^* \ar[rr]^{\p} \ar[dd]_{\p'}&&M \ar[ddll]^{\theta} \\ &&\\
M' &&
}$$
commutes. Whenever possible, to simplify notation, we omit the
morphism in the representation of $Y$-monoids, that is, we view $Y$ as
a subset of $M$ and $\p$ as canonical.

Clearly, $Y$-monoids and their morphisms constitute a category, and we
can consider expansions in the category of $Y$-monoids just as we did
for the category of monoids. We define now the expansion $\E_Y$ in
the category of $Y$-monoids, that can be described as the Rhodes
expansion {\em cut-down to the generators} $Y$.

Indeed, let $M$ be an $Y$-monoid. We remarked before that $\E(M)$
is generated (as a monoid) by the chains $m <_{\L} 1$. We define
$\E_Y(M)$ to be the submonoid of $\E(M)$ generated by the chains $y
<_{\L} 1$ $(y \in Y)$. 
It is shown in \cite{4BR} that $\E_Y$ defines an expansion of
$Y$-monoids. We omit the description of the expansion for morphisms.

The Rhodes expansion has many interesting properties that are
subsequently inherited by $\E_Y$, such as the
following:

\bp
\label{rex}
\cite{20Til, 4BR}
\bi
\item[(i)] The Rhodes expansion is aperiodic.
\item[(ii)] The Rhodes expansion preserves regularity.
\item[(iii)] $\forall \sigma \in \E(M)\; (\sigma \in E(\E(M)) \iff
  \sigma\eta_M \in E(M))$.
\ei
The expansion $\E_Y$ possesses analogous properties.
\ep

We introduce now another expansion with important properties.
Let $M$ be a semigroup and let $M^+$ denote the free semigroup on (the
set) $M$. Hence 
$$M^+ =  \{ (m_1, \ldots, m_k) \mid k \geq 1, \; m_i \in M \}.$$
Given $(m_1, \ldots, m_k) \in M^+$, let
$$F_3(m_1, \ldots, m_k) = \{ (m_1\ldots m_i, m_{i+1}\ldots m_j,
m_{j+1}\ldots m_k) \in M \times M \times M \mid 0 \leq i \leq j \leq k
\}.$$
Write $$\Phi_3(M) = \{ F_3(m_1, \ldots, m_k) \mid (m_1, \ldots, m_k)
\in M^+ \}.$$ We define a multiplication on $\Phi_3(M)$ by
$$F_3(m_1, \ldots, m_k)F_3(m'_1, \ldots, m'_l) = F_3(m_1, \ldots, m_k,
m'_1, \ldots, m'_l).$$
By \cite[Section 7.2]{4BR}, $M \to \Phi_3(M)$ is part of a semigroup expansion
(we omit here the expansion of morphisms). 
The surjective morphisms
$\eta_M: \Phi_3(M) \to M$ are defined by $$(F_3(m_1, \ldots,
m_k))\eta_M = m_1 \ldots
m_k.$$

We can also perform the cut-down to generators for this expansion
\cite{4BR}. Indeed, If $M$ is a $Y$-semigroup, we denote by
$\Phi_{3,Y}(M)$ the subsemigroup of  $\Phi_3(M)$ generated by the
elements of the form
$$F_3(y) = \{ (y,1,1),\; (1,y,1),\; (1,1,y)\} \quad (y \in Y).$$
Then the restriction of the morphism $\eta_M$ to $\Phi_{3,Y}(M)$ is
surjective and $\Phi_{3,Y}(M)$ is part of an expansion of
$Y$-semigroups \cite{4BR}.

We recall that a  semigroup $M$ is said to be {\em finite $\J$-above}
if $\{ y \in  
M\mid y\geq_{\J} x \}$ 
is finite for every $x \in M$. 

The following properties make the expansion $\Phi_3$ of great
interest. Note that, since $\Phi_{3,Y}(M)$ is a subsemigroup of $\Phi_{3}(M)$,
these properties generalize immmediately to $\Phi_{3,Y}$.

\bp
\label{fjap}
\cite[Propositions 7.8 and 7.9]{4BR}
\bi
\item[(i)] $\Phi_3(M)$ is finite $\J$-above for every semigroup $M$;
\item[(ii)] $\Phi_3$ is aperiodic.
\ei
The expansion $\Phi_{3,Y}$ possesses analogous properties.
\ep

\section{The Holonomy Theorem}

Given a semigroup $M$, we denote by $M^I$ the monoid obtained by
adjoining a new identity $I$ to $M$ ({\em even if $M$ is already a
  monoid}), see \cite[Chapter 1]{RS}. We shall consider the Rhodes 
expansion $\E(M^I)$ 
of the monoid $M^I$, consisting of all finite chains
of the form 
$$m_k <_{\L} \ldots <_{\L} m_1 <_{\L} m_0 = I$$
with $k \geq 0$ and $m_i \in M$ $(i = 1, \ldots, k)$.

We say that a mapping $f:M \to N$ between monoids is $\leq_{\J}$-{\em
  preserving} if
\bi
\item[(JP1)] $f(1) = 1$;
\item[(JP2)] $a \leq_{\J} b \Rw f(a) \leq_{\J} f(b)$ for all $a,b \in M$.
\ei
It follows that
$$a \J b \Rw f(a) \J f(b) \hspace{.5cm} \mbox{ for all }a,b \in M.$$

The important particular case arises for mappings $f:M^I \to \N$,
where we consider addition on $\N$. Note
that, for all $n,n' \in \N$,
$$n \leq_J n' \iff n \geq n'.$$
Thus $f:M^I \to \N$ is $\leq_J$-preserving if and only if $f(I) = 0$
and 
$$\forall m,m',m'' \in M^I, \; f(m'mm'') \geq f(m).$$
Since $\N$ is $\J$-trivial, note that
\beq
\label{perdiz}
a \J b \Rw f(a) = f(b) \hspace{.5cm} \mbox{ for all }a,b \in M^I.
\eeq


Let $$\sigma = (m_k <_{\L} \ldots <_{\L} m_1 <_{\L} m_0 = I), \quad
\tau = (m'_l <_{\L} \ldots <_{\L} m'_1 <_{\L} m'_0 = I)$$
be elements of $\E(M^I)$. The {\em maximum $\L$-point of agreement} of
$\sigma$ and $\tau$ is defined by $\sigma \wedge_{\L} \tau = m_r$, with
$$r = \max\{ i \in \{ 0, \ldots, \min\{ k,l \} \} \mid m_0 = m'_0,
\ldots, m_{i-1} = m'_{i-1}, m_i\L m'_i \}.$$

We present now the Holonomy Theorem in its most abstract version:

\bt
\label{hol}
{\em (Holonomy Theorem)} Let $M$ be a semigroup and let $f:M^I \to \N$
be $\leq_J$-preserving. Let $D:
\E(M^I) \times \E(M^I) \to \oo{\N}$ be defined by
$$D(\sigma,\tau) = \left\{
\begin{array}{ll}
f(\sigma \wedge_{\L} \tau)&\mbox{ if } \sigma \neq \tau\\
1 + {\rm sup}f&\mbox{ if } \sigma = \tau.
\end{array}
\right.$$
Then 
\bi
\item[(i)] $D$ is a strict  
length function for $\E(M^I)$;
\item[(ii)]
$D = D_{\chi}$ for some (unique up to isomorphism) strongly faithful elliptic
$\E(M^I)$-tree $\chi$.
\ei 
\et

\proof
We show that $D$ satisfies
axioms (L1) -- (L5).

(L1): Let $\sigma,\tau \in \E(M^I)$. Since $(\sigma \wedge_{\L} \tau) \L
(\tau \wedge_{\L} \sigma)$, we have $D(\sigma,\tau) = D(\tau,\sigma)$
in view of (\ref{perdiz}).

(L2) follows from the definition of $D$.

(L3): Let $\sigma,\tau,\rho \in \E(M^I)$. We start by showing that
\beq
\label{obib}
(\sigma\rho \wedge_{\L} \tau\rho) \leq_{\L} (\sigma \wedge_{\L} \tau).
\eeq

In view of
(\ref{gene}), we may assume that $\rho = (m <_{\L}
I)$.
Write 
$$\sigma = (m_k <_{\L} \ldots <_{\L} m_1 <_{\L} m_0 = I), \quad
\tau = (m'_l <_{\L} \ldots <_{\L} m'_1 <_{\L} m'_0 = I)$$ and assume
that $\sigma \wedge_{\L} \tau = m_r$. Then we may write
$$\sigma\rho = \lm(m_km \leq_{\L} \ldots \leq_{\L} m_1m \leq_{\L} m
<_{\L} I),$$
$$\tau\rho = \lm(m'_lm \leq_{\L} \ldots \leq_{\L} m'_1m \leq_{\L} m <_{\L}
I).$$ 
Clearly, $m_r \L m'_r$ yields $(m_rm) \L (m'_rm)$ and we also have
$m_im = m'_im$ for $ i \in \{ 0, \ldots, r-1\}$. Note that $m_rm$ may
not be in $\sigma\rho$, but some $m_sm \in \L_{m_rm}$ will, and
similarly for $m'_rm$. Hence (\ref{obib}) holds.

Back to checking (L3), we may assume that $\sigma\rho \neq 
\tau\rho$. Hence $\sigma \neq \tau$ as well. 
Since $f$ is $\leq_J$-preserving and 
$$(\sigma\rho \wedge_{\L} \tau\rho) \leq_{\J} (\sigma \wedge_{\L} \tau)$$
by (\ref{obib}), we get
$$D(\sigma\rho,\tau\rho) = f(\sigma\rho \wedge_{\L} \tau\rho) \geq
f(\sigma \wedge_{\L} \tau) = D(\sigma,\tau)$$
since $f$ is $\leq_J$-preserving and $\sigma\rho \neq \tau\rho$.  Thus
(L3) holds. 

(L4): Let $\sigma,\tau,\rho \in \E(M^I)$. We show that
\beq
\label{hol1}
(\sigma \wedge_{\L} \rho) \leq_{\L}
(\sigma \wedge_{\L} \tau) \hspace{.5cm} \vee \hspace{.5cm} (\sigma
\wedge_{\L}\rho)  \leq_{\L} (\tau \wedge_{\L} \rho).
\eeq
Write 
$$\sigma = (m_k <_{\L} \ldots <_{\L} m_1 <_{\L} m_0 = I), \quad
\tau = (m'_l <_{\L} \ldots <_{\L} m'_1 <_{\L} m'_0 = I),$$
$$\rho = (m''_p <_{\L} \ldots <_{\L} m''_1 <_{\L} m''_0 = I),$$
$$(\sigma \wedge_{\L} \rho) = m_r, \quad
(\sigma \wedge_{\L} \tau) = m_s, \quad  (\tau \wedge_{\L} \rho) =
m'_t.$$
Suppose that $m_r \not\leq_{\L} m_s$ and $m_r \not\leq_{\L}
m'_t$. Then $r < s$ and $m''_r \L m_r \not\leq_{\L} m'_t \L m''_t$
yields $r < t$ as well. For $i = 0, \ldots, r$, we have
$m_i = m'_i = m''_i$ since $r < s,t$. Moreover,
$m_{r+1} \L m'_{r+1} \L m''_{r+1}$ since $r+1 \leq s,t$, contradicting
$(\sigma \wedge_{\L} \rho) = m_r$. Therefore $m_r \leq_{\L} m_s$
or $m_r \leq_{\L} 
m'_t$, and so (\ref{hol1}) holds.

To prove (L4), we may assume that $\sigma,\tau,
\rho$ are all distinct. 
Without loss of generality,
  we may assume by (\ref{hol1}) that
$(\sigma \wedge_{\L} \rho) \leq_{\L}
(\sigma \wedge_{\L} \tau)$,
 hence $(\sigma \wedge_{\L}\rho) \leq_{\J}
(\sigma \wedge_{\L}\tau)$
and so $$D(\sigma,\rho) = f(\sigma \wedge_{\L} \rho) \geq
f(\sigma \wedge_{\L} \tau) = D(\sigma,\tau) \geq \min\{ D(\sigma,\tau),
D(\tau,\rho) \}$$ by (JP2). Thus (L4) holds.

(L5): Assume that $D(\tau,\rho) = D(\sigma,\sigma)$ for some 
$\sigma,\tau,\rho \in \E(M^I)$ with $\tau \neq \rho$. It follows that
$$1 + \sup f = D(\sigma,\sigma) = D(\tau,\rho) = f(\tau\wedge_{\L}
\rho) \in \N,$$
a contradiction. Thus $D(\tau,\rho) = D(\sigma,\sigma)$ imples $\tau =
\rho$ and (L5) holds.

Therefore $D$ is a strict 
length function for $\E(M^I)$ and so
$D = D_{\chi}$ for some (unique up to isomorphism) strongly faithful elliptic
$\E(M^I)$-tree $\chi$ by Corollary \ref{flftree}.
\qed

A preordered set $(X, \leq)$ is said to be {\em
  upper finite} if every subset of the form $[x) = \{ y \in X \mid y \geq x
\}$ is finite. In particular, if $(X, \leq)$ is
  upper finite, every nonempty subset of $X$
  must contain a maximal element.

To show how to obtain all $\leq_J$-preserving mappings $f:M^I
\to \N$ when $M$ is finite $\J$-above,
 we introduce the concept of weight
function in a more 
general setting. Given an upper finite partially ordered set
$(P,\leq)$ with 
maximum $I$, a 
{\em weight
function} $w:P \to \N$ is any function satisfying $w(I) = 0$. Given
$w:P \to \N$, let $h_w:P \to N$ be defined by
$$h_w(p) = \max\{ \sum_{i=0}^n w(p_i) \mid p = p_n < \ldots < p_1 <
p_0 = I\mbox{ is a chain in } P \}.$$
Since $(P,\leq)$ is upper finite, $h_w$ is
well defined. 

\bp
\label{dede}
{\em (Dedekind inversion)}. Let $(P,\leq)$ be an upper finite
partially ordered
set $(P,\leq)$ with 
maximum $I$. Then the correspondence $\mu:w \mapsto h_w$ defines a
bijection between the set of all weight functions $w:P \to \N$ and all
order-reversing mappings $h: P \to \N$ satisfying $h(I) = 0$.
\ep

\proof
Let $w:P \to \N$ be a weight function. Assume that $q \leq p$ in
$P$. Since any chain $p = p_n < \ldots < p_1 <
p_0 = I$ in $P$ extends to a chain
$$q \leq p = p_n < \ldots < p_1 <
p_0 = I,$$ we get $h_w(q) \geq h_w(p)$ and so $h_w$ is
order-reversing. Since the only ascending chain starting at $I$ is the
trivial chain and $w(I) = 0$, we have $h_w(I) = 0$. Thus $\mu$ is well
defined. 

Suppose that $w,w':P \to \N$ are distinct weight functions. Take a
maximal element $p$ from the set $$\{ x \in P \mid w(x) \neq w'(x)
\}.$$ Since $P$ is upper finite, there exist such maximal
elements. Assume that $w(p) < w'(p)$. Given a chain $p = p_n < \ldots < p_1 <
p_0 = I$, we have $w(p_i) = w'(p_i)$ for $i = 0, \ldots,n-1$ by
maximality of $p$, hence $$\sum_{i=0}^n w(p_i) = w(p) +
\sum_{i=0}^{n-1} w'(p_i) < \sum_{i=0}^n w'(p_i)$$ and so $h_w(p) <
h_{w'}(p)$. Thus $\mu$ is one-to-one.

Finally, take $h: P \to \N$ order-reversing satisfying $h(I) = 0$. We
define a weight function $w:P \to \N$ as follows. Given $p \in
P\setminus \{ I \}$, let 
$$\oo{p} = \{ q \in P \mid p < q \mbox{ and there exists no } r \in P
\mbox{ such that }p<r <q\}$$
denote the set of all elements of $P$ covering $p$. Since $P$ is
upper finite, $\oo{p}$ is nonempty. 
We define
$$w(p) = h(p) -\max\{ h(q) \mid q \in \oo{p} \}.$$
Since $h$ is order-reversing, $w(p) \geq 0$ and so $w$ is a
well-defined weight function. 
We show that $h = h_w$.

Let $p \in
P\setminus \{ I \}$. We show that
\beq
\label{botan}
h_w(q) = h(q) \mbox{ for every }q \in \oo{p} \Rw h_w(p) = h(p).
\eeq
Indeed, assume the hypothesis and let  $p = p_n < \ldots < p_1 <
p_0 = I$ be a chain in $P$ with
$h_w(p) = \sum_{i=0}^n w(p_i)$. By maximality of $\sum_{i=0}^n
w(p_i)$, we may assume that $p_{n-1} \in \oo{P}$. Moreover, 
$h_w(p_{n-1}) = \sum_{i=0}^{n-1} w(p_i)$ must be maximal among
$\{ h_w(q) \mid q \in \oo{p} \}$. It follows that
$$\begin{array}{lll}
h_w(p)&=&\sum_{i=0}^n w(p_i) = h_w(p_{n-1}) + w(p)\\
&=&\max\{ h_w(q) \mid q \in \oo{p} \} + w(p)\\
&=&\max\{ h(q) \mid q \in \oo{p} \} + w(p)\\
&=&h(p)
\end{array}$$
and so (\ref{botan}) holds. 

Suppose that $h \neq h_w$. Since $P$ is upper finite, we can take a
maximal element $p$ from the set $\{ x \in P \mid h(x) \neq h_w(x)
\}.$ Since $h_w(I) = 0 = h(I)$, we have $p \neq I$. By maximality of
$p$, we must have $h_w(q) = h(q)$ for every $q \in \oo{p}$. But then
$h_w(p) = h(p)$ by (\ref{botan}), a contradiction. Therefore $h_w = h$
and so $\mu$ is onto as required.
\qed

Throughout the paper, we consider the set ${M\slash\J}$ of all
$\J$-classes of a semigroup $M$ partially ordered by
$$\J_a \;\leq \;\J_b \mbox{ if } a \leq_{\J} b.$$

\bc
\label{jpfja}
Let $M$ be a finite $\J$-above semigroup. Then
the $\leq_J$-preserving 
mappings $f:M^I 
\to \N$ are defined by
$$f(m) = h_w(\J_m)$$
for some weight function $w:M^I \slash \J \to \N$.
\ec

\proof
Clearly, (\ref{perdiz}) implies that the $\leq_J$-preserving
mappings $f:M^I 
\to \N$ must be those of the form
$$f(m) = h(\J_m)$$
for some $\leq_J$-preserving
mapping $h:{M^I \slash \J}
\to \N$. Since ${M^I \slash \J}$ is an upper finite partially ordered
set, the claim follows from Proposition \ref{dede}.
\qed

The following is a straighforward corollary from Theorem \ref{hol} and
Corollary \ref{jpfja}.

\bc
\label{holwt}
Let $M$ be a finite $\J$-above semigroup and let $w:M^I \slash \J \to \N$
be a weight function. Let 
$f_w:M^I \to \N$
be defined by
$$f_w(m) = \max\{ \sum_{i=0}^n w(\J_{m_i}) \mid m = m_n <_{\J} \ldots
<_{\J} m_1 <_{\J} 
m_0 = I\mbox{ is a chain in } M^I \}.$$
Let $D:
\E(M^I) \times \E(M^I) \to \oo{\N}$ be defined by
$$D(\sigma,\tau) = \left\{
\begin{array}{ll}
f_w(\sigma \wedge_{\L} \tau)&\mbox{ if } \sigma \neq \tau\\
1+\sup f_w&\mbox{ if } \sigma = \tau.
\end{array}
\right.$$
Then 
\bi
\item[(i)] $D$ is a strict  
length function for $\E(M^I)$;
\item[(ii)]
$D = D_{\chi}$ for some (unique up to isomorphism) strongly faithful elliptic
$\E(M^I)$-tree $\chi$.
\ei 
\ec

\be
\cite[Example 2.8(a)]{Rho1}
Let $M$ be a finite $\J$-above semigroup and let $w:M^I \slash \J
\to \N$ be the null weight function. Then $f_w:M^I \to \oo{\N}$
is the null function 
and so the induced length function $D:\E(M^I)
\times \E(M^I) \to \N$ is induced by the strongly faithful elliptic $M$-tree
$\chi$ whose underlying tree can be depicted by
$$\xymatrix{
&r_0 \ar@{-}[dl] \ar@{-}[d] \ar@{-}[dr] \ar@{-}[drr]&&\\ 
{[1,\sigma_1]}  & {[1,\sigma_2]}  &{[1,\sigma_3]}
& {\ldots} 
}$$
if $\E(M^I) = \{\sigma_1, \sigma_2,
  \sigma_3, \ldots\}$
\ee


\be
Let $M$ be a finite $\J$-above monoid and let $w:M^I \slash \J
\to \N$ be the weight function defined by
$$w(\J_m) = \left\{
\begin{array}{ll}
1&\mbox{ if } m \J 1\\
0&\mbox{ otherwise}.
\end{array}
\right.$$
Then $f_w:M^I \to \oo{\N}$
is defined by
$$f_w(m) = \left\{
\begin{array}{ll}
0&\mbox{ if } m = I\\
1&\mbox{ otherwise}
\end{array}
\right.$$
and so the induced length function $D:\E(M^I)
\times \E(M^I) \to \oo{\N}$ is induced by the elliptic $M$-tree
$\chi$ whose underlying tree can be depicted by
$$
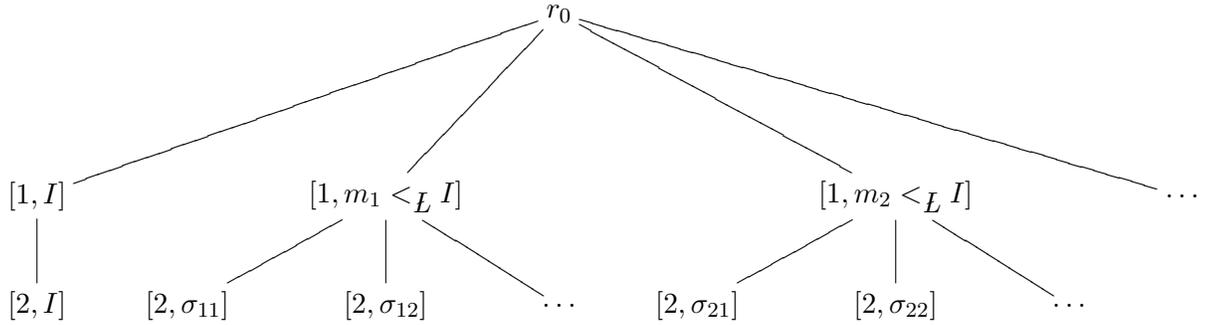
\xymatrix{
&&&r_0 \ar@{-}[ddlll] \ar@{-}[ddl] \ar@{-}[ddrr] \ar@{-}[ddrrrr] &&&&\\
&&&&&&&\\
{[1,I]} \ar@{-}[d] &&{[1,m_1  <_{\L} I]} \ar@{-}[dl] \ar@{-}[d]
\ar@{-}[dr] &&&
{[1,m_2  <_{\L} I]} \ar@{-}[dl] \ar@{-}[d] \ar@{-}[dr] &&{\ldots}\\
{[2,I]} &{[2,\sigma_{11}]} &{[2,\sigma_{12}]}
&{\ldots}&{[2,\sigma_{21}]}
&{[2,\sigma_{22}]}&{\ldots}&  
}$$
where $M\slash \L \; = \; \{ \L_{m_1}, \L_{m_2}, \ldots \} $ and
$\{ \sigma_{i1}, \sigma_{i2}, \ldots \}$ denotes the set of all $(
\ldots <_{\L} m'_i <_{\L} I) \in \E(M^I)$ with $m'_i \L m_i$. 
\ee

\proof
It is immediate that $f_w$ must be of the claimed form since all
$m \in M^I$ except $I$ satisfy $m \leq_{\J} 1$. Hence, for 
$$\sigma = (m_k <_{\L} \ldots <_{\L} m_1 <_{\L} m_0 = I), \quad
\tau = (m'_l <_{\L} \ldots <_{\L} m'_1 <_{\L} m'_0 = I),$$
we have
$$D(\sigma,\tau) = \left\{
\begin{array}{ll}
2&\mbox{ if }\sigma = \tau\\
1&\mbox{ if }\sigma \neq \tau \mbox{ and }k,l > 0\mbox{ and }m_1\L m'_1\\
0&\mbox{ otherwise}
\end{array}
\right.$$
since $D(\sigma,\tau) = 1$ if and only if $\sigma \neq \tau$ and $\sigma
\wedge_{\L} \tau \neq I$.

Following the Chiswell construction in the proof of Theorem
\ref{lftree}, the underlying tree $T$ of the elliptic $M$-tree $\chi$
induced by $\chi$ has vertex set $$\vert(T) = \{ r_0\} \cup \{
[i,\sigma];\; i = 1,2; \; \sigma \in \E(M^I) \}.$$ 
We have $[1,\sigma] = [1,\tau]$ if and only if $D(\sigma,\tau) \geq 1$,
hence $$\{ [1,I], [1,m_1  <_{\L} I], [1,m_2  <_{\L} I], \ldots \}$$
constitutes a full set of representatives for the classes
$[1,\sigma]$. Clearly, $[2,\sigma] = [2,\tau]$ if and only if $\sigma
= \tau$ and so the tree is the claimed one.
\qed

Given a finite $\J$-above
semigroup $M$, we define a mapping $h_{\J}: M \to \N$
by
$$h_{\J}(m) = \max\{ k \in \N : \mbox{ there exists a chain }
m = m_0 <_{\J} \ldots <_{\J} m_k \mbox{ in } M \}.$$
We say that $h_{\J}$ is the
(Dedekind) 
{\em $\J$-height function} of $M$ \cite{17Bir}. Since the $\J$-class of
$I$ contains only $I$ and lies above all the others, it is immediate that
$M^I$ has also a Dedekind 
  $\J$-height function $h'_{\J}$, satisfying
$$h'_{\J}(m) = \left\{
\begin{array}{ll}
h_{\J}(m)+1&\mbox{ if } m \in M\\
0&\mbox{ if } m = I
\end{array}
\right.$$

\bp
\label{hj}
Let $M$ be a finite $\J$-above semigroup and let $h_{\J}$ be the
$\J$-height function of 
$M^I$. Then:
\bi
\item[(i)] $h_{\J}$ is $\leq_{\J}$-preserving;
\item[(ii)] $h_{\J} = f_w$ for the
  weight function $w: M^I \slash \J
\to \N$ defined by 
$$w(\J_m) = \left\{
\begin{array}{ll}
1&\mbox{ if } m \in M\\
0&\mbox{ if } m = I.
\end{array}
\right.$$ 
\ei
\ep

\proof
(i) Immediate.

(ii) We have $h_{\J}(I) = 0 = h_w(I)$. Given $m \in M$,
$$\begin{array}{lll}
h_{\J}(m)&=&\max\{ k \in \N : \mbox{ there exists a chain }
m = m_0 <_{\J} \ldots <_{\J} m_k \mbox{ in } M^I \}\\
&=&\max\{ k \in \N : \mbox{ there exists a chain }
m = m_0 <_{\J} \ldots <_{\J} m_k = I\mbox{ in } M^I \}\\
&=&\max\{ \sum_{i=0}^k w(\J_{m_i}) : \mbox{ there exists a chain }
m = m_0 <_{\J} \ldots <_{\J} m_k = I\mbox{ in } M^I \}\\
&=&f_w(m).
\end{array}
$$
\qed

The mapping $h_{\J}$ will play the most important role as a
$\leq_{\J}$-preserving mapping in forthcoming sections.

We can use the expansion $\Phi_3$ to avoid the finite
$\J$-above requirement in Corollary \ref{holwt}:

\bc
\label{F3hol}
Let $M$ be a semigroup and let $w:(\Phi_3(M))^I \slash \J \to \N$
be a weight function. Let
$f_w:(\Phi_3(M))^I \to \N$
be defined by
$$f_w(x) = \max\{ \sum_{i=0}^n w(\J_{x_i}) \mid x = x_n <_{\J} \ldots
<_{\J} x_1 <_{\J} 
x_0 = I\mbox{ is a chain in } (\Phi_3(M))^I \}.$$
Let $D:
\E((\Phi_3(M))^I) \times \E((\Phi_3(M))^I) \to \oo{\N}$ be defined by
$$D(\sigma,\tau) = \left\{
\begin{array}{ll}
f_w(\sigma \wedge_{\L} \tau)&\mbox{ if } \sigma \neq \tau\\
1+\sup f_w&\mbox{ if } \sigma = \tau.
\end{array}
\right.$$
Then 
\bi
\item[(i)]
$D$ is a strict 
length function for $\E((\Phi_3(M))^I)$;
\item[(ii)]
$D = D_{\chi}$ for some (unique up to isomorphism) strongly faithful elliptic
$\E((\Phi_3(M))^I)$-tree $\chi$;
\item[(iii)] the canonical surjective morphism $\E((\Phi_3(M))^I) \to
  M$ is aperiodic.
\ei 
\ec

\proof
(i) and (ii) follow from Corollary \ref{holwt} since $\Phi_3(M)$ is
finite $\J$-above by Proposition \ref{fjap}(i).

For (iii), we can decompose the canonical morphism $\E((\Phi_3(M))^I) \to
  M$ as the composition
$$\E((\Phi_3(M))^I) \vvvvlongmapright{\eta_{(\Phi_3(M))^I}}
(\Phi_3(M))^I
\mapright{\p} \Phi_3(M) \longmapright{\eta'_M} M.$$
The morphisms $\eta_{(\Phi_3(M))^I}$ and $\eta'_M$ are aperiodic by
Propositions \ref{rex}(i) and \ref{fjap}(ii). Since $\p$ is trivially
aperiodic and the composition of aperiodic morphisms is aperiodic, the
result follows. 
\qed

\section{Stable monoids and the Zeiger encoding}

We start by introducing some well-known concepts and results. For
details, the reader is referred to \cite{CP,RS}.

A semigroup $M$ is said to be {\em stable} if the following
conditions hold for all $a,x \in M$:
\bi
\item[(S1)] $ax\,\J\, a \Rw ax\,\R\, a$;
\item[(S2)] $xa\,\J\, a \Rw xa\,\L\, a$.
\ei
It folows easily that \cite{CP,RS}
\beq
\label{stbl}
\mbox{if $M$ is stable, then } <_{\R} \;\subseteq\; <_{J} \mbox{ and }
<_{\L} \;\subseteq\; <_{J}.
\eeq 
The following lemma will turn out to be quite useful:

\bl
\label{gedi}
Let $M$ be stable and let $a,b,c \in M$ satisfy $a <_{\L} b \R
bc$. Then $a \R ac <_{\L} bc$.
\el

\proof
Clearly, $a <_{\L} b$ yields $ac \leq_{\L} bc$. Since 
$$ac \leq_{\J} a <_{\L} b \R bc,$$
it follows from (\ref{stbl}) that $ac <_{\J} bc$ and so $ac <_{\L}
bc$.

On the other hand, $b \R bc$ yields $b = bcx$ for some $x \in
M$. Since $a <_{\L} b$, we get $a = acx$ and so $a \R ac$.
\qed

Every finite $\J$-above semigroup is stable, a fact that will be
thoroughly used throughout the paper.

Assume that $M$ is stable. Then the Green relations $\J$ and $\D$
on $M$
coincide. Given a $\J$-class $J$ of a monoid, we can always define
a semigroup structure $(J^0,*)$ on $J^0 = J \cup \{ 0 \}$ by taking
$$a*b = \left\{
\begin{array}{ll}
ab&\mbox{ if } a,b,ab \in J\\
0&\mbox{ otherwise.}
\end{array}
\right.$$
If the monoid is stable, the semigroup $J^0$ defined above is
completely 0-simple and can thus be given a {\em Rees matrix
coordinatization}: there exist nonempty sets $A,B$, a group $G$ and a
$(B\times A)$-matrix $C$ with entries in $G \cup \{ 0 \}$ such that
$J^0 \cong M^0(G,A,B,C)$, where $M^0(G,A,B,C) = (A \times G \times B)
\cup \{ 0 \}$ is the semigroup with zero defined by
$$(a,g,b)(a',g',b') = \left\{
\begin{array}{ll}
(a,gC(b,a')g',b')&\mbox{ if } C(b,a') \in G\\
0&\mbox{ if }C(b,a') = 0.
\end{array}
\right.$$
The Green relations in $M^0(G,A,B,C)$ are characterized
by
$$\begin{array}{l}
(a,g,b) \,\R\, (a',g',b') \iff a = a',\\
(a,g,b) \,\L\, (a',g',b') \iff b = b'.
\end{array}$$

We shall need the detailed construction of the
Rees matrix semigroup, so we present it briefly. For more details, see
\cite{CP,RS}.

We fix a $\H$-class $H$ in $J$ and $h_0 \in H$. Let $A$ (respectively
$B$) be the set of $\R$-classes (respectively
 $\L$-classes) in $J$. For every $a \in A$, we fix $\wh{a} \in a \; \cap
 \L_{h_0}$. For every $b \in B$, we fix also $\wh{b} \in b \; \cap
 \R_{h_0}$. Finally, we fix $e_a,\oo{e}_a, f_b, \oo{f}_b 
 \in M$ such that
$$e_a\wh{a} = h_0, \quad \oo{e}_ah_0 = \wh{a}, \quad \wh{b}f_b = h_0,
\quad h_0\oo{f}_b = \wh{b}.$$
Let $\stab(H) = \{ x \in M \mid Hx = H\}$ and define an equivalence
relation on $\stab(H)$ by
$$[x] = [y] \quad \mbox{if } h_0x = h_0y.$$
Then the quotient $$G = \{ [x] \mid x \in \stab(H) \}$$ is the {\em
  Sch\"{u}tzenberger group} of $H$. For each $h \in H$, we fix some
$\wt{h} \in M$ such that $h_0\wt{h} = h$. By the well-known
Green's Lemma \cite{CP,RS}, $\wt{h} \in \stab(H)$. Then there exists some 
$(B\times A)$-matrix $C$ with entries in $G \cup \{ 0 \}$ such that
$$\begin{array}{rcl}
J^0&\to&M^0(G,A,B,C)\\
u&\mapsto&\left\{
\begin{array}{ll}
(\R_u = a, [\wt{e_auf_b}], \L_u = b)\;&\mbox{ if }u\in J\\
0\;&\mbox{ if }u = 0
\end{array}
\right.
\end{array}$$
is a semigroup isomorphism.

\medskip

Throughout this section, we assume that $M$ is a fixed stable
semigroup. Hence $M^I$ is a stable monoid. We fix a coordinatization
(assuming equality to simplify notation) 
$$\J_{m}^0 =
M^0(G_{m},A_{m},B_{m},C_{m})$$
for every $m \in M^I$. We assume that 1 denotes the identity
in every group and $1 \in A_m,B_{m}$ for every $m \in M^I$. If $m\J m'$,
we assume of course that $(G_{m},A_{m},B_{m},C_{m}) =
(G_{m'},A_{m'},B_{m'},C_{m'})$.

We fix mappings
$$\begin{array}{rclrcl}
M^I&\to&M^I\hspace{2cm}&M^I&\to&M^I\\
m&\mapsto&m^*&m&\mapsto&m^{\#}
\end{array}$$
defined as follows. If $m = (a,g,b) \in \; \J_{m}^0$, $m^*,m^{\#}$
satisfy
$$mm^* = (a,1,1), \quad (a,1,1)m^{\#} = m.$$
The existence of such elements follows from $(a,g,b) \,\R\, (a,1,1)$.
Note that $I^* = I^{\#} = I$.

\bl
\label{propss}
For all $m,m'\in M$,
\bi
\item[(i)] $mm^*m^{\#} = m$;
\item[(ii)] $m \R m' \iff mm^* = m'{m'}^*$.
\ei
\el

\proof
(i) is trivial. Assume now that $m \R m'$. Then we may write $m = (a,g,b),
\; m' = (a,g',b')$ as elements of the same $\J$-class. Thus
$mm^* = (a,1,1) = m'{m'}^*$. Conversely, assume that $mm^* =
m'{m'}^*$. By (i), we get $m = mm^*m^{\#} = m'{m'}^*m^{\#}$ and so $m
\leq_{\R} m'$. By symmetry, it follows that $m \R m'$ and so (ii)
holds.
\qed

We define $F_{\J}(M^I)$ to be the set of all finite chains of the form
$$n_k <_{\J} \ldots <_{\J} n_1 <_{\J} n_0 = I$$
with $k \geq 0$ and $n_i \in M^I$.

\bl
\label{eej}
Given $(m_k <_{\L} \ldots <_{\L} m_2 <_{\L} m_1 <_{\L} m_0 = I) \in
\E(M^I)$, let 
$$x_0 = I, \quad x_i = m_im_{i-1}^* \; (i
= 1, \ldots,k).$$
Then:
\bi
\item[(i)] $x_i\R m_i$ for $i = 0, \ldots, k$;
\item[(ii)] $x_i <_{\L} x_{i-1}x_{i-1}^*$ for $i = 1, \ldots, k$;
\item[(iii)] $x_i <_{\J} x_{i-1}$ for $i = 1, \ldots, k$;
\item[(iv)] $m_i = x_im_{i-1}^{\#}$ for $i = 1, \ldots, k$.
\ei
\el

\proof
We have $x_0 = I = m_0$. Let $i \in \{1, \ldots,k\}$. Since $m_{i}
<_{\L} m_{i-1}$, 
we may write 
$m_{i} = ym_{i-1}$ for some $y \in M^I$. Now $m_{i-1}\,\R\,
m_{i-1}m_{i-1}^*$ yields $$x_i = m_im_{i-1}^* = ym_{i-1}m_{i-1}^* \,\R\,
ym_{i-1} = m_i.$$ Thus (i) holds.

Since $m_{i}
<_{\L} m_{i-1}$, 
we have $m_{i}m_{i-1}^*
\leq_{\L} m_{i-1}m_{i-1}^*$. Since 
$$m_{i}m_{i-1}^* \leq_{\J} m_i <_{\J} m_{i-1} \,\R\,
m_{i-1}m_{i-1}^*$$ by (\ref{stbl}),
we obtain $m_{i}m_{i-1}^*
<_{\L} m_{i-1}m_{i-1}^*$ and so $$x_i = m_{i}m_{i-1}^*
<_{\L} m_{i-1}m_{i-1}^* = x_{i-1}x_{i-1}^*$$ by (i) and Lemma
\ref{propss}(ii). Thus (ii) holds.

Since $x_{i-1}x_{i-1}^* \,\R\, x_{i-1}$, (ii) implies (iii) in view of
(\ref{stbl}). Finally, $m_{i} = ym_{i-1}$ yields
$$m_i = ym_{i-1} = ym_{i-1}m_{i-1}^*m_{i-1}^{\#} =
m_{i}m_{i-1}^*m_{i-1}^{\#} = x_im_{i-1}^{\#}$$
and (iv) holds as well.
\qed

Lemmas \ref{propss} and \ref{eej} will be used so thoroughly for the
remainder of the paper that we shall often omit a specific reference to
them.

We can now define a mapping 
$\epsilon:\E(M^I) \to F_{\J}(M^I)$ by
$$\epsilon(m_k <_{\L} \ldots <_{\L} m_1 <_{\L} m_0 = I) =
(x_k <_{\J} \ldots <_{\J} x_1
<_{\J} x_0 = I)$$
taking $x_0 = I$ and $x_i = m_im_{i-1}^*$ for $i
= 1, \ldots,k$. 

Note that $\epsilon$ is sequential in the sense that if $$\epsilon(m_k
<_{\L} \ldots <_{\L} m_1 <_{\L} m_0 = I) = 
(x_k <_{\J} \ldots <_{\J} x_1
<_{\J} x_0 = I)$$ and $k > 0$, then $$\epsilon(m_{k-1} <_{\L} \ldots
<_{\L} m_1 <_{\L} m_0 = I) = 
(x_{k-1} <_{\J} \ldots <_{\J} x_1
<_{\J} x_0 = I).$$ The mapping $\epsilon$ is known as the
{\em Zeiger encoding map} and plays an essential role in the
next section to ensure the Zeiger property of the wreath product.

\bp
\label{inje}
The mapping $\epsilon:\E(M^I) \to F_{\J}(M^I)$ is one-to-one.
\ep

\proof
By definition, $\epsilon$ preserves chain length. Take $$\sigma = (m_k
<_{\L} \ldots <_{\L} m_1 <_{\L} m_0 = I), \quad \tau = (m'_k
<_{\L} \ldots <_{\L} m'_1 <_{\L} m'_0 = I)$$
such that $$\epsilon(\sigma) = (x_k <_{\J} \ldots <_{\J} x_1
<_{\J} x_0 = I) = \epsilon(\tau).$$
We show that $m_i = m'_i$ for $i = 0, \ldots,k$ by induction on
$i$. The case $i = 0$ being trivial, assume that $i \in \{
1,\ldots,k\}$ and $m_{i-1} = m'_{i-1}$. By Lemma \ref{eej}(ii) and
the induction hypothesis, we get $$m_i = x_im_{i-1}^{\#} =
x_i(m'_{i-1})^{\#} = m'_i.$$
It follows that $\sigma = \tau$ and so $\epsilon$ is
one-to-one.
\qed

The next result will reveal in the next section the adequacy of the
encoding map $\epsilon$ 
to deal with the product in $\E(M^I)$.

\bt
\label{epz}
Let $\sigma,\tau \in \E(M^I)$ with
$$\sigma = (m_k
<_{\L} \ldots <_{\L} m_1 <_{\L} m_0 = I), \quad \sigma\tau = (m'_l
<_{\L} \ldots <_{\L} m'_1 <_{\L} m'_0 = I),$$
$$\epsilon(\sigma) = (x_k <_{\J} \ldots <_{\J} x_1
<_{\J} x_0 = I),\quad
\epsilon(\sigma\tau) = (x'_l <_{\J} \ldots <_{\J} x'_1
<_{\J} x'_0 = I).$$
Assume that $m_k \R m'_l$. If 
$$\sigma' = (m_{k+p}
<_{\L} \ldots <_{\L} m_1 <_{\L} m_0 = I),\quad  
\epsilon(\sigma') = (x_{k+p} <_{\J} \ldots <_{\J} x_1
<_{\J} x_0 = I),$$
then $$\epsilon(\sigma'\tau) = (x_{k+p} <_{\J} \ldots <_{\J} x_{k+1}
<_{\J} x'_l <_{\J} \ldots <_{\J} x'_1 
<_{\J} x'_0 = I).$$
\et

\proof
Let $y \in M^I$ denote the leftmost term
in $\tau$. Then $m'_l = m_ky$. Since $m'_l \R m_k$, we can apply
successively Lemma \ref{gedi} to get 
$$\sigma'\tau = (m_{k+p}y <_{\L} \ldots <_{\L} m_{k+1}y <_{\L} m'_l
<_{\L} \ldots <_{\L} m'_0)$$
and
\beq
\label{gedi1}
m_{k+i}y \R m_{k+i} \quad (i = 0,\ldots, p).
\eeq
Since
$\epsilon$ is sequential, we obtain
$$\epsilon(\sigma'\tau) = (m_{k+p}y(m_{k+p-1}y)^* <_{\J} \ldots
<_{\J} m_{k+1}y(m_ky)^* 
<_{\J} x'_l <_{\J} \ldots <_{\J} x'_0).$$
Now (\ref{gedi1}) and Proposition \ref{propss}(ii) yield
$$m_{k+i-1}y(m_{k+i-1}y)^* = m_{k+i-1}m_{k+i-1}^* \quad (i = 1,\ldots,
p).$$
Thus $m_{k+i} <_{\L} m_{k+i-1}$ yields
$$m_{k+i}y(m_{k+i-1}y)^* = m_{k+i}m_{k+i-1}^* = x_{k+i}$$
for $i = 1,\ldots,
p$ and the lemma is proved.
\qed

We complete the section with a 
straightforward consequence of 
Green's Lemma \cite{CP,RS}, to be used in the next section.

\bp
\label{schutz}
Let $J$ be a $\J$-class of a stable monoid $M$ with $J^0 =
M^0(G,A,B,C)$. Let $u = (a,g,b) \in J$ and $v \in M$ be such that $uv
\R u$ and let
$$\begin{array}{rcl}
\p:G&\to&G\\
h&\mapsto&((a,h,b)v)\pi_2.
\end{array}$$
Then $\p$ is bijective and there exists some $g_0 \in G$ such that
$h\p = hg_0$ for every $h 
\in G$.
\ep

\proof
By Green's Lemma, we have a bijection
$$\begin{array}{rcl}
\H_u&\to&\H_{uv}\\
u'&\mapsto&u'v,
\end{array}$$
hence $\p$ is well defined.
Let $H$ be the fixed $\H$-class in the construction $M^0(G,A,B,C)$ and
consider all the distinguished elements introduced there. Our Rees
matrix representation restricts to bijections
$$\begin{array}{rclrcl}
\H_u&\to&G\hspace{1.5cm}&\H_{uv}&\to&G\\
u'&\mapsto&[\wt{e_au'f_b}]&u'v&\mapsto&[\wt{e_au'vf_{b'}}]
\end{array}$$
where $a = \; \R_{u} \; = \; \R_{uv}$, $b = \; \L_{u}$ and $b' = \; \L_{uv}$.
Hence we obtain a diagram
$$\xymatrix{
\H_u \ar[ddd] \ar[rrrr] &&&&\H_{uv} \ar[ddd]\\
&u' \ar[dd] \ar[rr]&&u'v \ar[dd]&\\
&&&&\\
G&[\wt{e_au'f_b}]&&[\wt{e_au'vf_{b'}}]&G
}$$
We must show that the mapping
$$\p:[\wt{e_au'f_b}] \mapsto [\wt{e_au'vf_{b'}}],$$
a composition of bijections, 
can be defined by right multiplication. 
We show that $h_0\oo{f}_bvf_{b'} \in H$ and
\beq
\label{schutz1}
[\wt{e_au'vf_{b'}}] = [\wt{e_au'f_{b}}][\wt{h_0\oo{f}_bvf_{b'}}]
\eeq
for every $u' \in \; \H_u$. Indeed, take $u_0 \in \; \H_u$ such that
$e_au_0f_b = h_0$. Since $e_au_0 =
e_au_0f_b\oo{f}_b$, we have
$$h_0\oo{f}_bvf_{b'} =  e_au_0f_b\oo{f}_bvf_{b'} = e_au_0vf_{b'} \in H.$$
Moreover, we have
\beq
\label{schutz2}
\forall h \in H\; \forall w \in M\; (hw \in H \Rw [\wt{hw}] =
[\wt{h}][\wt{h_0w}]). 
\eeq
Indeed, we have $h_0w \in H$ by Green's Lemma. Since $h_0\wt{h_0w} =
h_0w$, we obtain $h\wt{h_0w} =
hw$ and so (\ref{schutz2})
follows from  $h_0\wt{h}\wt{h_0w} = h\wt{h_0w} = hw = h_0\wt{hw}$.

Finally, making $h = e_au'f_b$ and $w = \oo{f}_bvf_{b'}$ in
(\ref{schutz2}), we get
$$[\wt{e_au'vf_{b'}}] = [\wt{e_au'f_b\cdot \oo{f}_bvf_{b'}}] = 
[\wt{e_au'f_{b}}][\wt{h_0\oo{f}_bvf_{b'}}]$$
and so (\ref{schutz2}) holds as claimed.
\qed

\section{From elliptic $M$-trees to wreath products}

Let $(r_0,T)$ be a rooted tree
and let $X$ be a nonempty
set. Given $v \in \vert(T)$, let $\son(v)$ denote the
(possibly empty)
set of sons of $v$. A mapping $f: \vert(T) \backslash \{ r_0 \} \to X$
is said to be 
{\em locally injective} if, for every $v \in \vert(T)$, $f|_{\son(v)}$
is injective.  

For $i = 1,2, \ldots$, let $\vert_i(T)$ denote the set of all $v \in
\vert(T)$ having depth $i$.

\bt
\label{embed}
Let $M$ be a semigroup and assume that $\theta:M^I \to \ell(r_0,T)$ is
a faithful elliptic action of $M^I$ on a uniform rooted tree
$(r_0,T)$. Let $f: (\vert(T))\setminus \{ r_0\}  \to \cup_{i\geq 1}
X_i$ be locally 
injective with $f(\vert_i(T)) \subseteq X_i$ for $i
\geq 1$. Then:
\bi
\item[(i)] if $(r_0,T)$ has finite depth $l$, then $M^I$ embeds in the
  wreath product $$(X_l,P(X_l))\circ \ldots \circ (X_2,P(X_2)) \circ
  (X_1,P(X_1));$$
\item[(ii)] if $(r_0,T)$ has infinite depth, then $M^I$ embeds in the
  infinite wreath product $$\ldots \circ (X_3,P(X_3)) \circ (X_2,P(X_2)) \circ
  (X_1,P(X_1)).$$
\ei
\et

\proof
We prove the finite depth case, the infinite case being analogous. 

Write $X = \cup_{i = 1}^l (X_i \times
\ldots \times 
X_1)$. Associating vertices with rays as usual,
we define a mapping 
$$\begin{array}{rcl} 
\psi:\ray(r_0,T)\setminus \{ (r_0)\} &\to&X\\
(v_i, \ldots, v_1,r_0)&\mapsto&(f(v_i), \ldots, f(v_1)).
\end{array}$$ 
Suppose that $(v_i, \ldots, v_1,v_0 = r_0),(v'_i, \ldots,
v'_1,v_0 = r_0) \in 
\ray(r_0,T)$ are distinct. Let 
$$k = \min\{ j \in \{ 1, \ldots,i\} : v_j \neq v'_j\}.$$
By minimality of $k$, $v_k$ and $v'_k$ must have $v_{k-1} = v'_{k-1}$
as their common father. 
Since $f$ is locally injective, it follows that $f(v_k) \neq f(v'_k)$,
thus $(v_i, \ldots, v_1,r_0)\psi \neq (v'_i, \ldots, v'_1,r_0)\psi$ and so
$\psi$ is one-to-one. 

Let $$\begin{array}{rcl} 
\Psi:M^I&\to&(X_l,P(X_l))\circ \ldots \circ   (X_1,P(X_1))\\
m&\mapsto&\Psi_m
\end{array}$$ 
be
defined by
$$x\Psi_m = \left\{
\begin{array}{ll}
x\psi\inv\theta_m\psi&\mbox{ if } m \in M\hspace{1.5cm} (x \in X).\\
x&\mbox{ if } m = I
\end{array}
\right.$$

Clearly, $\Psi_m \in P(X)$. 
To show that $\Psi_m \in (X_l,P(X_l))\circ \ldots \circ
(X_1,P(X_1))$, we only need to check that $\Psi_m$ 
is sequential. We may assume that $m \in M$. 
Since $\dom\Psi_m = \im\psi$, (SQ1) holds. Since $\theta$ is an
elliptic action, (SQ2) holds as well. 

Let $(x_j,\ldots,
x_1), (x'_k,\ldots, 
x'_1) \in \dom\Psi_m = \im\psi$ and suppose that $(x_j,\ldots,
x_1) \equiv_i (x'_k,\ldots, 
x'_1)$ with $1 \leq i \leq j,k$.  Write
$$\begin{array}{l}
(x_j,\ldots,
x_1) = (v_j,\ldots,
v_1,r_0)\psi = (f(v_j),\ldots, f(v_1))\\
(x'_k,\ldots,
x'_1) = (v'_k,\ldots,
v'_1,r_0)\psi = (f(v'_k),\ldots, f(v'_1)).
\end{array}$$
Since $(x_j,\ldots,
x_1) \equiv_i (x'_k,\ldots, 
x'_1)$, we have $f(v_i) = f(v'_i), \ldots, f(v_1) = f(v'_1)$. Since
$f$ is locally injective, we obtain successively $v_1 = v'_1, \ldots,
v_i = v'_i$.  Thus $(v_j,\ldots,
v_1,r_0)\theta_m \equiv_{i+1} (v'_k,\ldots,
v'_1,r_0)\theta_m$ and so $(v_j,\ldots,
v_1,r_0)\theta_m\psi \equiv_{i} (v'_k,\ldots,
v'_1,r_0)\theta_m\psi$, that is,
$$(x_j,\ldots,
x_1)\Psi_m = (x_j,\ldots,
x_1)\psi\inv\theta_m\psi \equiv_{i} (x'_k,\ldots,
x'_1)\psi\inv\theta_m\psi = (x'_k,\ldots,
x'_1)\Psi_m.$$
Thus (SQ3) holds. Therefore $\Psi_m$
is sequential and so $\Psi_m \in (X_l,P(X_l))\circ \ldots \circ
(X_1,P(X_1))$.

We show next that $\Psi$ is a monoid homomorphism. It suffices to show
that $\Psi_{mm'} = \Psi_{m}\Psi_{m'}$ for all $m,m' \in M$. Since
$\theta$ is an action and $\psi$ is injective, we obtain
$$\Psi_{mm'} = \psi\inv\theta_{mm'}\psi =
\psi\inv\theta_{m}\theta_{m'}\psi =
\psi\inv\theta_{m}\psi\psi\inv\theta_{m'}\psi = \Psi_{m}\Psi_{m'}.$$
Thus $\Psi$ is a monoid homomorphism.

It remains to show that $\Psi$ is one-to-one. Let $m,m' \in M^I$. We show that
\beq
\label{embed1}
m \neq m' \Rw \psi\inv\theta_{m}\psi \neq \psi\inv\theta_{m'}\psi.
\eeq

Indeed, since $\theta$ is one-to-one, $m \neq m'$ implies that there
exists some $v \in \vert(T)$ such that $v\theta_m \neq
v\theta_{m'}$. Taking the geodesic $(v = v_i,\ldots, 
v_1,r_0) \in \ray(r_0,T)$, it follows that $(v_i,\ldots,
v_1,r_0)\theta_m \neq (v_i,\ldots,
v_1,r_0)\theta_{m'}$. Let $$(x_i,\ldots,
x_1) = (v_i,\ldots,
v_1,r_0)\psi.$$
Then $$(x_i,\ldots,
x_1)\psi\inv\theta_m = (v_i,\ldots,
v_1,r_0)\theta_m \neq (v_i,\ldots,
v_1,r_0)\theta_{m'} = (x_i,\ldots,
x_1)\psi\inv\theta_{m'}.$$
Since $\psi$ is one-to-one, we obtain $(x_i,\ldots,
x_1)\psi\inv\theta_m\psi \neq (x_i,\ldots,
x_1)\psi\inv\theta_{m'}\psi$ and so (\ref{embed1}) holds.

Assume that $m \neq m'$. Now (\ref{embed1}) implies that 
$\Psi_m \neq \Psi_{m'}$ if $m,m' \in M$. For the remaining cases, we
may assume that $m = I$. If $\psi$ is
onto, then $\Psi_I = \psi\inv\theta_{I}\psi$ and so (\ref{embed1})
also implies that 
$\Psi_I \neq \Psi_{m'}$. Otherwise, we have $\Psi_I \neq \Psi_{m'}$
since $\dom\Psi_I = X \supset \im\psi = \dom\Psi_{m'}$. Therefore
$\Psi$ is one-to-one.
\qed 

By Theorem \ref{embed}, we know that
when a monoid $M^I$ acts faithfully by elliptic contractions on a uniform
rooted tree, then $M^I$ embeds into a (possibly infinite) wreath
product  $\ldots \circ (X_2,M_2) \circ
  (X_1,M_1)$ of partial transformation monoids. The question is {\em
    how small can the $M_i$'s be made} (where {\em small} is used in the
  sense of division). We start with a series of lemmas.

\bl
\label{kkk}
Let $\sigma,\tau,\rho \in \E(M^I)$. 
\bi
\item[(i)] If $h_{\J}(\sigma \wedge_{\L} 
\tau) = h_{\J}(\sigma \wedge_{\L} 
\rho)$, then $(\sigma \wedge_{\L} 
\tau) = (\sigma \wedge_{\L} 
\rho)$.
\item[(ii)] $(\rho\sigma \wedge_{\L} 
\rho\tau) \leq_{\L} (\sigma \wedge_{\L} 
\tau)$.
\ei
\el

\proof
(i) Assume that $h_{\J}(\sigma \wedge_{\L} 
\tau) = h_{\J}(\sigma \wedge_{\L} 
\rho)$. If $(\sigma \wedge_{\L} 
\tau) \neq (\sigma \wedge_{\L} 
\rho)$, we may assume that $(\sigma \wedge_{\L} 
\tau) <_{\L} (\sigma \wedge_{\L} 
\rho)$ and so  $(\sigma \wedge_{\L} 
\tau) <_{\J} (\sigma \wedge_{\L} 
\rho)$ by (\ref{stbl}), contradicting $h_{\J}(\sigma \wedge_{\L} 
\tau) = h_{\J}(\sigma \wedge_{\L} 
\rho)$. Therefore $(\sigma \wedge_{\L} 
\tau) = (\sigma \wedge_{\L} 
\rho)$.

(ii) Write 
$$\sigma = (m_k <_{\L} \ldots <_{\L} m_0 = I), \quad
\tau = (m'_l <_{\L} \ldots <_{\L} m'_0 = I).$$
In view of (\ref{gene}), we may assume that $\rho = (n <_{\L}
I)$. Hence
$$\rho\sigma = \lm(nm_k \leq_{\L} m_k <_{\L} \ldots <_{\L} m_0),$$
$$\rho\tau = \lm(nm'_l \leq_{\L} m'_l <_{\L} \ldots <_{\L} m'_0).$$
The claim follows at once.
\qed

Next we define 
$$\begin{array}{ll}
V(M^I) = \{&(\sigma,\tau) \in  \E(M^I) \times \E(M^I) \mid
\forall \rho \in  \E(M^I)\\
&(\rho\sigma \wedge_{\L} 
\rho\tau) \L (\sigma \wedge_{\L} 
\tau) \; \Rw \; (\rho\sigma \wedge_{\L} 
\rho\tau) \R (\rho\tau \wedge_{\L} 
\rho\sigma) \}.
\end{array}$$
Note that $(\sigma,\tau) \in V(M^I)$ implies in particular that
$(\sigma \wedge_{\L} \tau) \R (\tau \wedge_{\L} \sigma)$.

Next we define 
$$W(M^I) = \{ m \in M^I \mid \L_{m} \; = \; \H_m\}.$$

\bl
\label{WJclos}
$W(M^I)$ is a union of $\J$-classes of $M^I$.
\el

\proof
Let $m \in W(M^I)$. Since our monoid $M$ is finite $\J$-above, we have
$\J \; = \; \D$ and 
so it suffices to show that $\L_{m} \,\cup\, \R_m \; \subseteq W(M^I)$.

Assume that $m' \L m$. Then $m' \in \; \L_{m} \; = \; \H_m$ and so
$\L_{m'} \; = \; \L_{m} \; = \; \H_m \; = \; \H_{m'}$. Thus $m' \in
W(M^I)$.

Finally, assume that $m' \R m$ and take $u \in\;  \L_{m'}$. Write $m'
= mx$ and $m = m'y$. Then $m' \L u$ yields $m'y \,\L\,
uy$ and so $uy \in \; \L_{m} \; = \; \H_m$. By Green's Lemma, we
get $uyx \,\H \, mx = m'$. Since $u \L m' = m'yx$ yields $uyx = u$, we
obtain $u \H m'$ and so $m' \in
W(M^I)$.
\qed

The next lemma provides an alternative characterization of $V(M^I)$:

\bl
\label{charV}
Let $\sigma = (m_k <_{\L} \ldots <_{\L} m_0)$, $\tau = (m'_l <_{\L}
\ldots <_{\L} m'_0)$ and $(\sigma\wedge_{\L}\tau) = m_i$. Then
$(\sigma,\tau) \in V(M^I)$ if and only if $m_i \R m'_i$ and one of the
following conditions holds:
\bi
\item[(V1)] $i = k = l$;
\item[(V2)] $i < k, l$;
\item[(V3)] $i = k < l$ and $m_i \in W(M^I)$;
\item[(V4)] $i = l < k$ and $m'_i \in W(M^I)$;
\ei
\el

\proof
Since $(\sigma\wedge_{\L}\tau) = m_i$, we have
$(\tau\wedge_{\L}\sigma) = m'_i$. Taking $\rho = (I)$ in the condition
defining $V(M^I)$, it becomes 
clear that $m_i \R m'_i$ is a necessary condition for $(\sigma,\tau)
\in V(M^I)$. 

Assume that $m_i \R m'_i$ and one of conditions (V1)--(V4) holds. 
Write $\rho = (n_p <_{\L} \ldots <_{\L} n_0)$ and assume that 
$(\rho\sigma \wedge_{\L} 
\rho\tau) \L (\sigma \wedge_{\L} 
\tau)$. We have
$$\rho\sigma = \lm(n_pm_k \leq_{\L} \ldots \leq_{\L} n_1m_k
\leq_{\L} m_k <_{\L} \ldots <_{\L} m_0),$$
$$\rho\tau = \lm(n_pm'_l \leq_{\L} \ldots \leq_{\L} n_1m'_l
\leq_{\L} m'_l <_{\L} \ldots <_{\L} m'_0).$$

It should be clear that if (V2) holds, then $(\rho\sigma \wedge_{\L} 
\rho\tau) = m_i$ and $(\rho\tau \wedge_{\L} 
\rho\sigma) = m'_i$, hence $(\sigma,\tau) \in V(M^I)$.

Since $(\rho\sigma \wedge_{\L} 
\rho\tau) \L (\sigma \wedge_{\L} 
\tau)$ if and only if $(\rho\tau \wedge_{\L} 
\rho\sigma) \L (\tau \wedge_{\L} 
\sigma)$, it follows that $V(M^I)$ is a symmetric relation. Thus we
may assume that $i = k$. Let
$$j = \max\{ r \in \{ 0, \ldots, p\} \mid n_rm_k \,\L\, m_k\}.$$
Since $(\rho\sigma \wedge_{\L} 
\rho\tau) \L (\sigma \wedge_{\L} 
\tau) = m_k$, we have $(\rho\sigma \wedge_{\L} 
\rho\tau) = n_jm_k$. 

Assume first that $k < l$ (case (V3)). Then $m_k \in W(M^I)$ and so 
$(\rho\sigma \wedge_{\L} 
\rho\tau)  \,\L\, (\sigma \wedge_{\L} 
\tau) = m_k$ yields $$(\rho\sigma \wedge_{\L} 
\rho\tau)  \,\H\, m_k \R m'_k = (\rho\tau \wedge_{\L} 
\rho\sigma).$$
Hence $(\sigma,\tau)
\in V(M^I)$. 

It remains to be considered the case $k = l$ (case (V1)). Then $m'_k
\H m_k$ and so $n_jm_k\,\L\, m_k$ yields $n_jm'_k\,\L\, m'_k$. By
symmetry, we obtain $(\rho\tau \wedge_{\L} 
\rho\sigma) = n_jm'_k$. Since $m'_k
\H m_k$ implies $n_jm_k\,\R\, n_jm'_k$, it follows that $(\sigma,\tau)
\in V(M^I)$ also in this case.

To prove the converse implication, we assume that the necessary
condition $m_i \R m'_i$ holds
but none of the conditions (V1)--(V4) is satisfied. By symmetry, we
may assume that $i = k < l$ and $m_k \notin W(M^I)$. Then there exists
some $n \in \; \L_{m_k}\setminus \H_{m_k}$, say $n = xm_k$. Let $\rho = (x <_{\L}
I)$. It follows easily that $(\rho\sigma \wedge_{\L} 
\rho\tau) = xm_k = n$ and $(\rho\tau \wedge_{\L} 
\rho\sigma) = m'_k$. Since $m'_k \R m_k$ and $n \not\R m_k$, we get 
$(\rho\sigma \wedge_{\L} 
\rho\tau) \,\not\R \, (\rho\tau \wedge_{\L} 
\rho\sigma)$ and so $(\sigma,\tau)
\notin V(M^I)$ as required.
\qed

\bc
\label{charVi}
Let $\sigma,\tau \in \E(M^I)$ be such that $(\sigma\wedge_{\L}\tau)
\in W(M^I)$. Then 
$(\sigma,\tau) \in V(M^I)$.
\ec

\proof
Let $\sigma = (m_k <_{\L} \ldots <_{\L} m_0)$, $\tau = (m'_l <_{\L}
\ldots <_{\L} m'_0)$ and $(\sigma\wedge_{\L}\tau) = m_i$. Then
$(\tau\wedge_{\L}\sigma) = m'_i$. Since $m'_i \L m_i \in W(M^I)$, it
follows that $m'_i \H m_i$. Hence also $m'_i \in W(M^I)$ by Lemma
\ref{WJclos}. Now we obtain
$(\sigma,\tau) \in V(M^I)$ by Lemma \ref{charV}.
\qed

We define a mapping $H: \E(M^I) \times \E(M^I) \to \oo{\N}$ by
$$H(\sigma,\tau) = \left\{
\begin{array}{ll}
2\sup h_{\J}+2&\mbox{ if } \sigma = \tau\\
2h_{\J}(\sigma\wedge_{\L}\tau)+1&\mbox{ if } \sigma \neq \tau \mbox{ and }
(\sigma,\tau) \in V(M^I)\\
2h_{\J}(\sigma\wedge_{\L}\tau)&\mbox{ otherwise.}
\end{array}
\right.$$
If $M$ is an $Y$-semigroup, then $M^I$ is an $Y$-monoid. We denote by
$H_Y$ the restriction of $H$ to $\E_Y(M^I) \times \E_Y(M^I)$.

\bl
\label{hlf}
Let $M$ be a finite $\J$-above $Y$-semigroup. Then 
\bi
\item[(i)] $H$ is a 
strict length function for $\E(M^I)$;
\item[(ii)] $H_Y$ is a 
strict length function for $\E_Y(M^I)$;
\item[(iii)]
$H_Y = D_{\chi}$ for some (unique up to isomorphism) strongly faithful elliptic
$\E_Y(M^I)$-tree $\chi$.
\ei 
\el

\proof
(i) We show that $H$ satisfies
axioms (L1) -- (L5). By Corollary \ref{holwt} and Proposition \ref{hj}, we
may consider the length function 
$D: \E(M^I) \times \E(M^I) \to \oo{\N}$ defined by
$$D(\sigma,\tau) = \left\{
\begin{array}{ll}
\sup h_{\J}+1&\mbox{ if } \sigma = \tau\\
h_{\J}(\sigma\wedge_{\L}\tau)&\mbox{ otherwise.}
\end{array}
\right.$$
Clearly, $H'(\sigma,\tau) = 2D(\sigma,\tau)$ defines also a length
function for $\E(M^I)$. We shall make use of $H'$ and perform the necessary
adaptations. 
  
Axiom (L1) follows from $V(M^I)$ being a symmetric relation (see the
proof of Lemma \ref{charV}). 
Axioms (L2) and (L5) can be verified for $H$ straightforwardly
as in the proof of Theorem \ref{hol}. We concentrate our efforts on
(L3) and (L4).

(L3) Let $\sigma,\tau,\rho \in \E(M^I)$ and assume that $\sigma \neq
\tau$. By (\ref{gene}), we may assume that $\rho = (m <_{\L}
I)$. By (\ref{obib}), we have  
$(\sigma\rho \wedge_{\L} \tau\rho) \leq_{\J} (\sigma \wedge_{\L}
\tau)$. If
$(\sigma\rho \wedge_{\L} \tau\rho) <_{\J} (\sigma \wedge_{\L} \tau)$,
then
 $$H(\sigma\rho,\tau\rho) \geq
2h_{\J}(\sigma\rho \wedge_{\L} \tau\rho) > 2h_{\J}(\sigma \wedge_{\L}
\tau)+1 \geq H(\sigma,\tau).$$
Thus we may assume that 
\beq
\label{vlk2}
(\sigma\rho \wedge_{\L} \tau\rho) \J (\sigma \wedge_{\L}
\tau).
\eeq
It suffices to show that
\beq
\label{vlk1}
(\sigma,\tau) \in V(M^I) \Rw (\sigma\rho,\tau\rho) \in V(M^I).
\eeq
Assume that $(\sigma,\tau) \in V(M^I)$ and write $$\sigma = (m_k
<_{\L} \ldots <_{\L} m_1 <_{\L} m_0 
= I),\quad \tau = (m'_l <_{\L} \ldots <_{\L} m'_1 <_{\L} m'_0 = I).$$
Then
$$\sigma\rho = \lm(m_km
\leq_{\L} \ldots \leq_{\L} m_1m <_{\L} m <_{\L} I),$$
$$\tau\rho = \lm(m'_lm \leq_{\L} \ldots \leq_{\L} m'_1m <_{\L} m <_{\L} I).$$
We use Lemma \ref{charV}. In particular, we know that $(\sigma \wedge_{\L}
\tau) \H (\tau \wedge_{\L}
\sigma)$.

Suppose first that $(\sigma,\tau)$ satisfies (V1). Then $(\sigma \wedge_{\L}
\tau) = m_k \H m'_k = (\tau \wedge_{\L}
\sigma)$ and $m_km \,\L \, m'_km$ yields $(\sigma\rho \wedge_{\L}
\tau\rho) = m_km$ and $(\tau\rho \wedge_{\L}
\sigma\rho) = m'_km$. Now (\ref{vlk2}) yields $m_km \J m_k$ and
therefore $m_km \R m_k$ by (S1). Similarly, $m'_km \R m'_k$. It
follows that $m_km \R m'_km$ and $(\sigma\rho,\tau\rho)$ satisfies
(V1), thus (\ref{vlk1}) holds in this case.

Suppose next that $(\sigma,\tau)$ satisfies (V2). Then $(\sigma \wedge_{\L}
\tau) = m_i$ implies $(\sigma\rho \wedge_{\L}
\tau\rho) = m_im$. Indeed, It is clear that $(\sigma\rho \wedge_{\L}
\tau\rho) = m_jm$ for some $j \geq i$ since $m_im \,\L\, m'_im$ and
$m_rm = m'_rm$ for $r < i$. However, if $j > i$, then (\ref{stbl})
yields $m_jm \leq_{\J} m_j <_{\J} m_i$, contradicting
(\ref{vlk2}). Hence $(\sigma\rho \wedge_{\L}
\tau\rho) = m_im$. Similarly, $(\tau\rho \wedge_{\L}
\sigma\rho) = m'_im$. Similarly to the preceding case, we get
$m_im \,\R \, m_i \R m'_i \,\R \, m'_im$. Moreover,
$(\sigma\rho,\tau\rho)$ satisfies 
(V2), thus (\ref{vlk1}) holds in this case as well.  

Finally, we assume that $(\sigma,\tau)$ satisfies (V3) (the case (V4)
is dual). Similarly to the preceding cases, we get  $(\sigma\rho \wedge_{\L}
\tau\rho) = m_km$, $(\tau\rho \wedge_{\L}
\sigma\rho) = m'_km$ and $m_km \R m'_km$. Now (\ref{vlk2}) is
equivalent to $m_km \J m_k$ and so
$m_k \in W(M^I)$ yields
$m_km \in W(M^I)$ by Lemma \ref{WJclos}, hence $(\sigma\rho,\tau\rho)$
satisfies  
(V3). Thus (\ref{vlk1}) holds and (L3) is satisfied.

(L4): Let $\sigma,\tau,\rho \in \E(M^I)$. We may assume that
$\sigma,\tau,\rho$ are all distinct. 
Since $H'$ is a length function, we have
$$H'(\sigma,\rho) \geq \min\{ H'(\sigma,\tau),
H'(\tau,\rho) \}.$$
Since $H(x,y) = H'(x,y)$ or $H(x,y) = H'(x,y)+1$
for all $x,y \in \E(M^I)$, we may assume that 
\beq
\label{hlf4}
H'(\sigma,\rho) =
\min\{ H'(\sigma,\tau), 
H'(\tau,\rho) \}.
\eeq
By (L1), we may further assume that
\beq
\label{hlf3}
H'(\sigma,\rho) = H'(\sigma,\tau).
\eeq 
If $H(\sigma,\tau) =
H'(\sigma,\tau)$ we are done, hence assume also that $H(\sigma,\tau) =
H'(\sigma,\tau)+1$, that is,
$$(\sigma,\tau) \in V(M^I).$$
Similarly, we may assume that 
\beq
\label{hlf2}
H'(\sigma,\rho) = H'(\tau,\rho) \Rw
(\tau,\rho) \in V(M^I),
\eeq
otherwise $H(\sigma,\rho) \geq H(\tau,\rho)$.
In view of (\ref{hlf4}), to prove (L4) it suffices to show that
$(\sigma,\rho) \in V(M^I)$. 

Assume first that $H'(\sigma,\rho) = H'(\tau,\rho)$. By (\ref{hlf2}),
we have $(\tau,\rho) \in V(M^I)$. Moreover, $H'(\sigma,\rho) = H'(\sigma,\tau) =
H'(\tau,\rho)$ and Lemma \ref{kkk}(i) yield 
$$(\sigma \wedge_{\L}
\rho) = (\sigma \wedge_{\L} \tau) \,\R\, (\tau \wedge_{\L} \sigma) =
(\tau \wedge_{\L} \rho) \,\R\, (\rho \wedge_{\L} \tau) =  (\rho
\wedge_{\L} \sigma).$$  
We discuss now the cases (V1)--(V4).

If $(\sigma,\tau)$ satisfies (V1), then $(\tau,\rho)$
must satisfy either (V1) or (V3), and so $(\sigma,\rho)$
satisfies (V1) or (V3) accordingly in view of Lemma \ref{WJclos}.

If $(\sigma,\tau)$ satisfies (V2), then $(\tau,\rho)$
must satisfy either (V2) or (V4), and so $(\sigma,\rho)$
satisfies (V2) or (V4) accordingly.

If $(\sigma,\tau)$ satisfies (V3), then $(\tau,\rho)$
must satisfy either (V2) or (V4). In the first case, $(\sigma,\rho)$
satisfies (V3). In the latter, $(\sigma,\rho)$
satisfies (V1) or (V3).

Finally, if $(\sigma,\tau)$ satisfies (V4), then $(\tau,\rho)$
must satisfy either (V1) or (V3). In the first case, $(\sigma,\rho)$
satisfies (V4) by Lemma \ref{WJclos}. In the latter, $(\sigma,\rho)$
satisfies (V2). This completes the discussion of the case
$H'(\sigma,\rho) = H'(\tau,\rho)$. 

It remains to be considered the case $H'(\sigma,\rho) < H'(\tau,\rho)$. By
(\ref{hlf3}) and Lemma \ref{kkk}(i), we have 
\beq
\label{eli4}
(\sigma \wedge_{\L}\rho) =
(\sigma\wedge_{\L} \tau)
\eeq
and so
\beq
\label{eli2}
(\tau \wedge_{\L}
\sigma) \L (\rho 
\wedge_{\L} \sigma).
\eeq
Since $H'(\tau,\sigma) = H'(\sigma,\tau) =
H'(\sigma,\rho) < 
H'(\tau,\rho)$, $\tau \wedge_{\L} \sigma$ must be a term of 
$\rho$ with $\tau \wedge_{\L} \sigma >_{\L} \rho \wedge_{\L} 
\tau$. Hence (\ref{eli2}) yields
\beq
\label{eli3}
(\tau \wedge_{\L} \sigma) = (\rho
\wedge_{\L} \sigma).
\eeq
Since $(\sigma,\tau) \in V(M^I)$, it follows from (\ref{eli4}) and
(\ref{eli3}) that  
$$(\sigma \wedge_{\L} \rho)\,\R \, (\rho \wedge_{\L} \sigma).$$
We discuss now the cases (V1)--(V4).

Clearly, $H'(\sigma,\rho) < H'(\tau,\rho)$ implies that
$(\sigma,\tau)$ must satisfy either (V2) or (V3). It is easy to see
that $(\sigma,\rho)$ satisfies necessarily the same condition, hence
$(\sigma,\rho) \in V(M^I)$ and so (L4) holds. 

Therefore $H$ is a strict
length function for $\E(M^I)$.

(ii) Since $\E_Y(M^I)$ is a submonoid of $\E(M^I)$, the restriction of
$H$ to $\E_Y(M^I) \times \E_Y(M^I)$ is a length function for
$\E_Y(M^I)$.

(iii) We get
$H_Y = D_{\chi}$ for some (unique up to isomorphism) strongly faithful elliptic
$\E_Y(M^I)$-tree $\chi$ by Corollary \ref{flftree}.
\qed

Throughout the remaining part of this section, we assume that $M$, $H_Y
= D_{\chi}$ for $\chi = (r_0,T,\alpha,\theta)$ are fixed. Moreover, we 
may assume that $\chi$ is obtained by the Chiswell construction
according to the
proofs of Theorem \ref{lftree} and Corollary 
\ref{flftree}. 

We say that $[n, m_l <_{\L} \ldots <_{\L} m_0 = I] = v$ is a {\em
  minimal representation} of $v \in \vert(T)$ if $v \neq [n, m_i
<_{\L} \ldots <_{\L} m_0]$ for every $i < l$. 

The following lemma helps to establish that an $\L$-chain
belongs to $\E_Y(M^I)$:

\bl
\label{cutmemb}
 Let $(m_l <_{\L} m_{l-1} <_{\L} \ldots <_{\L} m_0) \in \E_Y(M^I)$.
\bi
\item[(i)]
For every $i < l$, $(m_i
<_{\L} \ldots <_{\L} m_0) \in \E_Y(M^I)$.
\item[(ii)]
If $m'_l \L
 m_l$, then $(m'_l <_{\L} m_{l-1} <_{\L} \ldots <_{\L} m_0) \in
 \E_Y(M^I)$.
\ei
\el

\proof
(i) If $(m_l
<_{\L} \ldots <_{\L} m_0) = (y_r <_{\L} I)\ldots (y_1 <_{\L} I)$, then 
\beq
\label{gedi2}
(m_i
<_{\L} \ldots <_{\L} m_0) = (y_s <_{\L} I)\ldots (y_1 <_{\L} I)
\eeq
for $s = \max\{ j < r \mid y_j\ldots y_1 \,\L\, m_i\}$.

(ii) Write $$\sigma = (m_l <_{\L} \ldots <_{\L} m_0),\quad \tau = (m'_l
<_{\L} m_{l-1} <_{\L} \ldots <_{\L} m_0).$$ 
Since $m'_l \L
 m_l$, we have $m'_l = y_r\ldots y_1m_l$
for some $y_1,\ldots,y_r \in Y$. Since $\sigma \in \E_Y(M^I)$, it
follows that 
$$\begin{array}{l}
\lm( y_r\ldots y_1m_l \leq_{\L} y_{r-1}\ldots y_1m_l \leq_{\L}
\ldots \leq_{\L} y_1m_l \leq_{\L}  m_l <_{\L} \ldots <_{\L} m_0
)\\
\hspace{1cm} = (y_r <_{\L} I)\ldots (y_1 <_{\L} I)\sigma \in \E_Y(M^I).
\end{array}$$
Since $m'_l = y_r\ldots y_1m_l \L m_l$, we obtain 
$$\lm( y_r\ldots y_1m_l \leq_{\L} y_{r-1}\ldots y_1m_l \leq_{\L}
\ldots \leq_{\L} y_1m_l \leq_{\L}  m_l <_{\L} \ldots <_{\L}
m_0) = \tau$$
and so $\tau \in \E_Y(M^I)$.
\qed

\bl
\label{norf}
\bi
\item[(i)] Let $v = [2k, m_l <_{\L} \ldots <_{\L} m_0 = I] \in
  \vert(T)$ and $i \in \{ 0, \ldots, l-1\}$. Then $v = [2k, m_i <_{\L}
  \ldots <_{\L} m_0]$ if and only if $h_{\J}(m_i) \geq k$.
\item[(ii)] Let $v = [2k+1, m_l <_{\L} \ldots <_{\L} m_0 = I] \in
  \vert(T)$ and $i \in \{ 0, \ldots, l-1\}$. Then $v = [2k+1, m_i <_{\L}
  \ldots <_{\L} m_0]$ if and only if $h_{\J}(m_i) > k$ or
$$h_{\J}(m_i) = k\mbox{ and } m_i \in W(M^I).$$
\ei
\el

\proof
Let $$\sigma = (m_l <_{\L} \ldots <_{\L} m_0),\quad \tau = (m_{i}
<_{\L} \ldots <_{\L} m_0).$$ Then $(\sigma \wedge_{\L} \tau) =
m_{i} = (\tau \wedge_{\L} \sigma)$.

(i) We have
$$[2k,\sigma] = [2k,\tau] \iff H_Y(\sigma,\tau) \geq 2k \iff
2h_{\J}(m_{i})\geq 2k \iff h_{\J}(m_{i}) \geq k.$$

(ii) Assume first that $m_i \in W(M^I)$. Then
$$H_Y(\sigma,\tau) = 2h_{\J}(\sigma \wedge_{\L} \tau) +1 =
2h_{\J}(m_{i})+1$$
and so
$$[2k+1,\sigma] = [2k+1,\tau] \iff H_Y(\sigma,\tau) \geq 2k+1 \iff
2h_{\J}(m_{i})+1 \geq 2k+1 \iff h_{\J}(m_{i}) \geq k.$$
If $m_i \notin W(M^I)$, then $H_Y(\sigma,\tau) = 
2h_{\J}(m_{i})$ and so
$$[2k+1,\sigma] = [2k+1,\tau] \iff H_Y(\sigma,\tau) \geq 2k+1 \iff
2h_{\J}(m_{i}) \geq 2k+1 \iff h_{\J}(m_{i}) > k.$$
\qed

We immediately obtain

\bc
\label{ferias}
\bi
\item[(i)] $v = [2k, m_l <_{\L} \ldots <_{\L} m_0 = I] \in
  \vert(T)$ is in minimal representation if and only if $h_{\J}(m_{l-1}) < k$.
\item[(ii)] $v = [2k+1, m_l <_{\L} \ldots <_{\L} m_0 = I] \in
  \vert(T)$ is in minimal representation if and only if
  $h_{\J}(m_{l-1}) < k$ or
$$h_{\J}(m_{l-1}) = k\mbox{ and } m_{l-1} \notin W(M^I).$$
\ei
\ec

Given $m \in M^I$, let 
$$Y_m = \{ y \in Y \mid ym <_{\L} m\}.$$
For every $y \in Y_m$, there exists a unique $b \in B_{ym}$ such that
$(1,1,b) \,\L \, ym$.
We denote by $Q_m$ the set of all such $b$ when $y$ takes values in
$Y_m$.

\bl
\label{Rcompat}
For every $m \in M^I$, $Y_m = Y_{mm^*}$ and
$$\begin{array}{rcl}
\p:Q_m&\to&Q_{mm^*}\\
b&\mapsto&((1,1,b)m^*)\pi_3
\end{array}$$
is a bijection.
\el

\proof
Let $m \in M^I$.
We always have $ym \leq_{\L} m$ and  $ymm^* \leq_{\L} mm^*$. Since $m
= mm^*m^{\sharp}$, it is
immediate that $ym \,\L\, m$ if and only if $ymm^* \,\L\, mm^*$, hence
$Y_m = Y_{mm^*}$. 

Let $b \in Q_m$. Then $(1,1,b) \,\L\, ym$ for some $y \in
Y_m$ and so $(1,1,b)m^* \,\L\, ymm^*$. Since $y \in Y_m = Y_{mm^*}$,
it follows that $\p(b) \in  Q_{mm^*}$. Thus $\p$ is well defined.

Suppose now that $\p(b) = \p(c)$. Then $(1,1,b)m^* \,\L\, (1,1,c)m^*$
and so $(1,1,b)m^*m^{\sharp} \,\L\, (1,1,c)m^*m^{\sharp}$. Since
$(1,1,b) <_{\L} m$, we get $(1,1,b)m^*m^{\sharp} =
(1,1,b)$. Similarly, $(1,1,c)m^*m^{\sharp} =
(1,1,c)$ and so $(1,1,b) \,\L\, (1,1,c)$. Thus $b = c$ and $\p$ is
one-to-one.

Finally, let $c \in Q_{mm^*}$. Then $(1,1,c) \,\L\, ymm^*$ for some $y
\in Y_{mm^*} = Y_m$. It follows that $(1,1,c)m^{\sharp} \,\L\,
ymm^*m^{\sharp} = ym$. Write $b = ((1,1,c)m^{\sharp})\pi_3$. Then $b
\in Q_m$. We show that $\p(b) = c$. It suffices to show that
$(1,1,b)m^* \,\L\,(1,1,c)$. Now $(1,1,b) \,\L\,(1,1,c)m^{\sharp}$
yields $(1,1,b)m^* \,\L\,(1,1,c)m^{\sharp}m^*$. Since $(1,1,c) \,\L\,
ymm^*$, we get $(1,1,c)m^{\sharp}m^* = (1,1,c)$ and so $(1,1,b)m^*
\,\L\,(1,1,c)$ as required. Thus $\p$ is onto and therefore a bijection.
\qed 

Given $m = (a,g,b) \in M^I$, define
$$A'_m = \{ a' \in A_m \mid Y_{(a',g,b)} \neq \emptyset\}.$$
For every $k \in \N$, let
\bi
\item[]
$U_0(k) = \{ m \in W(M^I) : h_{\J}(m) = k \mbox{ and } |A_m| > 1\},$
\item[]
$U_1(k) = \{ m \in M^I \setminus W(M^I) : h_{\J}(m) = k \mbox{ and }
|A_m|+|A'_m| > 1\},$ 
\item[]
$U_2(k) = \{ m \in W(M^I) : h_{\J}(m) = k \mbox{ and }
|G_m|(1+|Q_m|) > 1\},$
\item[]
$U_3(k) = \{ m \in M^I \setminus W(M^I) : h_{\J}(m) = k \mbox{ and }
|G_m| > 1\},$
\item[]
$U_4(k) = \{ m \in M^I \setminus W(M^I) : h_{\J}(m) = k \mbox{ and }
|G_m|\cdot |Q_m| > 1\}.$
\ei

\bl
\label{URclos}
$U_i(k)$ is a union of $\R$-classes of $M^I$ for $i = 0,2,3,4$.
\el

\proof
The claim follows from Lemmas \ref{WJclos} and \ref{Rcompat} and
\beq
\label{gedi3}
m \R m' \Rw A'_m = A'_{m'}.
\eeq
We prove that (\ref{gedi3}) holds. Indeed, assume that $m = (a,g,b)$
and $m' = (a,g',b')$ are $\R$-related. By Lemma \ref{Rcompat}, we have  
$Y_{(a',g,b)} = Y_{(a',g',b')}$ for every $a' \in A_m = A_{m'}$. Hence 
$A'_m = A'_{m'}$ and (\ref{gedi3}) holds as required.
\qed

We discuss now the cases when a vertex has more than one son. For every $m
\in M^I$ with $Y_m \neq \emptyset$, we fix an arbitrary element
$\gamma_m \in Y_mm$.

\bl
\label{norfA}
Let $v = [2k, m_l <_{\L} \ldots <_{\L} m_0 = I] \in
  \vert(T)$ be in minimal representation and $2k < \dep(r_0,T)$. Then
  $|\son(v)| > 1$ if and only if 
$m_l \in U_0(k) \cup U_1(k)$. 
In that case, if
  $m_l = (a,g,b)$, then
$$\son(v) = \left\{
\begin{array}{ll}
\son_1(v) &\mbox{ if }m_l \in U_0(k)\\
\son_1(v) \cup \son_2(v)&\mbox{ if }m_l \in U_1(k)
\end{array}
\right.$$ 
with 
$$\begin{array}{lll}
\son_1(v)&=&\{ [2k+1, (a',g,b) <_{\L} m_{l-1}
<_{\L} \ldots <_{\L} m_0];\; a' \in A_{m_l} \},\\
\son_2(v)&=&\{ [2k+1, \gamma_{(a',g,b)} <_{\L} (a',g,b) <_{\L} m_{l-1}
<_{\L} \ldots <_{\L} m_0];\; a' \in A'_{m_l} \}
\end{array}$$
and the represented elements are
all distinct in each case.
\el

\proof
Write $\sigma = (m_l <_{\L} \ldots <_{\L} m_0)$. Since $2k <
\dep(r_0,T)$, we have $|\son(v)| \geq 1$. It follows from
(\ref{lftree2}) that $|\son(v)| = 1$ if and only if
\beq
\label{norfA1}
[2k,\sigma] = [2k,\tau] \Rw [2k+1,\sigma] = [2k+1,\tau]
\eeq
for every $\tau \in \E_Y(M^I)$.

Suppose first that $h_{\J}(m_l) < k$. Since
$$2h_{\J}(\sigma \wedge_{\L} \tau) +1 \leq 2h_{\J}(m_l) +1 \leq 2(k-1)
+1 < 2k,$$
then $[2k,\sigma] = [2k,\tau]$ implies $\sigma = \tau$ and so
(\ref{norfA1}) holds.
Thus $|\son(v)| = 1$.

Suppose now that $h_{\J}(m_l) > k$. Assume that $[2k,\sigma] =
[2k,\tau]$ with $\sigma \neq \tau$. Then
$2h_{\J}(\sigma \wedge_{\L} \tau) +1  \geq H_Y(\sigma,\tau) \geq
2k$. Since $[2k,\sigma]$ is a minimal representation, it follows that
$h_{\J}(m_{l-1}) < k$ and so $(\sigma \wedge_{\L} \tau) = m_l$. Thus
$$H_Y(\sigma,\tau) \geq 2h_{\J}(\sigma \wedge_{\L} \tau) = 2h_{\J}(m_l)
\geq 2k+2$$
and so $[2k+1,\sigma] = [2k+1,\tau]$. Therefore (\ref{norfA1}) holds
and $|\son(v)| = 1$. 

Therefore we assume that $h_{\J}(m_l) = k$ and write $m_l =
(a,g,b)$.

Suppose first that $m_l \in W(M^I)$.  
We show 
that $\son(v) = \son_1(v)$.
For every $a' \in A_{m_l}$, we have $(a',g,b)
\L (a,g,b) = m_l$. Let $\tau = ((a',g,b) <_{\L} m_{l-1}
<_{\L} \ldots <_{\L} m_0)$. By Lemma \ref{cutmemb}(ii), $\tau \in
\E_Y(M^I)$. On the other hand, 
$(\sigma \wedge_{\L} \tau) =
m_l$ and so $H_Y(\sigma,\tau) \geq 2h_{\J}(m_{l}) = 2k$. Hence
$[2k,\tau] = v$ and we conclude by
(\ref{lftree2}) that $[2k+1,\tau] \in \son(v)$. 

Conversely, assume that $[2k+1,\zeta] \in \son(v)$ is in minimal
representation. We show 
that $[2k+1,\zeta]$ is of the claimed form and we may assume that
$\zeta \neq \sigma$. By (\ref{lftree2}), we have
$[2k,\sigma] = [2k,\zeta]$ and so $H_Y(\sigma,\zeta) \geq 2k$. Hence
$(\sigma \wedge_{\L} \zeta) = m_l \in W(M^I)$  and so $(\zeta
\wedge_{\L} \sigma) = \in W(M^I)$ by Lemma \ref{WJclos}. By
Corollary \ref{ferias}(ii), we get 
$\zeta = (m'_l <_{\L} m_{l-1}
<_{\L} \ldots <_{\L} m_0)$ with $m'_l \L m_l$. If $m'_l = (a',g',b')$,
it follows that $b' = b$ and we may (if $g' \neq g$) replace $g'$ by
$g$ to get $\zeta' 
= ((a',g,b) <_{\L} m_{l-1}
<_{\L} \ldots <_{\L} m_0)$
since
$(a',g,b)\R(a',g',b)$ and case (V1) of Lemma \ref{charV} imply
$(\zeta,\zeta') \in V(M^I)$. Hence 
$$H_Y(\zeta,\zeta') = 2h_{\J}((a',g',b))+1 = 2k+1$$
and so $[2k+1,\zeta] = [2k+1,\zeta']$. Thus $\son(v) = \son_1(v)$.

Finally, given
$\rho = ((a'',g,b) <_{\L} m_{l-1}
<_{\L} \ldots <_{\L} m_0)$ with $a'' \neq a'$, then $$(\tau
\wedge_{\L} \rho) = (a',g,b) \,\not\R\, (a'',g,b) = (\rho \wedge_{\L}
\tau)$$
and so $H_Y(\tau,\rho) = 2h_{\J}((a',g,b)) = 2k$. Thus $[2k+1,\tau] \neq
[2k+1,\rho]$ and so the elements in $\son_1(v)$ are
all distinct. 
In particular, $|\son(v)| = |A_{m_l}|$ and so $|\son(v)| > 1$ if and
only if $m_l \in U_0(k)$.

Assume now that $m_l \notin W(M^I)$. We show that  $\son(v) =
\son_1(v) \cup \son_2(v)$. We pass the inclusion $\son_1(v) \cup
\son_2(v) \subseteq \son(v)$, a straightforward adaptation of the
preceding case, and move straight to the converse inclusion. Let 
$[2k+1,\zeta] \in \son(v)$  and assume that
$\zeta \neq \sigma$. By (\ref{lftree2}), we have
$[2k,\sigma] = [2k,\zeta]$ and so $H_Y(\sigma,\zeta) \geq 2k$. Hence
$(\sigma \wedge_{\L} \zeta) = m_l$  and so by Corollary \ref{ferias}(ii)
we must have
$$\zeta = (m'_l <_{\L} m_{l-1}
<_{\L} \ldots <_{\L} m_0) \hspace{.7cm}\mbox{or} \hspace{.7cm}
\zeta = (m'_{l+1} <_{\L} m'_l <_{\L} m_{l-1}
<_{\L} \ldots <_{\L} m_0)$$
with $m_l \L m'_l = (\zeta \wedge_{\L} \sigma)$. The discussion of the
first case is analogous to the case $m_l \in W(M^I)$, hence we assume
that $\zeta = (m'_{l+1} <_{\L} m'_l <_{\L} m_{l-1}
<_{\L} \ldots <_{\L} m_0)$ and $m'_l = (a',g',b)$. Let
$$\zeta' 
= (\gamma_{(a',g,b)} <_{\L} (a',g,b) <_{\L} m_{l-1}
<_{\L} \ldots <_{\L} m_0).$$ Since $\zeta \in \E_Y(M^I)$, it follows
from the maximality of $s$ in (\ref{gedi2}) that $Y_{m'_l} \neq
\emptyset$. Since $m'_l \R (a',g,b)$, it follows from Lemma
\ref{Rcompat} that $Y_{(a',g,b)} \neq \emptyset$ and so $a' \in A'_m$. Thus 
$[2k+1,\zeta'] \in \son_2(v)$. Finally, either $h_{\J}(\zeta\wedge_{\L}\zeta') >
k$, or $(\zeta\wedge_{\L}\zeta') = m'_l$ and so  
$(\zeta,\zeta') \in V(M^I)$ through case (V2) of Lemma \ref{charV}. In
any case, it follows that 
$H_Y(\zeta,\zeta') \geq 2k+1$ 
and so $[2k+1,\zeta] = [2k+1,\zeta']\in \son_2(v)$. Thus $\son(v) =
\son_1(v)\cup \son_2(v)$.

For uniqueness, we only have to care about distinguishing
$[2k+1,\zeta]$ from $[2k+1,\zeta']$ for
$$\zeta 
= ((a',g,b) <_{\L} m_{l-1}
<_{\L} \ldots <_{\L} m_0),\quad
\zeta' 
= (\gamma_{(a',g,b)} <_{\L} (a',g,b) <_{\L} m_{l-1}
<_{\L} \ldots <_{\L} m_0),$$
the remaining cases following the same argument of the case $m_l \in
W(M^I)$.

Since $m_l \notin W(M^I)$, then $(a',g,b) \notin W(M^I)$ by Lemma
\ref{WJclos} and so $(\zeta,\zeta') \notin V(M^I)$ by Lemma
\ref{charV}. Hence
$$H(\zeta,\zeta') = 2h_{\J}(\zeta\wedge_{\L}\zeta') =
2h_{\J}((a',g,b)) = 2k$$
and so $[2k+1,\zeta]\neq [2k+1,\zeta']$. Thus the elements in
$\son_1(v)\cup \son_2(v)$ are
all distinct. 
In particular, $|\son(v)| = |A_{m_l}|+|A'_{m_l}|$ and so $|\son(v)| > 1$ if and
only if $m_l \in U_1(k)$. 
\qed

Note that  $v = [2k, m_l <_{\L} \ldots <_{\L} m_0 = I]$ with $m_l \in
U_0(k)\cup U_1(k)$ implies $2k < \dep(r_0,T)$ and so $|\son(v)| > 1$.

\bl
\label{norfB}
Let $v = [2k+1, m_l <_{\L} \ldots <_{\L} m_0 = I] \in
  \vert(T)$ be in minimal representation and $2k+1 < \dep(r_0,T)$. Then
  $|\son(v)| > 1$ if and only if $m_l \in U_2(k) \cup U_3(k)$ or $m_{l-1}
  \in U_4(k)$.  In that case,
  $$\son(v) = \left\{
\begin{array}{ll}
\son_1(v) \cup \son_2(v)&\mbox{ if }m_l \in U_2(k)\\
\son_1(v) &\mbox{ if }m_l \in U_3(k)\\
\son_3(v) &\mbox{ if }m_{l-1} \in U_4(k)
\end{array}
\right.$$ 
with 
$$\begin{array}{lll}
\son_1(v)&=&\{ [2k+2, m'_l <_{\L} m_{l-1}
<_{\L} \ldots <_{\L} m_0];\; m'_l \in \; \H_{m_l} \},\\
\son_2(v)&=&\{ [2k+2, (1,1,b') <_{\L} m'_l <_{\L} m_{l-1}
<_{\L} \ldots <_{\L} m_0];\; m'_l \in \; \H_{m_l},\; b' \in
Q_{m'_l}\}\\
\son_3(v)&=&\{ [2k+2, (1,1,b') <_{\L}  m'_{l-1} <_{\L}  m_{l-2}
<_{\L} \ldots <_{\L} m_0];\;  m'_{l-1} \in \; \H_{m_{l-1}},\; b' \in
Q_{m'_{l-1}}\} 
\end{array}$$
and the represented elements are
all distinct in each case.
\el

\proof
Write $\sigma = (m_l <_{\L} \ldots <_{\L} m_0)$. Since $2k+1 <
\dep(r_0,T)$, we have $|\son(v)| \geq 1$. 
By (\ref{lftree2}), $|\son(v)| = 1$ if and only if
\beq
\label{flip1}
[2k+1,\sigma] = [2k+1,\tau] \Rw [2k+2,\sigma] = [2k+2,\tau]
\eeq
for every $\tau \in \E_Y(M^I)$.

The case
$h_{\J}(m_l) < k$ is discussed analogously to the proof of Lemma
\ref{norfA}.

Assume next that $h_{\J}(m_l) = k$ and $m_l \in W(M^I)$.  We show 
that
$\son(v) = \son_1(v) \cup \son_2(v)$.

Let $m'_l \in \; \H_{m_l}$ and write $\tau = (m'_l <_{\L} m_{l-1}
<_{\L} \ldots <_{\L} m_0)$. By Lemma \ref{cutmemb}(ii), $\tau \in
\E_Y(M^I)$. Since $(\sigma \wedge_{\L} \tau) =
m_l \in W(M^I)$, we get $(\sigma, \tau) \in V(M^I)$ by
Corollary \ref{charVi}, hence $H_Y(\sigma,\tau) =
2h_{\J}(m_{l})+1 = 2k+1$.  
By
(\ref{lftree2}), we conclude that $[2k+2,\tau] \in \son(v)$.

Assume now that $b' \in Q_{m'}$. Then $(1,1,b')\,\L\, ym'_l$ for some
$y\in Y$ such that $ym'_l 
<_{\L} m'_l$. Write $$\rho = ((1,1,b')
<_{\L} m'_l <_{\L} m_{l-1} 
<_{\L} \ldots <_{\L} m_0), \quad \rho' = (ym'_l
<_{\L} m'_l <_{\L} m_{l-1} 
<_{\L} \ldots <_{\L} m_0).$$
It is immediate that $\rho' = (y <_{\L} I)\tau \in \E_Y(M^I)$. Since
$(1,1,b') \,\L\, ym'_l$, it follows from Lemma \ref{cutmemb}(ii) that
$\rho \in \E_Y(M^I)$ as well.
Now we have
$(\tau \wedge_{\L} \rho) = m'_l \in W(M^I)$ by Lemma \ref{WJclos} and
so Corollary \ref{charVi} yields
$$H_Y(\tau,\rho) = 2h_{\J}(m'_l)+1 = 2h_{\J}(m_l) +1 = 2k+1.$$
Hence $[2k+1,\rho] = [2k+1,\tau] = [2k+1,\sigma]$ and so  $[2k+2,\rho]
\in \son(v)$ as well.

Conversely, let $[2k+2,\zeta] \in \son(v)$ be a minimal representation. We show
that $[2k+2,\zeta]$ is of the claimed form and we may assume that
$\zeta \neq \sigma$. By (\ref{lftree2}), we have
$[2k+1,\sigma] = [2k+1,\zeta]$ and so $H_Y(\sigma,\zeta) \geq 2k+1$. 
It follows that $(\sigma
\wedge_{\L} \zeta) = m_l$. 
Let
$m'_l = (\zeta \wedge_{\L} \sigma)$. Then
$m_l \H m'_l$ since $m_l \in W(M^I)$, and $\zeta = (\ldots m'_l <_{\L} m_{l-1}
<_{\L} \ldots <_{\L} m_0)$. If $m'_l$ is the leftmost term of $\zeta$,
we are done. Otherwise, it follows from Corollary \ref{ferias}(i) and
$h_{\J}(m'_l) = k$ that $\zeta = (m'_{l+1} <_{\L} m'_l <_{\L} m_{l-1}
<_{\L} \ldots <_{\L} m_0)$ for some $m'_{l+1}$. 
Assume that $\zeta = (y_r <_{\L} I)\ldots (y_1 <_{\L} I)$ with
$y_1,\ldots,y_r \in Y$. Then
$$(m'_{l+1} <_{\L} m'_l <_{\L} m_{l-1}
<_{\L} \ldots <_{\L} m_0) = 
\lm( y_r\ldots y_1 \leq_{\L} y_{r-1}\ldots y_1 \leq_{\L}
\ldots \leq_{\L} y_1 <_{\L} I).$$
Let $s = \max\{ j < r\mid y_j\ldots y_1 = m'_l\}$. Then $y_{s+1}\ldots
y_1  <_{\L}  y_s\ldots y_1$ since otherwise, by maximality of $s$,
$m'_l$ would not be the
leftmost element in its $\L$-class. Moreover,
$y_{s+1}m'_l = y_{s+1}\ldots
y_1 \,\L\, y_{r}\ldots
y_1 = m'_{l+1}$.
Let $$\zeta' = (y_{s+1} <_{\L}
I)\ldots (y_1 <_{\L} I) = (y_{s+1}m'_l <_{\L} m'_l <_{\L} m_{l-1}
<_{\L} \ldots <_{\L} m_0)$$ 
and write $y_{s+1}m'_l = (a',g',b')$, $$\zeta'' = ((1,1,b') <_{\L}
m'_l <_{\L} m_{l-1} 
<_{\L} \ldots <_{\L} m_0).$$ 
Clearly, $\zeta' \in \E_Y(M^I)$ and so $\zeta'' \in \E_Y(M^I)$ by
Lemma \ref{cutmemb}(ii). Moreover, $m'_{l+1} \,\L\, y_{s+1}m'_l
\,\L\,(1,1,b')$ yields 
$$H_Y(\zeta,\zeta'') \geq 2h_{\J}(m'_{l+1}) \geq 2k+2$$
and so $[2k+2,\zeta] = [2k+2,\zeta'']$. Since $y_{s+1} \in Y_{m'_l}$
and $b' \in Q_{m'_l}$, this completes the proof of
$\son(v) = \son_1(v) \cup \son_2(v)$.

Finally, suppose that $[2k+2,\tau]$ and $[2k+2,\rho]$ are two sons of the
described form with $\tau \neq \rho$. Then $(\tau \wedge_{\L}
\rho) = m'_l$ for some $m'_l \H m_l$. It follows that
$$H_Y(\tau,\rho) \leq 2h_{\J}(m'_l) +1 = 2h_{\J}(m_l) +1 = 2k+1$$
and so $[2k+2,\tau] \neq [2k+2,\rho]$. Thus the claimed elements of
$\son(v)$ are 
all distinct. 

By Lemma \ref{Rcompat}, we have 
$|\son(v)| = |G_{m_l}|(1+|Q_{m_l}|)$ and so $|\son(v)| > 1$ if and
only if $m_l \in U_2(k)$. 

Assume next that $h_{\J}(m_l) = k$ and $m_l \notin W(M^I)$.  We show 
that
$\son(v) = \son_1(v)$.

Let $m'_l \in \; \H_{m_l}$ and write $\tau = (m'_l <_{\L} m_{l-1}
<_{\L} \ldots <_{\L} m_0)$. By Lemma \ref{cutmemb}(ii), $\tau \in
\E_Y(M^I)$. Since $(\sigma \wedge_{\L} \tau) = m_l \R m'_l = (\tau
\wedge_{\L} \sigma)$, we are in case (V1) of Lemma \ref{charV} and so
$(\sigma, \tau) \in V(M^I)$. Hence $H_Y(\sigma,\tau) =
2h_{\J}(m_{l})+1 = 2k+1$.  
By (\ref{lftree2}), we conclude that $[2k+2,\tau] \in \son(v)$.

Conversely, let $[2k+2,\zeta] \in \son(v)$. We show
that $[2k+2,\zeta]$ is of the claimed form and we may assume that
$\zeta \neq \sigma$. By (\ref{lftree2}), we have
$[2k+1,\sigma] = [2k+1,\zeta]$ and so $H_Y(\sigma,\zeta) \geq 2k+1$. 
It follows that $(\sigma
\wedge_{\L} \zeta) = m_l$ and $(\sigma,\zeta) \in V(M^I)$. Since $m_l
\notin W(M^I)$, it follows from Lemma \ref{charV} that
$(\sigma,\zeta)$ must be in case (V1),and so $\zeta = (m'_l <_{\L} m_{l-1}
<_{\L} \ldots <_{\L} m_0)$ with $m'_l \R m_l$. Since $m_l = (\sigma
\wedge_{\L} \tau) \L (\tau \wedge_{\L} \sigma) = m'_l$, we get
$[2k+2,\zeta] \in \son_1(v)$. 

Proving that the elements of $\son_1(v)$ are distinct is similar to
the preceding case. Therefore 
$|\son(v)| = |G_{m_l}|$ and so $|\son(v)| > 1$ if and
only if $m_l \in U_3(k)$. 

We consider now the case $h_{\J}(m_l) > k$. Suppose first that
$h_{\J}(m_{l-1}) < k$ and take $[2k+1,\sigma] = 
[2k+1,\tau]$ with $\sigma \neq \tau$.  then
$2h_{\J}(\sigma \wedge_{\L} \tau)+1 \geq H_Y(\sigma,\tau) \geq
2k+1$ and so
$(\sigma \wedge_{\L} \tau) = m_l$. Thus
$$H_Y(\sigma,\tau) \geq 2h_{\J}(\sigma \wedge_{\L} \tau) = 2h_{\J}(m_l)
\geq 2k+2$$
and so $[2k+2,\sigma] = [2k+2,\tau]$. Therefore (\ref{flip1}) holds
and $|\son(v)| = 1$. 

Since $v$ is in minimal
representation, we may assume now by Corollary \ref{ferias}(ii) that
$h_{\J}(m_{l-1}) = k$ and $m_{l-1} \notin W(M^I)$.
We show that
$\son(v) = \son_3(v)$.

Let $m'_{l-1} \in\; \H_{m_{l-1}}$ and $b' \in Q_{m'_{l-1}}$. Then
$(1,1,b')\,\L\, ym'_{l-1}$ for some 
$y\in Y$ such that $ym'_{l-1} 
<_{\L} m'_{l-1}$. Write $$\rho = ((1,1,b')
<_{\L} m'_{l-1} 
<_{\L} m_{l-2} 
<_{\L} \ldots <_{\L} m_0),$$
$$\rho' = (ym'_{l-1} 
<_{\L} m'_{l-1} <_{\L} m_{l-2} 
<_{\L} \ldots <_{\L} m_0).$$
It is immediate that $$\rho' = (y <_{\L} I)(m'_{l-1} <_{\L} m_{l-2} 
<_{\L} \ldots <_{\L} m_0) \in \E_Y(M^I).$$ Since
$(1,1,b') \,\L\, ym'_{l-1}$, it follows from Lemma \ref{cutmemb}(ii) that
$\rho \in \E_Y(M^I)$ as well. The case $(\sigma \wedge_{\L} \rho) = m_{l}$ is
straightforward, hence we assume that $(\sigma \wedge_{\L} \rho) = m_{l-1}$.
Thus
$(\rho \wedge_{\L} \sigma) = m'_{l-1} \H m_{l-1}$ and
we are in case (V2) of Lemma \ref{charV}, yielding $(\sigma,\rho) \in
V(M^I)$. It follows that
$H_Y(\sigma,\rho) = 2h_{\J}(m_{l-1})+1 = 2k+1$ and so
$[2k+1,\rho] = [2k+1,\sigma]$. Therefore $[2k+2,\rho]
\in \son(v)$.

Conversely, let $[2k+2,\zeta] \in \son(v)$ be in minimal
representation. We show 
that $[2k+2,\zeta] \in \son_3(v)$.
By (\ref{lftree2}), we have
$[2k+1,\sigma] = [2k+1,\zeta]$ and so $H_Y(\sigma,\zeta) \geq 2k+1$. 
It follows that $(\sigma
\wedge_{\L} \zeta) = m_l$ or else
\beq
\label{flip4}
(\sigma
\wedge_{\L} \zeta) = m_{l-1} \hspace{.7cm}\mbox{and}
\hspace{.7cm} (\sigma, \zeta) \in V(M^I).
\eeq
Suppose that (\ref{flip4}) holds. 
Let $m'_{l-1} = (\zeta \wedge_{\L} \sigma)$. Then
$m_{l-1} \H m'_{l-1}$ and $\zeta = (\ldots m'_{l-1} <_{\L} m_{l-2}
<_{\L} \ldots <_{\L} m_0)$. If $m'_{l-1}$ is the leftmost term of
$\zeta$, then $(\sigma, \zeta)$ would be in case (V4) of Lemma
\ref{charV} and so $m'_{l-1} \in W(M^I)$, contradicting $m_{l-1} \notin
W(M^I)$ in view of Lemma \ref{WJclos}.
On the other hand, since $[2k+2,\zeta]$ is in minimal
representation, it follows from Corollary \ref{ferias}(i) and
$h_{\J}(m'_{l-1}) = k$ that $\zeta = (m'_{l} <_{\L} m'_{l-1} <_{\L} m_{l-2}
<_{\L} \ldots <_{\L} m_0)$ for some $m'_{l}$. Now the proof that
$[2k+2,\zeta] \in \son_3(v)$ is completely analogous to the case
$h_{\J}(m_l) = k$ and $m_l \in W(M^I)$, and is therefore omitted. The
same arguments hold for the
case $(\sigma \wedge_{\L} \zeta) = m_l$, which is actually
simpler. Therefore $\son(v) = \son_3(v)$.

Proving that the elements of $\son_1(v)$ are distinct is similar to
the preceding case.
By Lemma \ref{Rcompat}, we have 
$|\son(v)| = |G_{m_l}|\cdot |Q_{m_l}|$ and so $|\son(v)| > 1$ if and
only if $m_{l-1} \in U_4(k)$.
\qed

Note that  $v = [2k+1, m_l <_{\L} \ldots <_{\L} m_0 = I]$ with $m_l \in
U_2(k)$ implies $2k+1 < \dep(r_0,T)$ and so $|\son(v)| > 1$.

\bl
\label{ferias1}
Let $v = [i, m_l <_{\L} \ldots <_{\L} m_0 = I] \in
  \vert(T)$ be in minimal representation and let $\sigma = (m'_p
  <_{\L} \ldots <_{\L} m_0)$ be such that $m_jm'_p \R m_j$ for some $j
  \in \{ 0, \ldots, l-1\}$. Then $v\sigma = [i, (m_l <_{\L} \ldots
  <_{\L} m_0)\sigma]$ is in minimal representation.
\el

\proof
By successive application of Lemma \ref{gedi}, we get
$$(m_l <_{\L} \ldots
  <_{\L} m_0)\sigma = (m_lm'_p <_{\L} \ldots
  <_{\L} m_jm'_p <_{\L} \ldots)$$
and $m_{l-1}m'_p \R m_{l-1}$. 
Hence $h_{\J}(m_{l-1}m'_p) = h_{\J}(m_{l-1})$. By
Lemma \ref{WJclos}, we also have $m_{l-1} \in W(M^I)$ if and only if
$m_{l-1}m'_p \in W(M^I)$. Thus $[i, (m_l <_{\L} \ldots
  <_{\L} m_0)\sigma]$ is in minimal representation by Corollary \ref{ferias}.
\qed

Assume that $\delta = \dep(r_0,T)$. For commodity, we assume for the
remaining part of this section that $\delta \in
\N$, the infinite case being 
absolutely similar. 
We take two new symbols $\downarrow,\ast$. 
For every $k \in \N$ such that $2k+1 \leq \delta$, let 
$$X_{2k+1} = \{ \downarrow\} \cup(\bigcup_{m \in U_0(k)\cup U_1(k)}
A_{m})
\cup(\bigcup_{m \in U_1(k)} (A'_{m} \times \{ \ast \})).$$
For every $k \in \N$ such that $2k+2 \leq \delta$, let
$$\begin{array}{lll}
X_{2k+2}&=&\{ \downarrow\} \cup(\bigcup_{m \in U_2(k)} (G_{m} \times
(\{ \ast \} 
\cup Q_{mm^*})))\\
&\cup&(\bigcup_{m \in U_3(k)} (G_{m} \times
\{ \ast \})) \cup (\bigcup_{m \in U_4(k)} 
(G_{m} \times Q_{mm^*})).
\end{array}$$
A very important remark: in view of Lemma \ref{Rcompat} and
(\ref{gedi3}), 
{\em we assume the union over} $m \in U_i(k)$
{\em to be disjoint over distinct $\R$-classes}, e.g.: if $m,m' \in
U_2(k)$ are $R$-related, i.e. $mm^* = m'(m')^*$, then $G_{m} \times
(\{ \ast \} 
\cup Q_{mm^*})) = G_{m'} \times
(\{ \ast \} 
\cup Q_{m'(m')^*}))$. Otherwise, they are disjoint. 

If $M$ is finitely generated, then the $X_i$ turn out to be finite:

\bl
\label{finiX}
If $Y$ is finite, then all $X_i$ are
finite.
\el

\proof
It is enough to show that each set
$$E_k = \{ m \in M^I \mid h_{\J}(m) = k\}$$
is finite. Since $M$ is finite $\J$-above, this follows easily by
induction on $k$ from  $E_0 = \{
I\}$ and
\beq
\label{flip6}
E_k \subseteq \bigcup_{i = 0}^{k-1} \bigcup_{x \in YE_{i-1}}  \J_{x}.
\eeq
Indeed, if $m = y_s\ldots y_1 \in E_k$ with $y_i \in Y$, take $$r =
\max\{ j \in \{ 0,\ldots,s\} \mid m <_{\J} y_j\ldots y_1\}.$$
Let $n = y_r\ldots y_1$. Then $n \in E_{i}$ for some $i \in \{
0,\ldots,k-1\}$ and $m \in
\;\J_{y_{r+1}n}$, hence (\ref{flip6}) holds and so does the lemma. 
\qed

In view of Lemmas \ref{norfA} and \ref{norfB}, we define a mapping
$f:(\vert(T))\setminus \{r_0\} \to \cup_{i=1}^{\delta} 
X_i$ as follows. Let 
$v \in \vert(T)$ and let $w \in \son(v)$. 
\bi
\item[(F1)]
If $\son(v) = \{ w \}$, let $f(w) = \downarrow$. 
\item[(F2)]
If $v = [2k, m_l <_{\L} \ldots
<_{\L} m_0]$ with $m_l = (a,g,b) \in U_0(k)\cup U_1(k)$ and $w = [2k+1, (a',g,b)
<_{\L} m_{l-1} <_{\L} \ldots
<_{\L} m_0]$, let $f(w) = a'$.
\item[(F3)]
If $v = [2k, m_l <_{\L} \ldots
<_{\L} m_0]$ with $m_l = (a,g,b) \in U_1(k)$ and $w = [2k+1, \gamma_{(a',g,b)}
<_{\L} (a',g,b)
<_{\L} m_{l-1} <_{\L} \ldots
<_{\L} m_0]$, let $f(w) = (a',\ast)$.  
\item[(F4)]
If $v = [2k+1,m_l <_{\L} \ldots <_{\L} m_0 = I]$ with $m_l = (a,g,b)
\in U_2(k)\cup U_3(k)$ and $w = [2k+2, (a,g',b)
<_{\L} m_{l-1} <_{\L} \ldots
<_{\L} m_0]$, write
$$\epsilon((a,g',b)
<_{\L} m_{l-1} <_{\L} \ldots
<_{\L} m_0) = (x_l <_{\J} \ldots <_{\J} x_0).$$
If $x_l = (a_1,g_1,b_1)$, let $f(w) = (g_1, \ast)$.
\item[(F5)]
If $v = [2k+1,m_l <_{\L} \ldots <_{\L} m_0 = I]$ with $m_l = (a,g,b)
\in U_2(k)$ and $w = [2k+2, (1,1,b') <_{\L} (a,g',b)
<_{\L} m_{l-1} <_{\L} \ldots
<_{\L} m_0]$, write
$$\epsilon((1,1,b') <_{\L} (a,g',b)
<_{\L} m_{l-1} <_{\L} \ldots
<_{\L} m_0) = (x_{l+1} <_{\J} \ldots <_{\J} x_0).$$
If $x_l = (a_1,g_1,b_1)$ and $x_{l+1} = (a_2,g_2,b_2)$, let $f(w) =
(g_1, b_2)$.
\item[(F6)]
If $v = [2k+1,m_l <_{\L} \ldots <_{\L} m_0 = I]$ with $m_{l-1} = (a,g,b)
\in U_4(k)$ and $w = [2k+2, (1,1,b') <_{\L} (a,g',b)
<_{\L} m_{l-2} <_{\L} \ldots
<_{\L} m_0]$, write
$$\epsilon((1,1,b') <_{\L} (a,g',b)
<_{\L} m_{l-2} <_{\L} \ldots
<_{\L} m_0) = (x_{l} <_{\J} \ldots <_{\J} x_0).$$
If $x_{l-1} = (a_1,g_1,b_1)$ and $x_{l} = (a_2,g_2,b_2)$, let 
$f(w) = (g_1,b_2)$.
\ei

Note that $w = [i,\sigma] \Rw f(w) \in X_i$ in all cases: this holds
trivially if $|\son(v)| = 1$. If $i = 2k+1$ and $v = [2k, m_l <_{\L} \ldots
<_{\L} m_0]$ with $m_l \in U_0(k)$, 
then $f(w) \in A_{m_l} \in X_{2k+1} = X_i$; if $m_l \in U_1(k)$, 
then $f(w) \in A_{m_l} \cup (A'_{m_l} \times \{ \ast \}) \subseteq
X_{2k+1} = X_i$.  
 
Finally, assume that $i = 2k+2$
and $v
= [2k+1,m_l <_{\L} \ldots <_{\L} m_0 = I]$ with $m_l \in U_2(k)$. If
$\epsilon(\sigma) = [x_l <_{\J} \ldots <_{\J} x_0]$, then $f(w) \in
G_{x_l} \times \{ \ast \} = G_{m_l} \times \{ \ast \}$ by Lemma
\ref{eej}(i). Thus 
$f(w) \in X_{2k+2} =
X_i$. 

Assume now that 
$\epsilon(\sigma) = [x_{l+1} <_{\J} \ldots
<_{\J} x_0 = I]$, 
$x_l = (a_1,g_1,b_1)$ and $x_{l+1} = (a_2,g_2,b_2)$.
Then $f(w) = (g_1,b_2)$. Clearly, $g_1 \in G_{x_l} = G_{m_l}$ by Lemma
\ref{eej}(i). We show that $b_2 \in Q_{m_lm_l^*}$.
By Lemma \ref{norfB}, we may assume that $$\sigma = ((1,1,b') <_{\L}
m'_l <_{\L} m_{l-1} 
<_{\L} \ldots <_{\L} m_0)$$ with $m'_l \in\; \H_{m_l}$ and $b' \in
Q_{m'_l}$. Hence $(1,1,b') \,\L\, ym'_l <_{\L} m'_l$ for some $y \in
Y$ and so $(1,1,b')(m'_l)^* \,\L\, ym'_l(m'_l)^* =
ym_lm_l^*$ by Lemma \ref{propss}(ii). Thus
$(1,1,b_2) \,\L\, x_{l+1} = (1,1,b')(m'_l)^* \,\L\,ym_lm_l^*$. Since
$m'_l\,\R\, m_lm_l^*$, 
$ym'_l <_{\L} m'_l$ implies $ym_lm_l^* <_{\L} m_lm_l^*$ by Lemma
\ref{Rcompat} and so $b_2 \in Q_{m_lm_l^*}$. Thus 
$f(w) \in X_{2k+2} = X_i$ as claimed.

The discussion of the cases arising from $U_3(k)$ and $U_4(k)$ is
analogous and can be omitted.

Clearly, for all $\sigma \in \E_Y(M^I)$ and $v \in \vert(T)$,
the elliptic action $\theta$ induces a  
mapping
$$\begin{array}{rcl}
\theta_{\sigma}^v:\son(v)&\to&\son(v\sigma)\\
w&\mapsto&w\sigma.
\end{array}$$ 

\bl
\label{ver1}
Let $v = [2k,m_l <_{\L} \ldots <_{\L} m_0 = I]$ with $m_l \in U_1(k)$ and
let $\sigma = (m'_p <_{\L} \ldots <_{\L} m'_0 = I) \in \E_Y(M^I)$.
Then
\bi
\item[(i)] $f|_{\son(v)}$ is one-to-one;
\item[(ii)] $f(\son(v)) = A_{m_l} \cup (A'_{m_l} \times \{ \ast \})$;
\item[(iii)] $|(\son(v))\sigma| > 1$ if and only if $m_lm'_p \,\R\,
  m_l$; in this case
  $f(w\sigma) = f(w)$ 
  for every $w \in \son(v)$ and $\theta_{\sigma}^v$ is a permutation;
\item[(iv)] $|(\son(v))\sigma| = 1$ if and only if $m_lm'_p <_{\J}
  m_l$; in this case $\theta_{\sigma}^v$ is constant.
\ei
\el

\proof
Writing $m_l =
(a,g,b)$, then
$$\begin{array}{lll}
\son(v)&=&\{ [2k+1, (a',g,b) <_{\L} m_{l-1}
<_{\L} \ldots <_{\L} m_0];\; a' \in A_{m_l} \}\\
&\cup&\{ [2k+1, \gamma_{(a',g,b)} <_{\L} (a',g,b) <_{\L} m_{l-1}
<_{\L} \ldots <_{\L} m_0];\; a' \in A'_{m_l} \}
\end{array}$$ by Lemma \ref{norfA}
and these elements are 
all distinct. Since 
$$f([2k+1, (a',g,b) <_{\L} m_{l-1}
<_{\L} \ldots <_{\L} m_0]) = a',$$
$$f([2k+1, \gamma_{(a',g,b)} <_{\L} (a',g,b) <_{\L} m_{l-1}
<_{\L} \ldots <_{\L} m_0]) = (a',\ast),$$ (i) and (ii) follow.

We may write
$$(m_l <_{\L} \ldots <_{\L} m_0)\sigma = (m_lm'_p <_{\L} n_t
<_{\L} \ldots <_{\L} n_0)$$ for some $n_0, \ldots,n_t \in M^I$. 
Since $(a',g,b)\,\L\,
m_l$, we get $(a',g,b)m'_p\,\L\,
m_lm'_p$ and so
$$[2k+1,(a',g,b) <_{\L} \ldots <_{\L} m_0]\sigma = [2k+1,(a',g,b)m'_p <_{\L} n_t
<_{\L} \ldots <_{\L} n_0].$$
Writing $\zeta = (\gamma_{(a',g,b)} <_{\L} (a',g,b) <_{\L} \ldots <_{\L}
m_0)$, we get also
$$[2k+1,\zeta]\sigma = [2k+1,\lm(\gamma_{(a',g,b)}m'_p \leq_{\L}
(a',g,b)m'_p <_{\L} n_t  
<_{\L} \ldots <_{\L} n_0)].$$

Suppose that $m_lm'_p \not\R m_l$. Since $m_lm'_p \leq_{\J} m_l$, it
follows from (\ref{stbl}) that $m_lm'_p <_{\J} m_l$ and so
$h_{\J}(m_lm'_p) > k$. Then $(a',g,b)m'_p \L m_lm'_p$ yields 
$h_{\J}((a',g,b)m'_p) > k$ and it follows easily that $|(\son(v))\sigma| =
1$.

Conversely, assume that $m_lm'_p \,\R\, m_l$. Since $(a',g,b)\,\L\,
m_l$, we get  
$(a',g,b)m'_p\,\R\,
(a',g,b)$ and so 
$$f( [2k+1,(a',g,b)
<_{\L} \ldots <_{\L} m_0]\sigma) = a' = 
f( [2k+1,(a',g,b) <_{\L} \ldots <_{\L} m_0]).$$ 
Moreover, if $a' \in A'_{m_l}$ and $w = [2k+1,\zeta]$, Lemma \ref{gedi} yields
$$w\sigma = [2k+1,\gamma_{(a',g,b)}m'_p <_{\L} (a',g,b)m'_p <_{\L} n_t 
<_{\L} \ldots <_{\L} n_0].$$
Assume that $m_lm'_p = (a,g',b')$ so that $v\sigma = [2k,
(a,g',b') <_{\L} n_t <_{\L} \ldots <_{\L} n_0)]$. Since $w\sigma \in
\son(v\sigma)$, it follows from Lemma \ref{norfA} that
\beq
\label{impie1}
w\sigma = [2k+1,(a',g',b') <_{\L} n_t <_{\L} \ldots <_{\L}
n_0]\mbox{ for some }a' \in A_{(a,g',b')}
\eeq
or
\beq
\label{impie2}
w\sigma = [2k+1,\gamma_{(a',g',b')} <_{\L} (a',g',b') <_{\L} n_t
<_{\L} \ldots <_{\L} 
n_0]\mbox{ for some }a' \in A'_{(a,g',b')}.
\eeq
If (\ref{impie1}) holds, then 
$$H_Y(\gamma_{(a',g,b)}m'_p <_{\L} (a',g,b)m'_p <_{\L} n_t 
<_{\L} \ldots <_{\L} n_0, (a',g',b') <_{\L} n_t <_{\L} \ldots <_{\L}
n_0) \geq 2k+1$$
and so this pair belongs to $V(M^I)$, yielding $(a',g',b') \in W(M^I)$
by Lemma \ref{charV}. Since
$$(a',g',b') \L (a,g',b') = m_lm'_p \R m_l \notin W(M^I),$$
this contradicts Lemma \ref{WJclos}. Hence
(\ref{impie2}) holds and so 
$f(w\sigma) = (a',\ast) = f(w)$.
Thus
$|(\son(v))\sigma| > 
1$ and also $f(w\sigma) = f(w)$
  for every $w \in \son(v)$. Since $A_{m_l} = A_{m_lm'_p}$ and
  $A'_{m_l} = A'_{m_lm'_p}$ by (\ref{gedi3}), we have a
  commutative diagram
$$\xymatrix{
\son(v) \ar@{->}^{\theta_{\sigma}^v}[rrrr]
\ar@{->}_{f_1}[ddrr]
&&&&\son(v\sigma) \ar@{->}^{f_2}[ddll]\\
&&&&\\
&& A_{m_l} \cup (A'_{m_l} \times \{ \ast \}) &&
}$$
where $f_1$ and $f_2$ are the corresponding restrictions of $f$. Since $f_1$ and
$f_2$ are bijective by (i) and (ii), 
$\theta_{\sigma}^v$ must be bijective as well. Thus (iii) holds.

We have $m_lm'_p \leq_{\J}
  m_l$. By (iii) and (S1), $|(\son(v))\sigma| = 1$ if and only if
  $m_lm'_p \,\not\J\, 
  m_l$ and therefore $m_lm'_p <_{\J}
  m_l$. It is straightforward to check that $\theta_{\sigma}^v$ is constant.
\qed

The proof of the  following lemma is a simplification of the preceding
one and is therefore omitted.

\bl
\label{ver1A}
Let $v = [2k,m_l <_{\L} \ldots <_{\L} m_0 = I]$ with $m_l \in U_0(k)$ and
let $\sigma = (m'_p <_{\L} \ldots <_{\L} m'_0 = I) \in \E_Y(M^I)$.
Then
\bi
\item[(i)] $f|_{\son(v)}$ is one-to-one;
\item[(ii)] $f(\son(v)) = A_{m_l}$;
\item[(iii)] $|(\son(v))\sigma| > 1$ if and only if $m_lm'_p \,\R\,
  m_l$; in this case
  $f(w\sigma) = f(w)$ 
  for every $w \in \son(v)$ and $\theta_{\sigma}^v$ is a permutation;
\item[(iv)] $|(\son(v))\sigma| = 1$ if and only if $m_lm'_p <_{\J}
  m_l$; in this case $\theta_{\sigma}^v$ is constant.
\ei
\el

\bl
\label{ver2}
Let $v = [2k+1,m_l <_{\L} \ldots <_{\L} m_0 = I]$ with $m_l \in U_2(k)$ and
let $\zeta = (m'_p <_{\L} \ldots <_{\L} m'_0 = I) \in \E(M^I)$.
Then
\bi
\item[(i)] $f|_{\son(v)}$ is one-to-one;
\item[(ii)] $f(\son(v)) = G_{m_l} \times
(\{ \ast \} \cup Q_{m_lm_l^*})$;
\item[(iii)] $|(\son(v))\zeta| > 1$ if and only if $m_lm'_p \,\R\,
  m_l$; in this case $f(\son(v\zeta)) = f(\son(v))$
  and $\theta_{\zeta}^v$ is a permutation;
\item[(iv)] $|(\son(v))\zeta| = 1$ if and only if $m_lm'_p <_{\J}
  m_l$; in this case $\theta_{\zeta}^v$ is constant.
\ei
\el

\proof
Writing $m_l = 
(a,g,b)$, it follows from Lemma \ref{norfB} that $\son(v) = \son_1(v)
\cup \son_2(v)$ with
$$\begin{array}{lll}
\son_1(v)&=&\{ [2k+2, r <_{\L} m_{l-1}
<_{\L} \ldots <_{\L} m_0];\; r \in \; \H_{m_l} \},\\
\son_2(v)&=&\{ [2k+2, (1,1,b') <_{\L} r <_{\L} m_{l-1}
<_{\L} \ldots <_{\L} m_0];\; r \in \; \H_{m_l},\; b' \in Q_{r}\}
\end{array}$$
and these elements are
all distinct. 

Let $$\sigma = (r <_{\L} m_{l-1}
<_{\L} \ldots <_{\L} m_0) \in \son_1(v).$$
If $\epsilon(\sigma) = (x_l <_{\J} \ldots <_{\J} x_0 = I)$ and $x_l =
(a_1,g_1,b_1)$, then $f([2k+2,\sigma]) = (g_1,\ast) \in G_{x_l} \times
\{ \ast \}$. Note that $x_l 
\R r \H m_l$ by 
Lemma \ref{eej}(i) and so $f([2k+2,\sigma]) \in G_{m_l} \times \{ \ast \}
\subseteq X_{2k+2}$. By Green's Lemma, the mapping
$$\begin{array}{rcl}
\H_{m_l}&\to&\H_{m_lm_{l-1}^*}\\
r&\mapsto&rm_{l-1}^*
\end{array}$$
is a bijection and so $f|_{\son_1(v)}$ is one-to-one and
\beq
\label{bay1} 
f(\son_1(v)) = G_{m_l} \times \{ \ast \}.
\eeq 

Next let
$$\tau = ((1,1,b') <_{\L} r <_{\L} m_{l-1}
<_{\L} \ldots <_{\L} m_0) \in \son_{2}(v)$$
with $r \in \; \H_{m_l}$ and $b' \in Q_{r}$. 
If $\epsilon(\tau) = (x_{l+1} <_{\J} \ldots <_{\J} x_0 = I)$, $x_l =
(a_1,g_1,b_1)$ and $x_{l+1} = (1,1,b')r^* = (a_2,g_2,b_2)$, then
$f([2k+2,\tau]) = 
(g_1,b_2)$. 
We fix $r \in\; \H_{m_l}$ and write
$$\son_{2,r}(v) = \{ [2k+2, (1,1,b') <_{\L} r <_{\L} m_{l-1}
<_{\L} \ldots <_{\L} m_0];\; b' \in Q_{r} \}.$$
In view of the preceding case, to complete the proof of (i)
and (ii) it suffices
to show that $f|_{\son_{2,r}(v)}$ is one-to-one and
\beq
\label{bay2}
f(\son_{2,r}(v)) = \{ g_1 \} \times Q_{m_lm_l^*}.
\eeq 

Indeed, since $rr^* = m_lm_l^*$ 
by Proposition \ref{propss}(ii), the mapping
$$\begin{array}{rcl}
\p:Q_{r}&\to&Q_{m_lm_l^*}\\
b&\mapsto&((1,1,b)r^*)\pi_3
\end{array}$$
is a bijection by Lemma \ref{Rcompat}. 
Thus (\ref{bay2}) holds and so
\beq
\label{bay3} 
f(\son_2(v)) = G_{m_l} \times Q_{m_lm_l^*}.
\eeq 

In view of (\ref{bay1}), (\ref{bay2}) and the partial injectivity
results obtained, (i) and (ii) hold.

Assume now that $|(\son(v))\zeta| >
1$. Suppose that $h_{\J}(m_lm'_p) > k$. Then
$h_{\J}(rm'_p)$\linebreak
$>
k$ for every $r \in \;\H_{m_l}$ due to $rm'_p \,\L\,
m_lm'_p$. Since the $rm'_p$ would then be all $\L$-equivalent,
we would get 
$|(\son(v))\zeta| = 1$, a contradiction. Thus
$h_{\J}(m_lm'_p) = k = h_{\J}(m_l)$ and so $m_lm'_p \,\J\, m_l$. By
(S1), we get $m_lm'_p \,\R\, m_l$.

Conversely, assume that $m_lm'_p \,\R\, m_l$. Then $m_l = m_lm'_pz$
for some $z \in M^I$. 
Taking a minimal representation
$$v\zeta = [2k+1, m_lm'_p <_{\L} n_t
<_{\L} \ldots <_{\L} n_0]$$ for some $n_0,
\ldots,n_t \in M^I$,
it follows easily from $m_l = m_lm'_pz$ that the elements
$[2k+2, rm'_p <_{\L} n_t
<_{\L} \ldots <_{\L} n_0]$ and 
$[2k+2, (1,1,b')m'_p <_{\L} rm'_p <_{\L} n_t
<_{\L} \ldots <_{\L} n_0]$
of $(\son(v))\zeta$ are all distinct, hence $|(\son(v))\zeta| =
|\son(v)| > 1$.
Moreover, applying (i) and (ii) to $v$ and $v\zeta$, we have 
$$|\son(v)| = |G_{m_l}| \cdot (1 + |Q_{m_lm_l^*}|),$$
$$|\son(v\zeta)| = |G_{m_lm'_p}| \cdot (1 +
|Q_{m_lm'_p(m_lm'_p)^*}|).$$ Since $m_lm'_p \,\R\, m_l$, we get
$G_{m_lm'_p} = G_{m_l}$  and also
$m_lm'_p(m_lm'_p)^* = m_lm_l^*$ by Proposition \ref{propss}(ii). Thus
$|\son(v\zeta)| = |\son(v)| > 1$.
Still applying (i) and (ii) to $v$ and 
$v\zeta$, we get $$f(\son(v\zeta)) = G_{m_l} \times
(\{ \ast \} 
\cup Q_{m_lm_l^*}) = f(\son(v)).$$ Furthermore, we
have a
  commutative diagram
$$\xymatrix{
\son(v) \ar@{->}^{\theta_{\zeta}^v}[rrrr]
\ar@{->}_{f_1}[ddrr]
&&&&\son(v\zeta) \ar@{->}^{f_2}[ddll]\\
&&&&\\
&& G_{m_l} \times
(\{ \ast \} 
\cup Q_{m_lm_l^*})&&
}$$
where $f_1$ and $f_2$ are the corresponding restrictions of $f$. Since $f_1$ and
$f_2$ are bijective by (i) and (ii), 
$\theta_{\zeta}^v$ must be bijective as well.

The proof of (iv) is analogous to the proof of Lemma \ref{ver1}(iv).
\qed

The proofs of the following two lemmas constitute straightforward
adaptations of the proof of Lemma \ref{ver2} and can therefore be omitted.

\bl
\label{ver3}
Let $v = [2k+1,m_l <_{\L} \ldots <_{\L} m_0 = I]$ with $m_l \in U_3(k)$ and
let $\zeta = (m'_p <_{\L} \ldots <_{\L} m'_0 = I) \in \E(M^I)$.
Then
\bi
\item[(i)] $f|_{\son(v)}$ is one-to-one;
\item[(ii)] $f(\son(v)) = G_{m_l} \times \{ \ast \}$;
\item[(iii)] $|(\son(v))\zeta| > 1$ if and only if $m_lm'_p \,\R\,
  m_l$; in this case $f(\son(v\zeta)) = f(\son(v))$
  and $\theta_{\zeta}^v$ is a permutation;
\item[(iv)] $|(\son(v))\zeta| = 1$ if and only if $m_lm'_p <_{\J}
  m_l$; in this case $\theta_{\zeta}^v$ is constant.
\ei
\el

\bl
\label{ver4}
Let $v = [2k+1,m_l <_{\L} \ldots <_{\L} m_0 = I]$ with $m_{l-1} \in U_4(k)$ and
let $\zeta = (m'_p <_{\L} \ldots <_{\L} m'_0 = I) \in \E(M^I)$.
Then
\bi
\item[(i)] $f|_{\son(v)}$ is one-to-one;
\item[(ii)] $f(\son(v)) = G_{m_{l-1}} \times Q_{m_{l-1}m_{l-1}^*}$;
\item[(iii)] $|(\son(v))\zeta| > 1$ if and only if $m_{l-1}m'_p \,\R\,
  m_{l-1}$; in this case $f(\son(v\zeta)) = f(\son(v))$
  and $\theta_{\zeta}^v$ is a permutation;
\item[(iv)] $|(\son(v))\zeta| = 1$ if and only if $m_{l-1}m'_p <_{\J}
  m_{l-1}$; in this case $\theta_{\zeta}^v$ is constant.
\ei
\el

Given a set $X$, we write
\bi
\item[] $S(X) = \{ \p \in M(X) : \p$ is a permutation of $X\}$.
\item[] $K(X) = \{ \p \in P(X) : |X\p| \leq 1 \}$.
\ei
It is immediate that
both $S(X) \cup K(X)$ and $\{ \id_X \} \cup K(X)$ constitute submonoids
of $P(X)$.

In the main result of the paper, we construct an embedding 
$$\begin{array}{rcl}
\p:\E(M^I)&\to&\Pi_{i=1}^{\delta}
(X_i,M_i) = \ldots \circ (X_2,M_2) \circ (X_1,M_1)\\
\sigma&\mapsto&\p_{\sigma}
\end{array}$$
into an iterated
wreath product of partial transformation semigroups
where $M_{2k+1}$ is a submonoid of $\{ \id_{X_{2k+1}}\} \cup K(X_{2k+1})$ and
$M_{2k+2}$ is a submonoid of $S(X_{2k+2}) \cup K(X_{2k+2})$.
Furthermore, we shall prove that this embedding has the {\em Zeiger
  property}: 
if $$(\cdot,x_{2k+1}, \ldots,x_1)\p_{\sigma}\pi_{2k+2} \in
S(X_{2k+2})\setminus K(X_{2k+2}),$$ then any local mapping of the form
$(\cdot,x_{q-1}, \ldots,x_1)\p_{\sigma}\pi_{q}$ for $2k+2 \leq q-1 <
\delta$ must be the identity.

\bt
\label{pcm}
Let $M$ be a finite $\J$-above $Y$-semigroup and let $\delta =
2+2\sup \{ 
h_{\J}(m) \mid m \in M \} \in \oo{\N}$. Then there exists an embedding
$\p$ of $\E_Y(M^I)$ into the iterated
wreath product of partial transformation semigroups
$\Pi_{i=1}^{\delta}
(X_i,M_i) = \ldots \circ (X_2,M_2) \circ (X_1,M_1)$
such that:
\bi
\item[(i)] $M_{2k+1}$ is a submonoid of $\{ \id_{X_{2k+1}}\} \cup
  K(X_{2k+1})$ for $2k+1 \leq \delta$.
\item[(ii)] $M_{2k+2}$ is a submonoid of $S(X_{2k+2}) \cup
  K(X_{2k+2})$  for $2k+2 \leq \delta$; 
if $\{ R_{\lambda} \mid \lambda \in \Lambda \}$ is the set of all
  $\R$-classes of $M$ contained in $U_2(k) \cup U_3(k) \cup U_4(k)$, then
\beq
\label{reg1}
M_{2k+2} \cap S(X_{2k+2}) \cong \oplus_{\lambda \in \Lambda} G'_{\lambda},
\eeq 
where $G'_{\lambda}$ is a subgroup of $G_{\lambda}$.
\item[(iii)] $\p$ has the Zeiger property.
\ei
Moreover, if $Y$ is finite, then the
$X_i$ (and consequently the 
$M_i$) are all finite.
\et

\proof
For commodity, we assume that $\delta \in \N$, the infinite case being
absolutely similar. 

We consider the length function $H_Y$ and we assume that $H_Y
= D_{\chi}$ for $\chi = (r_0,T,\alpha,\theta)$, $\chi$ being obtained
by the Chiswell construction. Let $X_i$ and $f$ be defined as before for $i =
1, \ldots, \delta$. Write $X = \bigcup_{i = 1}^{\delta} (X_i \times \ldots \times
X_1)$. By Theorem \ref{embed} and Lemmas \ref{ver1}--\ref{ver4}(i),
there exists an injective monoid homomorphism 
$$\begin{array}{rcl} 
\Psi:\E_Y(M^I)&\to&(X_{\delta},P(X_{\delta}))\circ \ldots \circ   (X_1,P(X_1))\\
\sigma&\mapsto&\Psi_{\sigma}
\end{array}$$ 
defined by
$$x\Psi_{\sigma} = \left\{
\begin{array}{ll}
x\psi\inv\theta_{\sigma}\psi&\mbox{ if } \sigma \neq (I)\hspace{1.5cm}
(x \in X).\\ 
x&\mbox{ if } \sigma = (I),
\end{array}
\right.$$
where
$$\begin{array}{rcl} 
\psi:\ray(r_0,T) &\to&X\\
(v_i, \ldots, v_1,r_0)&\mapsto&(f(v_i), \ldots, f(v_1)).
\end{array}$$ 

Given $\sigma \in \E_Y(M^I) \setminus \{ I\}$, we extend
$\Psi_{\sigma}$ to a mapping $\p_{\sigma} \in P(X)$ by taking
$$\begin{array}{ll}
\dom\p_{\sigma} = \im\psi \cup(\cup_{i = 1}^{\delta} \{& (x_i,f(v_{i-1}),
\ldots, f(v_1)) 
\in X: (v_{i-1}, \ldots, v_1) \in \ray(r_0,T)\\
&\mbox{and }|(\son(v_{i-1}))\sigma| > 1\})
\end{array}$$
and
$$(x_i,f(v_{i-1}), \ldots, f(v_1))\p_{\sigma} =
(x_i,(f(v_{i-1}), \ldots, f(v_1))\Psi_{\sigma})$$
if $(x_i,f(v_{i-1}), \ldots, f(v_1)) \notin \im\psi$.
Since $\psi$ is one-to-one, $\p_{\sigma}$ is well-defined. Being an extension of
$\Psi_{\sigma}$, it is easy to see that $\p_{\sigma}$ inherits some of
its properties, namely being sequential. Moreover, it follows from
Lemmas \ref{ver1}--\ref{ver4}(iii) that 
\beq
\label{bekk1}
(\dom\p_{\sigma} \setminus \im\psi)\p_{\sigma} \;\cap \; \im\psi =
\emptyset.
\eeq

Taking $\p_{I} = \Psi_{I} = \id_X$, we define
$$\begin{array}{rcl} 
\p:\E_Y(M^I)&\to&P(X)\\
\sigma&\mapsto&\p_{\sigma}.
\end{array}$$ 
We show that $\p$ is a monoid homomorphism. 

Since $\p_{I}$ is the identity and $\dom\p_{\sigma\tau} \subseteq
\dom\p_{\sigma}$, we only 
have to take $\sigma,\tau \in \E_Y(M^I) \setminus \{ I \}$ and show
that 
\beq
\label{pcm3}
x\p_{\sigma}\p_{\tau} = x\p_{\sigma\tau}
\eeq
holds for every $x \in \dom\p_{\sigma}$. Since $\Psi_{\sigma} \subseteq
\p_{\sigma}$ is a homomorphism, (\ref{pcm3}) holds for $x \in
\im\psi$. Assume now that 
$$x = (x_i,f(v_{i-1}), \ldots, f(v_1)) \in \dom\p_{\sigma} \setminus
\im\psi.$$ Hence $|(\son(v_{i-1}))\sigma| > 1$.
Write $(f(v_{i-1}), \ldots, f(v_1))\Psi_{\sigma} = (f(v'_{i-1}),
\ldots, f(v'_1))$. In particular, $v_{i-1}\sigma = v'_{i-1}$. 

Assume
first that $|(\son(v'_{i-1}))\tau| \leq 
1$. Then $v_{i-1}\sigma = v'_{i-1}$ yields $(\son(v_{i-1}))\sigma
\subseteq \son(v'_{i-1})$ since the action is elliptical and so
$|(\son(v_{i-1}))\sigma\tau| \leq 
1$ as well. Thus $x \notin \dom\p_{\sigma\tau}$. On the other hand,
$x\p_{\sigma} = (x, f(v'_{i-1}),
\ldots, f(v'_1)) \notin \im\psi$ by (\ref{bekk1}) and so $x\p_{\sigma}
\notin \dom\p_{\tau}$. Thus (\ref{pcm3}) holds in this case.

Finally, assume
that $|(\son(v'_{i-1}))\tau| > 
1$. Write $$(f(v_{i-1}), \ldots, f(v_1))\Psi_{\sigma}\Psi_{\tau} = (f(v'_{i-1}),
\ldots, f(v'_1))\Psi_{\tau} = (f(v''_{i-1}),
\ldots, f(v''_1)).$$
Then 
$$x\p_{\sigma}\p_{\tau} = (x, f(v'_{i-1}),
\ldots, f(v'_1))\p_{\tau} = (x, f(v''_{i-1}),
\ldots, f(v''_1))$$
by (\ref{bekk1}). On the other hand, in view of Lemmas
\ref{ver1}--\ref{ver4}(iii), $|(\son(v_{i-1}))\sigma| > 1$ and 
$|(\son(v'_{i-1}))\tau| > 1$ together yield
$|(\son(v_{i-1}))\sigma\tau| > 1$. Since $x \notin \im\psi$ and $\Psi$
is a homomorphism, we obtain
$$x\p_{\sigma\tau} = (x,(f(v_{i-1}),
\ldots, f(v_1))\Psi_{\sigma\tau}) = (x,f(v''_{i-1}),
\ldots, f(v''_1)) = x\p_{\sigma}\p_{\tau}$$
and so (\ref{pcm3}) holds as well in this case. Thus $\p$ is a monoid
homomorphism. 

We show next that $\p$ is one-to-one. Given distinct $\sigma,\tau
\in \E_Y(M^I) \setminus \{I\}$, we have $\Psi_{\sigma} \neq \Psi_{\tau}$ by
Theorem
\ref{embed}. Since $\dom\Psi_{\sigma} = \im\psi = \dom\Psi_{\tau}$, it follows
that $\p_{\sigma} \neq \p_{\tau}$ as well. To show that $\p_{\sigma}
\neq \p_{I}$, it suffices now to show that $\Psi_{\sigma}$ is not 
one-to-one. Indeed, using the Chiswell construction and by Lemma
\ref{norfB}, we have 
$|\son(v)| > 1$ for $v = [1,I]$ since $h_{\J}(I) = 0$ and $Q_{I} \neq
\emptyset$. However, for $\sigma = (n_p <_{\L} \ldots <_{\L} n_0)$
with $p > 0$, we have $In_p = n_p <_{\J} I$ and so
$\theta_{\sigma}^v$ is constant by Lemma \ref{ver2}(iv). Thus
$\Psi_{\sigma}$ is not  
one-to-one and so $\p$ is indeed one-to-one.

We proceed now to discuss the local mappings.
Let $(x_{i-1}, \ldots, x_1) \in X_{i-1} \times \ldots \times X_1$ and
write $\xi = (\cdot, x_{i-1}, \ldots, x_1)\p_{\sigma} \in P(X_i)$. We
assume $\sigma \neq (I)$. 
Assume that $\xi \notin K(X_i)$. In particular, $\xi$ is not the empty
map and so $(x_{i-1}, \ldots, x_1) = (f(v_{i-1}), \ldots, f(v_1))$ for
some $(v_{i-1}, \ldots, v_1) \in \ray(r_0,T)$. Let $\xi' =
\xi|_{\im\psi}$. It follows from the definition of $\p_{\sigma}$ that
$|\son(v_{i-1})\sigma| > 1$, otherwise $\xi = \xi' \in K(X_i)$. By
Lemmas \ref{ver1}--\ref{ver4}, it follows that 
$\xi' \in S(X'_i)$ for some $X'_i \subset X_i$ and so $\xi \in
S(X_i)$ by definition of $\p_{\sigma}$.

If $i$ is odd, then $\xi = \id_{X_i}$ by Lemmas \ref{ver1}(iii) and
\ref{ver1A}(iii), thus
we can take $M_i$ to be a submonoid of $\{ \id_{X_i}\} \cup
K(X_i)$ and (i) holds.

Assume now that $i = 2k+2$ is even. We can take $M_i$ to be the submonoid of
$S(X_i)  \cup
K(X_i)$ generated by the local mappings $\xi$. Write $\sigma = (m'_p
<_{\L} \ldots  <_{\L} m'_0 = I)$ 
with $v_{i-1} = [2k+1,m_l
<_{\L} \ldots  <_{\L} m_0]$ in minimal representation. Since $\xi
\notin K(X_i)$, then $|\son(v_{i-1})| \; > 1$ and so, by Lemma \ref{norfB},
either $m_l \in U_2(k) \cup U_3(k)$ or $m_{l-1}
  \in U_4(k)$. 

We consider first the case $m_l \in U_3(k)$.
By Lemma \ref{ver3}(iii), $\xi'$ 
permutes $G_{m_l} \times \{ \ast \}$. 
We show that there exists some $g_0 \in
G_{m_l}$ such that 
\beq
\label{bja1}
(h,\ast)\xi = (hg_0,\ast) \mbox{ for every } h \in G_{m_l}.
\eeq
Indeed, by Lemma \ref{eej}(i) we may write $m_l = (a,g,b)$ and $x_l =
(a,g_1,b_1)$. Write also
\beq
\label{bja2}
(m_l <_{\L} \ldots  <_{\L} m_0)\sigma = (m_lm'_p <_{\L} n_t <_{\L}
\ldots <_{\L} n_0). 
\eeq
Given $h \in G_{m_l}$, take $r = (a,h,b_1) m_{l-1}^{\#}$, 
$\tau = (r <_{\L} m_{l-1} <_{\L} \ldots <_{\L} m_0)$
and $w = [2k+2,\tau]$.
We claim that 
\beq
\label{rust}
w \in \son(v_{i-1}) \hspace{.7cm}\mbox{ and }\hspace{.7cm}
f(w) = (h,\ast).
\eeq
Indeed,
$(a,h,b_1) \H x_l$ yields $r \L x_lm_{l-1}^{\#} = m_l$ by Lemma
\ref{eej}(iv). On the other hand, $(a,h,b_1) m_{l-1}^{\#} \L m_l \R
x_l \H (a,h,b_1)$ 
yields $r \R (a,h,b_1)$ by (S1) and so $r \R x_l \R m_l$. Thus $r \in
\; \H_{m_l}$. Moreover, $$rm_{l-1}^* = (a,h,b_1) m_{l-1}^{\#}m_{l-1}^*
= (a,h,b_1)$$ since $(a,h,b_1) \L x_l$ and $x_lm_{l-1}^{\#}m_{l-1}^* =
m_lm_{l-1}^* = x_l$ by Lemma
\ref{eej}(iv). Thus (\ref{rust}) holds.
Now, since
$m_lm'_p \,\L\, rm'_p$, it follows from (\ref{bja2}) that  
$\tau\sigma = (rm'_p <_{\L}
n_t <_{\L} \ldots <_{\L} n_0)$  
and so 
\beq
\label{vli1}
\epsilon(\tau\sigma) =
(rm'_pn_t^* <_{\J} \ldots ).
\eeq  
Thus
$$(h,\ast)\xi = (  (rm'_pn_t^*)\pi_2,\ast) = ( ((a,h,b_1)
m_{l-1}^{\#} m'_pn_t^*)\pi_2,\ast).$$ 
Let $y = m_{l-1}^{\#} m'_pn_t^*$. Since
$$(a,g_1,b_1)m_{l-1}^{\#} m'_pn_t^* = x_lm_{l-1}^{\#}
m'_pn_t^* = m_lm'_pn_t^*$$ and 
$m_lm'_pn_t^* \R m_lm'_p \R m_l \R x_l = (a,g_1,b_1)$ by Lemma \ref{eej}(i) and
(\ref{vli1}), it follows from Proposition \ref{schutz} that there
exists some $g_0 \in G_{m_l}$ such that
$$\forall h \in G_{m_l},\; ((a,h,b_1)y)\pi_2 = hg_0.$$
Thus 
(\ref{bja1}) holds.

We consider next the case $m_{l-1} \in U_4(k)$.
By Lemma \ref{ver4}(iii), $\xi'$
permutes $G_{m_{l-1}} \times Q_{m_{l-1}m_{l-1}^*}$. 
We show that there exists some $g_0 \in G_{m_{l-1}}$ such that
\beq
\label{bja3}
(h,c)\xi' = (hg_0,c) \mbox{ for all } h \in G_{m_{l-1}} \mbox{ and } c \in
Q_{m_{l-1}m_{l-1}^*}. 
\eeq
Since $m_{l-1}m'_p \R m_{l-1}$ by Lemma \ref{ver4}(iii), $m_l <_{\L}
m_{l-1}$ yields $m_{l}m'_p \R m_{l}$ by Lemma \ref{gedi} and so we may
assume that  
\beq
\label{bja2new}
(m_l <_{\L} \ldots  <_{\L} m_0)\sigma = (m_lm'_p <_{\L} m_{l-1}m'_p
<_{\L} n_t <_{\L}
\ldots <_{\L} n_0). 
\eeq
Let $h \in G_{m_{l-1}}$ and $c \in
Q_{m_{l-1}m_{l-1}^*}$. Let $r = (a,h,b_1)m_{l-2}^{\sharp}$ and
  $(1,1,b') \L (1,1,c)r^{\sharp}$. 
Let
$$\tau = ((1,1,b') <_{\L} r <_{\L} m_{l-2} <_{\L} \ldots <_{\L} m_0)$$
and $w = [2k+2,\tau]$. We claim that 
\beq
\label{rust1}
w \in \son(v_{i-1}) \hspace{.7cm}\mbox{ and }\hspace{.7cm}
f(w) = (h,c).
\eeq
Indeed, the proof of (\ref{rust}) can be easily adapted to show that
$r \H m_{l-1}$ and $f(w) = (h,\ldots)$ (if indeed $w \in \son(v_{i-1})$). 
Since $c \in Q_{m_{l-1}m_{l-1}^*}$, we have $(1,1,c) \L
ym_{l-1}m_{l-1}^*$ for some $y \in Y_{m_{l-1}m_{l-1}^*}$. Hence
$$(1,1,b') \L (1,1,c)r^{\sharp} \L ym_{l-1}m_{l-1}^*r^{\sharp} =
yrr^*r^{\sharp} = yr.$$ Since $yr \L r$ would imply $ym_{l-1}m_{l-1}^*
\L m_{l-1}m_{l-1}^*$ in view of $r \R m_{l-1}m_{l-1}^*$, contradicting
 $y \in Y_{m_{l-1}m_{l-1}^*}$, we get $yr <_{\L} r$ and so $y \in
 Y_r$. Thus $b' \in Q_r$ and so $w \in \son(v_{i-1})$ by Lemma \ref{norfB}. Now
$$(1,1,b')r^* \L yrr^* = ym_{l-1}m_{l-1}^* \L (1,1,c),$$
hence $f(w) = (h,c)$ and so (\ref{rust1}) holds.

Now $(1,1,b') <_{\L} r \L m_{l-1}$ yields $(1,1,b') m'_p <_{\L} r
m'_p$ by Lemma \ref{gedi}. 
Similarly to the preceding case, it follows easily from
(\ref{bja2new}) that
$$\tau\sigma = ((1,1,b')m'_p <_{\L} rm'_p <_{\L} n_{t} <_{\L}
\ldots <_{\L} n_0)$$ 
and so $$\epsilon(\tau\sigma) = ((1,1,b')m'_p (rm'_p)^* <_{\J}
rm'_p n_{t}^* <_{\J} \ldots <_{\J} I).$$ 
Since $(1,1,b') <_{\L} r$, we may write $(1,1,b') = zr$ for some $z \in M$.  
Since $m_{l-1}m'_p \,\R \, m_{l-1}$ and $r \L m_{l-1}$, we get $r
m'_p \,\R \, r$ and so Lemma \ref{propss}(ii) yields  
$$(1,1,b')m'_p (rm'_p)^* = zrm'_p (rm'_p)^* =
zrr^* = (1,1,b')r^*.$$  
Hence the leftmost term in $\epsilon(\tau\sigma)$ is the same as in
$\epsilon(\tau)$ and so $(h,c)\xi' =
(\ldots,c)$. A straightforward adaptation of the proof of (\ref{bja1})
completes the proof of (\ref{bja3}).
 
Similarly, in the case $m_l \in U_2(k)$ we show that there exists some
$g_0 \in G_{m_{l}}$ such that 
\beq
\label{rust2}
(h,c)\xi' = (hg_0,c) \mbox{ for all } h \in G_{m_{l}} \mbox{ and } c \in
\{ \ast\} \cup Q_{m_{l}m_{l}^*}. 
\eeq
Indeed, by (\ref{bay1}) and (\ref{bay3}), $\xi'$ is the (disjoint)
union of a permutation $\xi'_1$ of $G_{m_l} \times \{ \ast \}$ with a
permutation $\xi'_2$ of $G_{m_l} \times Q_{m_lm_l^*}$. A
straightforward combination of the two preceding cases yields (\ref{rust2}).

Write
$$K = \left\{ 
\begin{array}{ll}
G_{m_l} \times (\{ \ast \} \cup Q_{m_lm_l^*})&\mbox{ if }m_l \in U_2(k)\\
G_{m_l} \times \{ \ast \} &\mbox{ if }m_l \in U_3(k)\\
G_{m_{l-1}} \times Q_{m_{l-1}m_{l-1}^*}&\mbox{ if }m_{l-1} \in U_4(k).
\end{array}
\right.$$
By (\ref{bja1}), (\ref{bja3}) and (\ref{rust2}),
each local map $\xi \in M_i \cap S(X_i)$ can be decomposed as a
disjoint union of permutations $\xi = \xi' \cup \xi''$ where 
$$\begin{array}{rcl}
\xi': K&\to&K\\ 
(h,c)&\mapsto&(hg_0,c),
\end{array}$$
for some $g_0 \in G_{m_l}$ ($G_{m_{l-1}}$ if $m_{l-1} \in U_4(k)$),
and $\xi''$ is the identity mapping on $X_i\setminus K$.

For every $\lambda \in \Lambda$, take $m \in R_{\lambda}$
and 
$$K_{\lambda} = \left\{ 
\begin{array}{ll}
G_{m} \times (\{ \ast \} \cup Q_{mm^*})&\mbox{ if }m \in U_2(k)\\
G_{m} \times \{ \ast \} &\mbox{ if }m \in U_3(k)\\
G_{m} \times Q_{mm^*}&\mbox{ if }m \in U_4(k).
\end{array}
\right.$$
Note that $K_{\lambda}$ is well defined in view of Lemmas \ref{WJclos},
\ref{Rcompat} and \ref{URclos}. 
Write
$$S_{\lambda}(X_{2k+2}) = \{ \p \in S(X_{2k+2}) \mid \; \p|_{ X_{2k+2}\setminus
  K_{\lambda}} = \id
\}.$$
Since the sets $K_{\lambda}$ are disjoint subsets of $X_{2k+2}$, we can view
$S_{\lambda}(X_{2k+2})$ as a direct sum of its subgroups
$S_{\lambda}(X_{2k+2})$. 
We show that
\beq
\label{espswe}
M_{2k+2} \cap S(X_{2k+2}) = \oplus_{\lambda \in \Lambda} \; (M_{2k+2} \cap
S_{\lambda}(X_{2k+2})).
\eeq
Indeed, the union 
$$\begin{array}{lll}
X_{2k+2}&=&\{ \downarrow\} \cup(\bigcup_{m \in U_2(k)} (G_{m} \times
(\{ \ast \} 
\cup Q_{mm^*})))\\
&\cup&(\bigcup_{m \in U_3(k)} (G_{m} \times
\{ \ast \})) \cup (\bigcup_{m \in U_4(k)} 
(G_{m} \times Q_{mm^*})).
\end{array}$$
is supposed to be disjoint over distinct
$\R$-classes, and the  decomposition $\xi = \xi' \cup \xi''$ shows
that every local map $\xi$ belongs indeed to a unique
$S_{\lambda}(X_{2k+2})$. Since $M_{2k+2}$ is by definition generated by the
local maps $\xi$, it follows that $M_{2k+2} \cap S(X_{2k+2}) \subseteq
\oplus_{\lambda 
  \in \Lambda} \; (M_{2k+2} \cap 
S_{\lambda}(X_{2k+2}))$. The opposite inclusion is trivial, hence
(\ref{espswe}) holds.

It follows from the decomposition $\xi = \xi' \cup \xi''$,
(\ref{bja1}), (\ref{bja3}) and (\ref{rust2}) that we can take
$M_{2k+2} \cap S_{\lambda}(X_{2k+2}) \cong G'_{\lambda}$ for some subgroup
$G'_{\lambda}$ of $G_{\lambda}$, hence  
$$M_{2k+2} \cap S(X_{2k+2}) \cong  \oplus_{\lambda \in \Lambda} G'_{\lambda}$$
and (ii) holds.  

Finally, we prove that $\p$ has the Zeiger property. Let $\sigma = (m'_p
<_{\L} \ldots <_{\L} m'_0) \in
\E_Y(M^I)$.
We may assume that $p > 0$. 
Suppose that
$\xi = (\cdot,f(v_{2k+1}), \ldots,f(v_1))\p_{\sigma}\pi_{2k+2} \in
S(X_{2k+2})\setminus K(X_{2k+2})$ and $(v_{q-1}, \ldots, v_1,r_0) \in
\ray(r_0,T)$ with $2k+1 < q-1 < \delta$. Let $\xi' = (\cdot,f(v_{q-1}),
\ldots,f(v_1))\p_{\sigma}\pi_{q}$. We show that $\xi'$ is the identity
mapping by induction on $q$. Assume the claim holds for $q'$ whenever
$2k+1 < q'-1 < q-1$.   

Let $v_{i-1} = [2k+1,\tau]$ in minimal representation, with $\tau = (m_l
<_{\L} \ldots  <_{\L} m_0)$. Since $\xi \in
S(X_{2k+2})\setminus K(X_{2k+2})$, we have either $m_l \in
U_2(k)\cup
U_3(k)$ or $m_{l-1} \in
U_4(k)$ by Lemma \ref{norfB}. Let
$$d = \left\{
\begin{array}{ll}
l&\mbox{ if } m_l \in
U_2(k)\cup
U_3(k)\\
l-1&\mbox{ if }m_{l-1} \in
U_4(k)
\end{array}
\right.$$ 

By Lemmas
\ref{ver2}--\ref{ver4}(iii), we have $m_dm'_p \R m_d$. Write $v_{q-1} = [q-1,
\rho]$ in minimal 
representation. Since $v_{q-1}$ must be a descendant of
$v_{2k+1}$, it follows from (\ref{lftree2}) that
$H(\tau \wedge_{\L} \rho) \geq 2k+1$ and so either $h_{\J}(\tau \wedge_{\L}\rho)
> k$ or $(\tau,\rho) \in V(M^I)$. Hence
\beq
\label{vli3}
\rho = (n_{l'} <_{\L} \ldots <_{\L} n_d
<_{\L} m_{d-1} <_{\L} \ldots 
<_{\L} m_0)
\eeq
for some $n_j$. By Lemma \ref{charV}, we have
$n_d = (\rho \wedge_{\L} \tau) \H (\tau \wedge_{\L} \rho) = 
m_d$. Since $m_dm'_p \R m_d$, Lemma \ref{gedi} yields
\beq
\label{vli5}
n_{j}m'_p \R n_{j}\quad (j =
d, \ldots,l'),
\eeq 
\beq
\label{grerus}
\rho\sigma = (n_{l'}m'_p <_{\L} \ldots <_{\L}
n_{d}m'_p <_{\L} \ldots).
\eeq 

Assume first that $|\son(v_{q-1})| = 1$. Then we must have $\dom\xi' = \{
\downarrow \}$. Suppose that $|\son(v_{q-1}\sigma)|\;  > 1$. Suppose
further that $l' =
d$. Then
$$v_{q-1}\sigma = [q-1, \rho\sigma] = [q-1, n_dm'_p <_{\L} \ldots].$$
Since $h_{\J}(n_dm'_p) = h_{\J}(n_d) = h_{\J}(m_d) = k$ and $2k+2 \leq
q-1$, it follows from Lemmas \ref{norfA} and \ref{norfB} that
$|\son(v_{q-1}\sigma)| = 1$, a contradiction. Hence $l' >
d$. By (\ref{vli5}), (\ref{grerus}) and Lemma \ref{ferias1},
$v_{q-1} = [q-1, \rho]$ being in minimal
representation implies that so it is
$v_{q-1}\sigma = [q-1, \rho\sigma]$. It follows that for $q$ odd
(respectively even) we have $n_{l'}m'_p \in U_0(k')\cup U_1(k')$ for $k' =
\frac{q-1}{2}$ (respectively $n_{l'}m'_p \in U_2(k')\cup U_3(k')$ or
$n_{l'-1}m'_p \in U_4(k')$ for $k' =
\frac{q-2}{2}$).

Suppose first that $n_{l'}m'_p \in U_0(k')$. Then $h_{\J}(n_{l'}m'_p)
= k'$ and $|A_{n_{l'}m'_p}| >
1$. Since $n_{l'}m'_p \R n_{l'}$ by (\ref{vli5}), we get $h_{\J}(n_{l'})
= k'$ and $|A_{n_{l'}}| >
1$ and so $n_{l'} \in U_0(k')$. By Lemma \ref{norfA}, this contradicts
$|\son(v_{q-1})| = 1$. The case $n_{l'}m'_p \in U_1(k')$ is analogous.

Assume now that 
$n_{l'}m'_p \in U_2(k')$. Then $h_{\J}(n_{l'}m'_p)
= k'$ and $|G_{n_{l'}m'_p}|(1+|Q_{n_{l'}m'_p}|) >
1$. 
By Lemma \ref{Rcompat}, we have
\beq
\label{laus}
|Q_{n_{l'}m'_p}| =
|Q_{n_{l'}m'_p(n_{l'}m'_p)^*}|,
\quad
|Q_{n_{l'}}| = |Q_{n_{l'}n_{l'}^*}|.
\eeq
Since $n_{l'}m'_p \R n_{l'}$ by (\ref{vli5}), we get
$n_{l'}m'_p(n_{l'}m'_p)^* = n_{l'}n_{l'}^*$, hence 
$h_{\J}(n_{l'})
= k'$ and (\ref{laus}) yields $|G_{n_{l'}}|(1+|Q_{n_{l'}}|) >
1$ and thus $n_{l'} \in U_2(k')$. By Lemma \ref{norfB}, this contradicts
$|\son(v_{q-1})| = 1$ as well.

The cases $n_{l'}m'_p \in U_3(k')$ and
$n_{l'-1}m'_p \in U_4(k')$ are analogous and can be omitted.
Therefore we may conclude that $|\son(v_{q-1}\sigma)| = 1$
and so $\downarrow\xi' = \downarrow$. 

We assume now
that
$|\son(v_{q-1})| \neq 1$. Since $q-1 < \delta$, it follows that
$|\son(v_{q-1})| > 1$.
Clearly, if $l' = d$, then $h_{\J}(n_{l'}) = h_{\J}(m_d) = k$ and so,
since $q-1 \geq 2k+2$, 
$v_{q-1}$ has
a unique son by Lemmas \ref{norfA} and \ref{norfB}, a contradiction. 
Therefore $l' > d$. 
Now (\ref{vli5}) yields $n_{l'}m'_p \R n_{l'}$, which implies
$|(\son(v_{q-1}))\sigma| > 1$ by Lemmas
\ref{ver1}--\ref{ver4}(iii). Thus $\xi' \in S(X_i)$ by definition of 
$\p_{\sigma}$. If $q$ is odd, we obtain $\xi' = \id$ by Lemmas
\ref{ver1}(iii) and \ref{ver1A}(iii), hence we may assume that $q =
2k'+2$ with $k < 
k'$.

Since $\xi'$ is the identity anyway for all the other cases, it
suffices to prove that 
$$(f(v_q), \ldots, f(v_1))\p_{\sigma} = (f(v_q), \ldots, f(v_1))$$
whenever $(v_q, \ldots, v_1) \in \ray(r_0,T)$, that is,
$$(v_q, \ldots, v_1)\sigma\psi = (f(v_q), \ldots, f(v_1)).$$
By the induction hypothesis, we have
$$(f(v_{q-1}), \ldots, f(v_1))\p_{\sigma} = (f(v_{q-1}), \ldots,
f(v_1)),$$
hence it is enough to show that 
\beq
\label{vli2}
f(v_q\sigma) = f(v_q).
\eeq
Since $|\son(v_{q-1})| > 1$, it follows from Lemma \ref{norfB} that
either $n_{l'} \in U_2(k') \cup U_3(k')$ or $n_{l'-1}
  \in U_4(k')$. 

We consider first the case $n_{l'}
  \in U_2(k')$. Since $h_{\J}(n_{l'}) = k'$, we may
replace in (\ref{vli3}) $n_{l'}$ by any element in its
$\H$-class. Indeed, if
$$\eta = (r <_{\L} n_{l'-1} <_{\L} \ldots <_{\L} n_d
<_{\L} m_{d-1} <_{\L} \ldots 
<_{\L} m_0)$$ with $r \H n_{l'}$, then $(\rho,\eta) \in V(M^I)$ by
Lemma \ref{charV} (case (V1)) and so $H(\rho,\eta) =
2h_{\J}(n_{l'})+1 = 2k'+1$ yields $[q-1,\rho] = [q-1,\eta]$. 

Thus we may assume by Lemma \ref{norfB} that either
$$v_q = [2k'+2,\rho] \hspace{.5cm} \mbox{ or }  \hspace{.5cm}
v_q = [2k'+2,\rho']$$
with
$$\rho' = (n_{l'+1} <_{\L} n_{l'} <_{\L} \ldots <_{\L} n_d
<_{\L} m_{d-1} <_{\L} \ldots 
<_{\L} m_0).$$
Write
$$\epsilon(\rho) = (x_{l'} <_{\J} \ldots <_{\J} x_0), \quad
\epsilon(\rho\sigma) = (x'_{s} <_{\J} \ldots <_{\J} x'_0).$$

Assume first that $v_q = [2k'+2,\rho]$. Writing $n_{l'} = zn_{l'-1}$,
it follows from $l' > d$, (\ref{grerus}) and (\ref{vli5}) that
$$x'_s = (n_{l'}m'_p)(n_{l'-1}m'_p)^* = zn_{l'-1}m'_p(n_{l'-1}m'_p)^*
= zn_{l'-1}n_{l'-1}^* = n_{l'}n_{l'-1}^* =
x_{l'},$$ hence
$f(v_q\sigma) = (g,\ast) = f(v_q)$ for the same $g \in
G_{m_{l'}}$.
Assume now that 
$v_q = [2k'+2,\rho']$. Since $\epsilon$ is sequential, we may write
$$\epsilon(\rho') = (x_{l'+1} <_{\J} x_{l'} <_{\J} \ldots <_{\J}
x_0)$$ for some $x_{l'+1} \in M$. Since $n_{l'} \R n_{l'}m'_p$ by (\ref{vli5}),
it follows from Theorem \ref{epz} that
$$\epsilon(\rho'\sigma) = (x_{l'+1} <_{\J} x'_s <_{\J} \ldots <_{\J}
x'_0).$$ 
Since $x'_s = x_{l'}$ as before, it follows that
$f(v_q\sigma) = (g,b) = f(v_q\sigma)$ for the same $g \in
G_{m_{l'}}$ and $b \in Q_{m_{l'}m_{l'}^*}$. Therefore 
(\ref{vli2}) holds in this case.

The case $n_{l'} \in U_3(k')$ being actually a simplification of the
preceding case, we may assume now that $n_{l'-1}
  \in U_4(k')$. 
Since $h_{\J}(n_{l'-1}) = k'$, we may
replace in (\ref{vli3}) $n_{l'-1}$ by any element in its
$\H$-class. Indeed, if
$$\eta = (n_{l'} <_{\L} r <_{\L} n_{l'-2} <_{\L} \ldots <_{\L} n_l
<_{\L} m_{l-1} <_{\L} \ldots 
<_{\L} m_0)$$ with $r \H n_{l'}$, then $(\rho,\eta) \in V(M^I)$ by
Lemma \ref{charV} (case (V2)) and so $H(\rho,\eta) =
2h_{\J}(n_{l'})+1 = 2k'+1$ yields $[q-1,\rho] = [q-1,\eta]$. 

Thus we may assume by Lemma \ref{norfB} that 
$v_q = [2k'+2,\rho']$
with
$$\rho' = (r <_{\L} n_{l'-1} <_{\L} \ldots <_{\L} n_d
<_{\L} m_{d-1} <_{\L} \ldots 
<_{\L} m_0).$$
Let $\rho'' = (n_{l'-1} <_{\L} \ldots <_{\L} n_d
<_{\L} m_{d-1} <_{\L} \ldots 
<_{\L} m_0)$ and
$$\epsilon(\rho'') = (x_{l'-1} <_{\J} \ldots <_{\J} x_0), \quad
\epsilon(\rho''\sigma) = (x'_{s} <_{\J} \ldots <_{\J} x'_0).$$
Since $h_{\J}(n_{l'-1}) = k' > k = h_{\J}(n_d)$, we have $l'-1 > d$.

Similarly to the preceding case, we have 
$$
n_{j}m'_p \R n_{j}\quad (j =
d, \ldots,l'-1),
$$ 
$$\rho''\sigma = (n_{l'-1}m'_p <_{\L} n_{l'-2}m'_p <_{\L} \ldots <_{\L}
n_{d}m'_p <_{\L} \ldots).
$$
and
$$\rho'\sigma = (rm'_p <_{\L} n_{l'-1}m'_p  <_{\L} n_{l'-2}m'_p <_{\L} \ldots <_{\L}
n_{d}m'_p <_{\L} \ldots). 
$$
Now we get $x'_s = x_{l'-1}$ as in the preceding case. Since
$\epsilon$ is sequential, we now repeat the argument of the preceding
case to reach (\ref{vli2}) as well. Therefore
(iii) is proved.

The final claim follows from Lemma \ref{finiX}. 
\qed

We can show that, by computing the length function naturally
associated by Proposition \ref{turq} to the wreath product in Theorem
\ref{pcm}, we recover the original length function $H_Y$. We need a
further lemma.

\bl
\label{ferias2}
For all $\sigma,\tau,\rho  \in \E(M^I)$, $H(\rho\sigma,\rho\tau) \geq
H(\sigma,\tau)$.
\el

\proof
Let $\sigma,\tau,\rho  \in \E(M^I)$ and assume that $\sigma \neq
\tau$. By Lemma \ref{kkk}(ii), we have $(\rho\sigma \wedge_{\L} 
\rho\tau) \leq_{\L} (\sigma \wedge_{\L} 
\tau)$. If $(\rho\sigma \wedge_{\L} 
\rho\tau) <_{\L} (\sigma \wedge_{\L} 
\tau)$, then
$$H(\rho\sigma,\rho\tau) \geq H'(\rho\sigma,\rho\tau) >
H'(\sigma,\tau)$$ yields $H(\rho\sigma,\rho\tau) \geq
H(\sigma,\tau)$. Hence we may assume that 
\beq
\label{ferias3}
(\rho\sigma \wedge_{\L} 
\rho\tau) \L (\sigma \wedge_{\L} 
\tau)
\eeq
It suffices to show that
\beq
\label{ferias4}
(\sigma,\tau) \in V(M^I) \Rw (\rho\sigma,\rho\tau) \in V(M^I).
\eeq
Indeed, let $\mu \in \E(M^I)$ and assume that $(\mu\rho\sigma \wedge_{\L} 
\mu\rho\tau) \L (\rho\sigma \wedge_{\L} 
\rho\tau)$. Then $(\mu\rho\sigma \wedge_{\L} 
\mu\rho\tau) \L (\sigma \wedge_{\L} 
\tau)$ by (\ref{ferias3}).
Since $(\sigma,\tau) \in V(M^I)$, it follows that
$(\mu\rho\sigma \wedge_{\L}  
\mu\rho\tau) \R (\mu\rho\tau \wedge_{\L} 
\mu\rho\sigma)$ and so $(\rho\sigma,\rho\tau) \in V(M^I)$. Thus
(\ref{ferias4}) holds and so does the lemma.
\qed

\bc
\label{recover}
Let $D: \Pi_{i=1}^{\delta}
(X_i,M_i) \times \Pi_{i=1}^{\delta}
(X_i,M_i) \to \oo{\N}$ be the length function defined by
$$D(\mu,\nu) = \sup\{ j \mid \mu|_{X_j \times \ldots \times X_1}
= \nu|_{X_j \times \ldots \times X_1}\}.$$
Then $D(\p_{\sigma},\p_{\tau}) =
  H_Y(\sigma,\tau)$ for all
$\sigma,\tau \in \E_Y(M^I)$.
\ec

\proof
Write $\wh{X_j} = X_j \times \ldots \times X_1$. Note that 
\beq
\label{recover2}
\wh{X_j} = (\wh{X_j} \cap \im\psi) \cup
(\wh{X_j} \setminus \im\psi).
\eeq 
We show by induction on $j$ that
\beq
\label{recover1}
\p_{\sigma}|_{\wh{X_j}} = \p_{\tau}|_{\wh{X_j}} \iff
\p_{\sigma}|_{\wh{X_j}\cap \im\psi} = \p_{\tau}|_{\wh{X_j}\cap \im\psi} 
\eeq
holds for all $\sigma,\tau
\in \E_Y(M^I)$ and $j \in \N$. 
The case $j = 0$ being trivial, assume that
$j > 0$ and (\ref{recover1}) holds for $j-1$. Let $x = (x_j,\ldots,x_1)
\in \wh{X_j} \setminus \im\psi$ and assume that
\beq
\label{recover3}
\p_{\sigma}|_{\wh{X_j}\cap \im\psi} = \p_{\tau}|_{\wh{X_j}\cap
  \im\psi}.
\eeq
We must show that either $x \notin \dom\p_{\sigma} \cup
\dom\p_{\tau}$ or else $x\p_{\sigma} =
x\p_{\tau}$. 

Suppose first that $x \in \dom\p_{\sigma} \setminus \im\psi$. Then
$(x_{j-1},\ldots,x_1) = (f(v_{j-1}),\ldots,f(v_1))$ for some
$(v_{j-1},\ldots,v_1,r_0) \in \ray(r_0,T)$ such that
$|(\son(v_{j-1}))\sigma| > 1$. Then $(f(v_{j-1}),\ldots,f(v_1)) \in
\im\psi$ and since $\p_{\sigma}, \p_{\tau}$ are sequential,
(\ref{recover3})
yields $\p_{\sigma}|_{\wh{X_{j-1}}\cap \im\psi} = \p_{\tau}|_{\wh{X_{j-1}}\cap
  \im\psi}$ and \linebreak
$(\son(v_{j-1}))\sigma = (\son(v_{j-1}))\tau$. Hence
$|(\son(v_{j-1}))\tau| = |(\son(v_{j-1}))\sigma| > 1$ and 
$$x\p_{\tau} = (x_j,(x_{j-1},\ldots,x_1)\p_{\tau}) =
(x_j,(x_{j-1},\ldots,x_1)\p_{\sigma}) = x\p_{\sigma}.$$
The case $x \in \im\psi$ follows directly from (\ref{recover3}). By
symmetry, we get $\p_{\sigma}|_{\wh{X_j}} =
\p_{\tau}|_{\wh{X_j}}$. Thus (\ref{recover1}) holds.

Now it suffices to show that
\beq
\label{recover4}
\p_{\sigma}|_{\wh{X_j}} = \p_{\tau}|_{\wh{X_j}}
\iff  H_Y(\sigma,\tau) \geq j.
\eeq
Indeed, $\p_{\sigma}|_{\wh{X_j}\cap \im\psi} = \p_{\tau}|_{\wh{X_j}\cap
  \im\psi}$ if and only if $(v_j,\ldots,v_1,r_0)\theta_{\sigma}\psi =
(v_j,\ldots,v_1,r_0)\theta_{\tau}\psi$ for every
$(v_{j},\ldots,v_1,r_0) \in \ray(r_0,T)$. Since $\psi$ is one-to-one,
this is equivalent to
\beq
\label{recover5}
\forall (v_{j},\ldots,v_1,r_0) \in \ray(r_0,T)\hspace{.7cm} 
(v_j,\ldots,v_1,r_0)\theta_{\sigma} = 
(v_j,\ldots,v_1,r_0)\theta_{\tau}.
\eeq
The vertices of $T$ with depth $j$ are precisely those of the form 
$[j,\rho]$ with $\rho \in \E_Y(M^I)$. Since $\theta_{\sigma}$ and
$\theta_{\tau}$ are sequential, (\ref{recover5}) is
equivalent to
$$\forall \rho \in \E_Y(M^I) \; [j,\rho]\sigma = [j,\rho]\tau$$
and so to
$$\forall \rho \in \E_Y(M^I) \; H_Y(\rho\sigma,\rho\tau) \geq j.$$
By Lemma \ref{ferias2}, the latter is equivalent to $H_Y(\sigma,\tau)
\geq j$ and so (\ref{recover4}) holds as required.
\qed

We present now some further corollaries of Theorem \ref{pcm}.

\bc
\label{moura}
Let $M$ be a $Y$-semigroup and let $\delta = \sup \{
h_{\J}(u) \mid u \in \Phi_{3,Y}(M) \} \in \oo{\N}$. Then there exists
an embedding 
$\p$ of $\E((\Phi_{3,Y}(M))^I)$ into an iterated
wreath product of full transformation semigroups
$\Pi_{i=1}^{\delta}
(X_i,M_i) = \ldots \circ (X_2,M_2) \circ (X_1,M_1)$
such that:
\bi
\item[(i)] $M_{2k+1}$ is a submonoid of $\{ \id_{X_{2k+1}}\} \cup
  K(X_{2k+1})$ for $2k+1 \leq \delta$.
\item[(ii)] $M_{2k+2}$ is a submonoid of $S(X_{2k+2}) \cup
  K(X_{2k+2})$  for $2k+2 \leq \delta$; 
if $\{ R_{\lambda} \mid \lambda \in \Lambda \}$ is the set of all
  $\R$-classes of $\Phi_{3,Y}(M)$ contained in $U_2(k) \cup U_3(k)
  \cup U_4(k)$, then 
$$M_{2k+2} \cap S(X_{2k+2}) \cong \oplus_{\lambda \in \Lambda} G'_{\lambda},
$$
where $G'_{\lambda}$ is a subgroup of $G_{\lambda}$.
\item[(iii)] $\p$ has the Zeiger property.
\ei
Furthermore, if $Y$ is finite, then the
$X_i$ (and consequently the 
$M_i$) are all finite, and the canonical morphism $\eta:
\E((\Phi_{3,Y}(M))^I) \to M$ is 
aperiodic. 
\ec

\proof
The
existence of $\p$ and its properties follow from Proposition 
\ref{fjap}(i) and Theorem \ref{pcm}. The aperiodicity of $\eta$
follows from Propositions \ref{rex}(i) and
\ref{fjap}(ii) since the composition of aperiodic morphisms is clearly
aperiodic.
\qed

Let $G = \langle A \rangle$ be an infinite group generated by $A = A
\cup A\inv$. The
{\em Cayley graph} $\Gamma(G,A)$ is the directed labeled graph defined
by
\bi
\item[] $V(\Gamma(G,A)) = G$;
\item[] $E(\Gamma(G,A)) = \{ (g,a,h) \in G \times A \times G \mid ga =
  h\}$.
\ei
The {\em Munn-Margolis-Meakin expansion} $M_3(G,A)$ (see
\cite{12MM,10Rho,RS} is defined by
$$M_3(G,A) = \{ (\gamma,g); \; \gamma\mbox{ is a finite connected
  subgraph of }\Gamma(G,A)\mbox{ and }1,g \in \gamma\}.$$
With the binary operation
$$(\gamma,g)(\gamma',g') = (\gamma \cup g\gamma',gg'),$$
$M_3(G,A)$ is a E-unitary inverse $A$-monoid
\cite{12MM}. Moreover, the morphism
$$\begin{array}{rcl}
\alpha: M_3(G,A)&\to&G\\
(\gamma,g)&\mapsto&g
\end{array}$$
provides the maximal group homomorphic image of $M_3(G,A)$.

Since a finite graph can have only finitely many
subgraphs, it is easy to see that $M_3(G,A)$ is finite $\J$-above as
well.

We recall that a semigroup $M$ is {\em orthodox} if it is regular and
the subset $E(M)$ of all idempotents of $M$ constitutes a subsemigroup
of $M$.
A monoid $M$ is said to be an {\em orthodox covering} of a group $G$
if $M$ is orthodox and there exists an onto homomorphism $\p:M \to G$
such that $1\pinv = E(M)$.

\bc
\label{ortho}
Let Let $G = \langle A \rangle$ be an infinite group. Then
$\E_A(M_3(G,A))$ is an orthodox covering of $G$ and there
exists an embedding 
$\p$ of $\E_A(M_3(G,A))$ into an iterated
wreath product of full transformation semigroups
$\Pi_{i=1}^{\infty}
(X_i,M_i) = \ldots \circ (X_2,M_2) \circ (X_1,M_1)$
such that:
\bi
\item[(i)] $M_i $ is a finite submonoid of $\{ 1_{X_i}\} \cup K(X_i)$ for
  $i$ odd.
\item[(ii)] $M_i$ is a finite submonoid of $S(X_i) \cup K(X_i)$ for $i$ even;
the local groups are then finite subgroups of $G$.
\item[(iii)] $\p$ has the Zeiger property.
\ei
\ec

\proof
Note that $(M_3(G,A))\setminus \{ (\{1\},1) \}$ is an $A$-semigroup and
$(M_3(G,A)) \cong ((M_3(G,A))\setminus \{ (\{1\},1) \})^I$. Since $G$
is infinite, it follows easily that $M_3(G,A)$ has arbitrarily long
$\J$-chains and so
$\sup \{
h_{\J}(u) \mid u \in M_3(G,A) \} = \omega$. Since $M_3(G,A)$ is finite
$\J$-above, the existence of $\p$
and its properties follow from Theorem \ref{pcm} and its proof, since any local
group must be the Sch\"utzenberger group of some $\J$-class and
therefore a (group) $\H$-class since $M_3(G,A)$ is inverse. It follows
that such a group must be a finite subgroup of $G$ (see \cite{12MM}
for more details).

By Proposition \ref{rex}(ii), $\E_A(M_3(G,A))$ is regular. We consider
the canonical morphisms $\eta: \E_A(M_3(G,A)) \to M_3(G,A)$ and $\alpha:
M_3(G,A) \to G$. Clearly, $1\alpha\inv = E(M_3(G,A))$. By  Proposition
\ref{rex}(iii), 
$$1(\eta\alpha)\inv = 1\alpha\inv\eta\inv = (E(M_3(G,A)))\eta\inv =
E(\E_A(M_3(G,A)))$$ and so 
$\E_A(M_3(G,A))$ is an orthodox covering of $G$.
\qed

\section{Free Burnside monoids}

Given $p,q \geq 1$, let $\B(p,q)$ denote the variety of semigroups
defined by the identity $x^{p+q} = x^p$. Given a set $X$, we denote by
$B_X(p,q)$ the free $\B(p,q)$-semigroup on $X$. Clearly, $B_X(p,q)$
can be defined by the semigroup presentation 
\beq
\label{pres1}
\langle X \mid u^{p+q} = u^p \; (u \in X^+)\rangle.
\eeq
We say that $B_X(p,q)$ is a {\em free Burnside semigroup}. 
The corresponding {\em free Burnside monoid} $B^I_X(p,q)$ can be obtained by
adjoining an identity to $B_X(p,q)$.
For details on $B_X(p,q)$, the reader is referred to \cite{13McC, GST,
NF},
\cite{dLV} and \cite{Gub}.

\bl
\label{sameexp}
For all $p,q \geq 1$, $B^I_X(p,q) \cong \E_X(B^I_X(p,q))$.
\el

\proof
Take the canonical surjective morphism $\eta: \E_X(B_X^I(p,q)) \to
B_X^I(p,q)$. Since $B_X(p,q)$ is presented by (\ref{pres1}), it suffices
to show that $\sigma^{p+q} = \sigma^p$ for every $\sigma \in
\E(B_X^I(p,q))$. Let
$$\sigma = (m_l <_{\L} \ldots <_{\L} m_0 = I)$$
and write $\sigma^p = (m_l^p <_{\L} n_t <_{\L} \ldots <_{\L}
n_0)$. Then
$$\sigma^{p+q} = \lm(m_l^{p+q} \leq_{\L} \ldots  \leq_{\L} m_l^p  <_{\L} n_t
<_{\L} \ldots <_{\L} n_0).$$ Since $m_l^{p+q} = m_l^{p}$, it follows
that $\sigma^{p+q} = \sigma^p$ as required.
\qed

\bp
\label{abeli}
\cite{13McC}
For all $p,q \geq 1$, $B_X(p,q)$ is finite $\J$-above and its maximal
subgroups are cyclic.
\ep

Clearly, if $|X| \leq 1$ then $B_X(p,q)$ is finite. From now on, we
assume that $|X| \; > 1$. Then $B_X(p,q)$ has infinite $\J$-chains
\cite{13McC,NF}. Now Theorem \ref{pcm} yields

\bt
\label{appcm}
Let $p,q \geq 1$ and $X$ be a finite set with $|X| \; > 1$.
Then there exists an embedding
$\p$ of $B_X^I(p,q)$ into an iterated
wreath product of finite partial transformation semigroups
$\Pi_{i=1}^{\infty}
(X_i,M_i) = \ldots \circ (X_2,M_2) \circ (X_1,M_1)$
such that:
\bi
\item[(i)] $M_{2k+1}$ is a finite submonoid of $\{ \id_{X_{2k+1}}\} \cup
  K(X_{2k+1})$ for every $k$.
\item[(ii)] $M_{2k+2}$ is a finite submonoid of $S(X_{2k+2}) \cup
  K(X_{2k+2})$  for every $k$; 
if $\{ R_{\lambda} \mid \lambda \in \Lambda \}$ is the set of all
  $\R$-classes of $M$ contained in $U_2(k) \cup U_3(k) \cup U_4(k)$, then
\beq
\label{apreg1}
M_{2k+2} \cap S(X_{2k+2}) \cong \oplus_{\lambda \in \Lambda} G'_{\lambda},
\eeq 
where $G'_{\lambda}$ is a subgroup of $G_{\lambda}$. Therefore $M_{2k+2}
\cap S(X_{2k+2})$ is a finite Abelian group.
\item[(iii)] $\p$ has the Zeiger property.
\ei
\et

A future  paper will apply the
results of this paper to elliptic actions of the free Burnside semigroups.

\section*{Acknowledgment}

The first author is 
grateful to J.  Elston by pointing some corrections to \cite{Rho1}.

The second author acknowledges support from C.M.U.P., financed by F.C.T.
(Portugal) through the programmes POCTI and POSI, with national and
E.U. structural funds.


\begin{thebibliography}{99}
\bibitem{2AB} R.~Alperin and H.~Bass, Length functions of group
  actions on $\Lambda$-trees, in: {\em Combinatorial Group Theory and
    Topology}, S.~M.~Gersten and J.~R.~Stallings (eds.), Ann. of
  Math. Stud. 111, Princeton University Press, 1987, 265--378.
\bibitem{1Arb} M.~Arbib (ed.), {\em Algebraic Theory of Machines, Languages and
  Semigroups}, Academic Press, New York, 1968. 
\bibitem{17Bir} J.-C.~Birget, The synthesis theorem for finite regular
  semigroups and its generalization, {\em J. Pure Appl. Algebra} 55
  (1988), 1--80.
\bibitem{4BR} J.-C.~Birget and J.~Rhodes, Almost finite expansions of
  arbitrary semigroups, {\em J. Pure Appl. Algebra} 32
  (1984), 239--287.
\bibitem{5Chi} I.~M.~Chiswell, Abstract length functions in groups,
  {\em Math. Proc. Cambridge Phil. Soc.} 80 (1976), 451--463.
\bibitem{CP} A.~H.~Clifford and G.~B.~Preston, {\em The Algebraic
    Theory of Semigroups} vol. I, Mathematical Surveys Nr. 7, American
  Mathematical Society, 1961.
\bibitem{dLV} A.~de~Luca and S.~Varricchio, On noncounting regular
  classes, {\em Theoret. Comp. Sci.} 100  (1992), 67--104.
\bibitem{6Eil} S.~Eilenberg, {\em Automata, Languages and Machines},
  Vol. B, Academic Press, New York, 1976.
\bibitem{7Gro} M.~Gromov, Hyperbolic groups, in: {\em Essays in group
    theory}, S.~M.~Gersten (ed.), MSRI 8, Springer-Verlag, Berlin,
  1987, 75--263.
\bibitem{Gub} V.~S.~Guba, The word problem fopr the relatively free
  semigroup satisfying $T^m = T^{m+n}$ with $m \geq 4$ or $m = 3$, $n
  = 1$, {\em Internat. J. Algebra Comput.} 3  (1993), 125--140.
\bibitem{8LHR} S.~Lazarus, K.~Henckell and J.~Rhodes, Prime
  decomposition theorem for arbitrary semigroups: general holonomy
  decomposition and synthesis theorem, {\em J. Pure Appl. Algebra} 55
  (1988), 127--172.
\bibitem{9Lyn} R.~C.~Lyndon, Length functions in groups, {\em
    Math. Scand.} 12 (1963), 209--234.
\bibitem{12MM} S.~W.~Margolis and J.~Meakin, E-unitary inverse monoids
  and the Cayley graph of a group presentation, {\em J. Pure Appl. Algebra} 58
  (1989), 45--76.
\bibitem{11MM} S.~Margolis and J.~Meakin, Graph immersions
  and inverse monoids, in: {\em Monoids and Semigroups with
    Applications}, J.~Rhodes (ed.), World Scientific, 1991, 144--158.
\bibitem{13McC} J.~McCammond, The solution to the word problem for the
  relatively free semigroup satisfying $T^a = T^{a+b}$ with $a \geq
  b$, {\em Internat. J. Algebra Comput.}
  1.1 (1991), 1--32.
\bibitem{NF} J.~McCammond, Normal forms for free aperiodic semigroups,   
{\em Internat. J. Algebra Comput.} 11.5 (2001), 581--625.
\bibitem{GST} J.~McCammond, J.~Rhodes and B.~Steinberg, Geometric
  semigroup theory, soon to appear in arXiv.
\bibitem{14Rei} N.~Reilly, The Rhodes expansion and free objects in
  varieties of completely regular semigroups, {\em J. Pure Appl. Algebra} 69
  (1990), 89--109.
\bibitem{15Rho} J.~Rhodes, Infinite iteration of matrix semigroups,
  Part II, {\em J. Algebra} 100 (1986), 25--137.
\bibitem{16Rho} J.~Rhodes, A short proof that $\widehat{S}_A^+$ is
  finite if $S$ is finite, {\em J. Pure Appl. Algebra} 55
  (1988), 197--198.
\bibitem{10Rho} J.~Rhodes, Survey of global semigroup theory, in: {\em
    Lattices, Semigroups and Universal Algebra}, J.~Almeida,
  G.~Bordalo and P.~Dwinger (eds.), Plenum Press, 1990, 243--270.
\bibitem{Rho1} J.~Rhodes, Monoids acting on trees: elliptic and wreath
  products and the Holonomy Theorem for arbitrary monoids with
  applications to infinite groups, {\em Internat. J. Algebra Comput.}
  1.2 (1991), 253--279; with Erratum to diagram p. 274.
\bibitem{RS} J. Rhodes and B. Steinberg, {\em The q-theory of finite
semigroups}, Springer Monographs in Mathematics, 2009.
\bibitem{18Ser} J.-P.~Serre, {\em Trees}, Springer-Verlag, New York, 1980.
\bibitem{19Sta} J.~R.~Stallings, Topology of finite graphs, {\em
    Invent. Math.} 71 (1983), 551--565.
\bibitem{20Til} B.~Tilson, Chapters XI and XII in \cite{6Eil}.
\end{thebibliography}
\end{document}